\documentclass[11pt,english,a4paper]{smfart}
\usepackage[english]{babel}
\usepackage{xcolor}
\usepackage{pbsi}
\usepackage[T1]{fontenc}
\usepackage{stmaryrd}
\usepackage{mathabx}
\usepackage{dsfont}
\usepackage[mathscr]{euscript}
\usepackage{cancel}
\usepackage{array}
\usepackage{float}
\usepackage{placeins}
\usepackage[pagebackref,colorlinks=true,linkcolor=blue,citecolor=blue,urlcolor=red, hypertexnames=true]{hyperref}

 \usepackage[abbrev,backrefs,nobysame]{amsrefs}
 \usepackage{caption}
 \usepackage{enumerate}

\usepackage{amsmath,amssymb,bm,amsfonts}
\usepackage{nicematrix}
\usepackage[all]{xy}
\entrymodifiers={+!!<0pt,\fontdimen22\textfont2>}

\usepackage{mathrsfs}

\usepackage{graphicx}

\usepackage{color}

\newtheorem{theorem}{Theorem}[section]
\newtheorem{lemma}[theorem]{Lemma}
\newtheorem{observation}[theorem]{Observation}
\newtheorem{corollary}[theorem]{Corollary}
\newtheorem{proposition}[theorem]{Proposition}

\newcounter{hypothesis}
\renewcommand{\thehypothesis}{(H$_{\arabic{hypothesis}}$)}

\theoremstyle{remark}
\newtheorem{case}{Case}

\newtheorem*{remark*}{\it Remark}

\newcommand{\flba}[2]{
\xymatrix@C15pt{#1\ar@{|->}[r]&#2}}
\newcommand{\flcourte}[2]{
\xymatrix@C12pt{#1\ar[r]&#2}}

{\theoremstyle{definition}}

{\theoremstyle{definition}\newtheorem{example}[theorem]{Example}}

\theoremstyle{definition}

\newtheorem{question}[theorem]{Question}
\newtheorem{fact}[theorem]{Fact}

\newtheorem{problem}[theorem]{Open Problem}

{\theoremstyle{definition}\newtheorem{remark}[theorem]{Remark}}

\def\D{\ensuremath{\mathbb D}}

\def\T{\ensuremath{\mathbb T}}
\def\R{\ensuremath{\mathbb R}}
\def\Z{\ensuremath{\mathbb Z}}

\def\C{\ensuremath{\mathbb C}}
\def\Q{\ensuremath{\mathbb Q}}
\def\N{\ensuremath{\mathbb N}}

\def\bth{\begin{theorem}}
\def\blm{\begin{lemma}}
\def\bpr{\begin{proposition}}
\def\bpf{\begin{proof}}
\def\epf{\end{proof}}
\def\epr{\end{proposition}}
\def\elm{\end{lemma}}
\def\eth{\end{theorem}}
\def\bco{\begin{corollary}}
\def\eco{\end{corollary}}
\def\be{\begin{enumerate}}
\def\ee{\end{enumerate}}
\def\bea{\begin{enumerate}[\rm (a)]}
\def\beun{\begin{enumerate}[\rm (1)]}
\def\bei{\begin{enumerate}[\rm (i)]}

\renewcommand{\Im}{\operatorname{Im}}
\renewcommand{\Re}{\operatorname{Re}}

\newcommand{\AAA}{\mathcal{A}}

\newcommand{\vertiii}[1]{{\left\vert\kern-0.25ex\left\vert\kern-0.25ex\left\vert #1 
    \right\vert\kern-0.25ex\right\vert\kern-0.25ex\right\vert}}

\newcommand{\gd}{G_{\delta }}

\newcommand{\bx}{{\mathcal B}(X)}

\newcommand{\lan}{\langle}
\newcommand{\ran}{\rangle}

\newcommand{\omg}{\omega}
\newcommand{\Omg}{\Omega}
\newcommand{\ta}{\tau}
\newcommand{\indic}{\mathds{1}}
\newcommand{\vph}{\varphi}

\numberwithin{equation}{section}

\author[Valentin Gillet]{Valentin Gillet}
\address[Valentin Gillet]{Univ. Lille, CNRS, UMR 8524 - Laboratoire Paul Painlevé, F-59000 Lille, France}
\email{valentin.gillet@univ-lille.fr}

\thanks{This work was supported in part by the project COMOP of the French National Research Agency (grant ANR-24-CE40-0892-01) and by the Labex CEMPI (ANR-11-LABX-0007-01). The author also acknowledges the support of the CDP C2EMPI, as well as the French State under the France-2030 programme, the University of Lille, the Initiative of Excellence of the University of Lille, and the European Metropolis of Lille for their funding and support of the R-CDP-24-004-C2EMPI project.}

\begin{document}

\title[Linear dynamics of random products of operators]{Linear dynamics of random products of operators}


\begin{abstract} 
We study the linear dynamics of the random sequence $(T_n(.))_{n \geq 1}$ of the operators $T_n(\omg) := T(\tau^{n-1}\omg) \dotsm T(\tau \omg) T(\omg)$, $n \geq 1$. These products depend on an ergodic measure-preserving transformation $\tau : \T \to \T$ on the probability space $(\T,m)$ and on a strongly measurable map $T : \T \to \bx$, where $X$ is a separable Fréchet space. We will be focusing on the case where $T(\omg)$ is equal to an operator $T_1$ on $X$ for every $\omg \in A_1$ and equal to an operator $T_2$ on $X$ for every $\omg \in A_2$, where $A_1, A_2$ are two disjoint Borel subsets of $[0,1)$ such that $A_1 \cup A_2 = [0,1)$ and $m(A_k) > 0$ for $k = 1,2$. More precisely, we will be focusing on the case where the operators $T_1$ and $T_2$ are adjoints of multiplication operators on the Hardy space $H^2(\D)$, as well as the case where $T_1$ and $T_2$ are entire functions of exponential type of the derivation operator on the space of entire functions. Finally, we will study the linear dynamics of a case of a random product $T_n(\omg)$ for which the operators $T(\tau^i \omg), i \geq 0$, do not commute. We will give particular importance to the case where the ergodic transformation is an irrational rotation or the doubling map on $\T$.
\end{abstract}

\maketitle

\par\bigskip
\section{Introduction}\label{Section introduction}

\smallskip

\subsection{Objectives of the article and motivation}
The aim of this article is to study the linear dynamics of random sequences of products of operators on a separable Fréchet space. That is, we would like to determine if the sequence $(T_n(.))_{n \geq 1}$ of the random operators
\begin{align} \label{formeproduitopintro}
    T_n(\omega) = T(\tau^{n-1} \omega) \dotsm T(\tau \omega) T(\omg)
\end{align}
is universal, or topologically weakly mixing, or topologically mixing, or none of it, in the sense that the sequence $(T_n(\omega))_{n \geq 1}$ has this property for almost every $\omega \in \T := \R / \Z$, or none of these properties. The probability space for our study will be the set $\T$, equipped with its Borel subsets and with the normalized Lebesgue measure $m$. The map $\tau : \T \to \T$ will be a measure-preserving transformation on $\T$, which means that 
$$
m(\tau^{-1}(B)) = m(B)
$$
for every Borel subset $B$ of $\T$, and the map $T : \T \to \bx$ will be given by
\begin{equation}  \notag
    T(\omega) = \left\{
    \begin{array}{ll}
        T_1 & \mbox{if} \quad \omega \in A_1 \\
        T_2 & \mbox{if} \quad \omega \in A_2
    \end{array},
\right.
\end{equation}
with $T_1, T_2$ two operators on a separable Fréchet space $X$, and with $A_1, A_2$ two disjoint Borel subsets of $[0,1)$ such that $A_1 \cup A_2 = [0,1)$ and $m(A_k) > 0$ for $k = 1,2$. Such types of products (\ref{formeproduitopintro}) appear when we consider linear cocycles $F : \T \times X \to \T \times X, \, (\omega, x) \mapsto (\tau \omg, T(\omega) x)$. In this case, we have $F^n(\omega, x) = (\tau^n x, T_n(\omega)x)$. Some work on the dynamical properties of linear cocycles, which is somewhat close to our study, are considered in \cite{BCH10}, \cite{CH08} and \cite{MMPR}.

\smallskip

The motivation of our work comes from the study of the asymptotic behavior of random products of matrices, which come into play, for example, in considering solutions to differential or difference equations with random coefficients. 
Bellman was the first to study random products of $2 \times 2$ matrices with strictly positive entries (\cite{Bell}). Roughly speaking, he proved that the weak law of large numbers holds for these kinds of random products.
A few years later, Furstenberg and Kesten also considered random products of matrices (\cite{FK}). In particular, they extended Bellman's result by proving that the strong law of large numbers is applicable for random products of $2\times 2$ matrices with strictly positive entries, under some assumptions on the coefficients of these matrices.

\smallskip

While random products of matrices and operators have been extensively studied from a probabilistic and statistical viewpoint - most notably via Lyapunov exponents (see, for instance, \cite{VIANA}) - the linear dynamics of such random products remains little explored. This work aims to explore the dynamical properties (e.g., orbit structure, universality) of random products of operators of the form (\ref{formeproduitopintro}). 

\smallskip

In particular, our study includes the case of random products of independent variables on $\T$ with a common distribution. Indeed, let $A_1 = [0,1/2)$ and $A_2 = [1/2, 1)$. Let $\tau $ be the doubling map on $\T$ defined by $\tau(x) = 2x$. Let $T : \T \to \bx$ be defined by
\begin{equation}  \notag
    T(\omega) = \left\{
    \begin{array}{ll}
        T_1 & \mbox{if} \quad \omega \in A_1 \\
        T_2 & \mbox{if} \quad \omega \in A_2
    \end{array},
\right.
\end{equation}
with $T_1, T_2$ two bounded operators on a separable Fréchet space $X$.
For each $n \geq 1$ and $\omg \in \T$, let $T_n(\omg) := T(\tau^{n-1} \omg) \dotsm T(\tau \omg) T(\omg)$. Let us denote by $\varphi : \{ 0,1 \}^{\N_0} \to \T$ the map defined by $\varphi((x_n)_{n \geq 0}) = \displaystyle \sum_{n \geq 0} \frac{x_n}{2^{n+1}}$ in $\T$, and by $\sigma : \{ 0 , 1 \}^{\N_0} \to \{ 0,1 \}^{\N_0}$ the backward shift defined by $\sigma((x_n)_{n \geq 0}) = (x_{n+1})_{\geq 0}$. The map $\varphi$ is a bijection from the set $\Omega_0$ of sequences in $\{ 0 , 1\}^{\N_0}$ not eventually equal to $1$ onto $\T$, and $\tau \varphi = \varphi \sigma$.
For every $j \geq 0$, let $\hat{X}_j : u = (u_k)_{k \geq 0} \in \{ 0 , 1 \}^{\N_0} \mapsto u_j  $ be the $j$-th canonical projection. The random variables $\hat{X}_j, j \geq 0,$ are independent on the probability space $(\{ 0,1 \}^{\N_0}, \mathcal{C}, \nu^{\otimes \N_0})$ and have the same distribution $\nu := \frac{1}{2}(\delta_0 + \delta_1) $, where $\mathcal{C}$ denotes the $\sigma$-algebra generated by the cylinders and $\nu^{\otimes \N_0}$ denotes the product measure of $\nu$.
For every $\omg = \varphi(u) \in \T$ with $u \in \Omg_0$, and for every $j \geq 0$, we have $T(\tau^j \omg) = T(\varphi(\sigma^j u))$, which is equal to $T_1$ if and only if $\varphi(\sigma^ju) \in A_1$, that is, if and only if $\hat{X_j}(u) = 0$, and is equal to $T_2$ if and only if $\hat{X_j}(u) = 1$. In particular, the random variables $T(\tau^j \omg), j \geq 0,$ are independent on $\T$ and have the same distribution, namely $\frac{1}{2}(\delta_{T_1} + \delta_{T_2})$.

\smallskip

\subsection{Linear dynamics} \label{subsection lin dynamics}
We now introduce the notions of linear dynamics that we will use for our study. For background on linear dynamics, we refer to Bayart's and Matheron's book \cite{BM}, and to Grosse-Erdmann's and Peris Manguillot's book \cite{GEPM}. A good introduction to the notion of universality is given in \cite{GE}.

\smallskip

A sequence of continuous maps $(T_n)_{n \geq 1}$ from a metric space $X$ to itself is said to be universal if there is an element $x \in X$ such that its orbit under $(T_n)_{n \geq 1}$,
$$
\textrm{orb}(x,(T_n)) := \{ T_n x : n \geq 1 \},
$$
is dense in $X$. Fekete was the first to give an example of a universal family. He showed that there exists a real power series $\displaystyle \sum_{n \geq 1} a_n \, x^n$ on $[-1,1]$ that diverges at every point $x \ne 0$ in the worst possible way, in the sense that to every continuous function $g$ on $[-1,1]$ with $g(0) = 0$, there exists an increasing sequence of positive integers such that $\displaystyle \sum_{n = 1}^{n_k} a_n \, x^n \underset{k\to \infty}{\longrightarrow} g$ uniformly on $[-1,1]$. In this case, $T_n$ is the map that associates to $f$ its partial Taylor series of order $n$ at $0$. This example of universality shows two aspects of universality: an aspect of maximal divergence, and the aspect of the existence of a single object which allows us to approximate a maximal class of objects.

\smallskip

The notion of topological transitivity also plays an important role in topological and linear dynamics.
We say that a sequence of continuous maps $(T_n)_{n \geq 1}$ from a space $X$ to itself is topologically transitive if, for any pair $U,V \subset X$ of nonempty open sets, there exists an integer $n \geq 1$ such that $T_n(U) \cap V \ne \emptyset$. If this property holds for sufficiently large $n \geq 1$, then we say that the sequence $(T_n)_{n \geq 1}$ is topologically mixing. Finally, we say that the sequence $(T_n)_{n \geq 1}$ is topologically weakly mixing if, for any nonempty open sets $U_1, U_2, V_1, V_2$ of $X$, there is some $n \geq 1$ such that $T_n(U_1) \cap V_1 \ne \emptyset$ and $T_n(U_2) \cap V_2 \ne \emptyset$. Every topologically mixing sequence is obviously topologically weakly mixing, and every topologically weakly mixing sequence is topologically transitive. If the sequence $(T_n)_{n \geq 1}$ is topologically transitive and if the maps $T_n$ are invertible, then the sequence $(T_n^{-1})_{n \geq 1}$ is also topologically transitive. 

\smallskip

For a single continuous map $T$ on $X$, we say that $T$ is topologically transitive if the sequence $(T^n)_{n \geq 1}$ is topologically transitive. When the space $X$ is a separable complete metric space without isolated points, a continuous map $T$ is topologically transitive if and only if there is a point $x \in X$ with a dense orbit under $T$. In fact, in this situation, the map $T$ has a dense $\gd$-set of points with a dense orbit under $T$. This is Birkhoff's transitivity theorem (\cite[Theorem 1.2]{BM} and \cite[Theorem 1.16]{GEPM}). For a sequence of continuous and commuting maps $T_n$ on $X$, it is possible to generalize this result in the following way (\cite[Exercise 1.6.2]{GEPM}).

\smallskip

\bpr \label{thbirkhoffsequence}
Let $T_n : X \to X$, $n \geq 1$, be commuting continuous maps on a separable complete metric space $X$. Suppose that each $T_n$, $n \geq 1$, has a dense range. The following assertions are equivalent:
\begin{enumerate}[(i)]
    \item $(T_n)_{n \geq 1}$ is topologically transitive.
    \item there is a point $x \in X$ whose orbit under $(T_n)_{n \geq 1}$ is dense in $X$.
\end{enumerate}
In this situation, the set of points in $X$ with a dense orbit under $(T_n)_{n \geq 1}$ is a dense $\gd$-set.
\epr

\smallskip

A topologically transitive continuous linear map on a separable Fréchet space $T$ is also called a hypercyclic operator. In the linear setting, we have a more useful criterion than Birkhoff's transitivity theorem to prove that a sequence of operators on a separable Fréchet space $X$ is universal.

\smallskip

\bpr[{\cite[Corollary 1.4]{GodShapiro}}] \label{critunivforsequence}
Let $X$ be a separable Fréchet space and $(T_n)_{n \geq 1}$ a sequence of operators on this space. Suppose that there are dense subsets $\mathcal{D}_1 $ and $\mathcal{D}_2$ of $X$, a strictly increasing sequence $(n_k)_{k \geq 1}$ of positive integers, and maps $S_{n_k} : \mathcal{D}_2 \to X$, $k \geq 1$, such that
\begin{enumerate}[(i)]
    \item $T_{n_k} x \underset{k\to \infty}{\longrightarrow} 0\quad \textrm{for every} \; x \in \mathcal{D}_1$.
    \item $S_{n_k} y \underset{k\to \infty}{\longrightarrow} 0\quad \textrm{for every} \; y \in \mathcal{D}_2$.
    \item $T_{n_k} S_{n_k} y \underset{k\to \infty}{\longrightarrow} y \quad \textrm{for every} \; y \in \mathcal{D}_2$.
\end{enumerate}
Then the sequence $(T_n)_{n \geq 1}$ is topologically weakly mixing, and in particular universal. If moreover this criterion is satisfied with respect to the full sequence $(n)$, then the sequence $(T_n)_{n \geq 1}$ is topologically mixing. 
\epr

\smallskip

The criterion of Proposition \ref{critunivforsequence} is called the Universality Criterion. This criterion will be a powerful tool for our study. Observe that the maps $S_n$ are not required to be linear nor continuous.

\smallskip

An important example of operators in this paper will be the adjoint of the multiplication operator $M_\phi$ on the Hardy space $H^2(\D)$. Its adjoint $(M_\phi)^*$ is hypercyclic on $H^2(\D)$ if and only if $\phi(\D) \cap \T \ne \emptyset$, when $\phi$ is a nonconstant bounded holomorphic function on $\D$ (\cite[Theorem 4.5]{GodShapiro}). We will also consider entire functions of exponential type of the derivation operator $D$ on the space of entire functions. These operators are exactly the operators on the space of entire functions that commute with $D$, and they are always topologically mixing and hence hypercyclic, provided they are not a scalar multiple of the identity (\cite[Theorem 5.1]{GodShapiro}).

\subsection{Main results}
We study the linear dynamics of random sequences $(T_n(.))_{n \geq 1}$ when the operators $T(\tau^i \omega), i \geq 0$, are adjoints of multiplication operators on the Hardy space $H^2(\D)$, or entire functions of exponential type of the derivation operator $D$ on the space of entire functions. We also study a case of random products of operators on the space of entire functions, where the operators $T(\tau^i \omega)$ do not commute. The ergodic properties of the transformation $\tau : \T \to \T$ will play an important role in examining the universal behavior of the sequence $(T_n(\omega))_{n \geq 1}$, for almost every $\omega \in \T$. We put together here some of our main results. The notations will be defined in Subsection \ref{sectionnotetdefinitions}.  

\smallskip

Thanks to the fact that the adjoint of a multiplication operator has plenty of eigenvalues, we obtain a first sufficient criterion which gives the universality of the random sequence $(T_n(.))_{n \geq 1} $ in the case of adjoints of multiplication operators. 

\smallskip

\bth \label{propunivvaleursproprescocycleintroduction}
Let $\tau$ be an ergodic measure-preserving transformation on $(\T,m)$. Let $A_1, A_2$ be two disjoint Borel subsets of $[0,1)$ such that $A_1 \cup A_2 = [0,1)$ and $m(A_k) > 0$ for $k = 1,2$.  

Let $\phi_1, \phi_2 \in H^\infty(\D)$ be nonconstant on $\D$ and suppose that there exist $\lambda,\mu \in \D$ such that
\begin{align} 
   \label{eq1critvpscocyclemesnonegalesintro} &\lvert \phi_1(\lambda) \rvert^{m(A_1)} \lvert \phi_2(\lambda)\rvert^{m(A_2)} < 1 \quad \textrm{and}\\
 \label{eq2critvpscocyclemesnonegalesintro}   &\lvert \phi_1(\mu) \rvert^{m(A_1)} \lvert \phi_2(\mu)\rvert^{m(A_2)} > 1.
\end{align}

Suppose that, for every $\omega \in \T$, the operator $T(\omega)$ is given by
\begin{equation} \notag
    T(\omega) = \left\{
    \begin{array}{ll}
        (M_{\phi_1})^* & \mbox{if} \quad \omega \in A_1 \\
        (M_{\phi_2})^* & \mbox{if} \quad \omega \in A_2
    \end{array}.
\right.
\end{equation}
Then, the random sequence $(T_n(.))_{n \geq 1}$ is topologically mixing, that is, the sequence $(T_n(\omg))_{n \geq 1}$ is topologically mixing for almost every $\omega \in \T$.
\eth

\smallskip
In the case $m(A_1) = m(A_2)$, Conditions (\ref{eq1critvpscocyclemesnonegalesintro}) and (\ref{eq2critvpscocyclemesnonegalesintro}) are equivalent to $(\phi_1 \phi_2 )(\D) \cap \T \ne \emptyset$ and $\phi_1 \phi_2$ nonconstant, by the open mapping theorem. This leads to two limit cases when $m(A_1) = m(A_2)$: the case $(\phi_1 \phi_2)(\D) \subset \D$, and the case $(\phi_1 \phi_2)(\D) \subset \C \setminus \overline{\D}$. We will see that the study of these limit cases is not that easy.

\smallskip

Under certain conditions, we establish that the random sequence of operators $(T_n(.))_{n \geq 1}$ is universal when the operators $T(\tau^i \omega)$ are entire functions of exponential type of the derivation operator $D$. This will follow from the fact that the derivation operator on $H(\C)$ has suitable eigenvectors, namely the exponential functions $e_\lambda : z \mapsto e^{\lambda z}$, with $\lambda \in \C$.

\smallskip

\bth \label{propositionfoncenttypeexpointroduction}
Let $\ta : \T \to \T$ be an ergodic measure-preserving transformation on $(\T,m)$. Let $A_1, A_2$ be two disjoint Borel subsets of $[0,1)$ such that $A_1 \cup A_2 = [0,1)$ and $m(A_k) > 0$ for $k = 1,2$. Let $\vph_1$ and $\vph_2$ be two nonconstant entire functions of exponential type. Consider the map
$$
T(\omega) := \left\{ \begin{array}{ll}
    \varphi_1(D) & \textrm{if $\omega \in A_1$}   \\
    \varphi_2(D) & \textrm{if $\omega \in A_2$}
\end{array}.
\right. 
$$
Suppose that there exist $z, w \in \C$ such that
\begin{align}
    \lvert \varphi_1(z)\rvert^{m(A_1)} \lvert \varphi_2( z)\rvert^{m(A_2)} &< 1 \label{introeqfonctypeexpo1} \\
    \textrm{and} \quad  \lvert \varphi_1(w)\rvert^{m(A_1)} \lvert \varphi_2( w)\rvert^{m(A_2)} &> 1. \label{introeqfonctypeexpo2}
\end{align}
Then, the random sequence $(T_n(.))_{n \geq 1}$ is topologically mixing.
\eth

\smallskip

In the case where $m(A_1) = m(A_2) = 1/2$, Conditions (\ref{introeqfonctypeexpo1}) and (\ref{introeqfonctypeexpo2}) are equivalent to the fact that $\vph_1 \vph_2$ is nonconstant, when $\vph_1$ and $\vph_2$ are two nonconstant entire functions of exponential type. This comes from the fact that every nonconstant entire function has a dense range in $\C$. We will also study a situation where $\vph_1 \vph_2$ is constant, showing in particular that the random sequence $(T_n(.))_{n \geq 1}$ can be non-universal even if $T(\omg)$ is hypercyclic for every $\omg \in \T$. 

Until now, all the random products of operators $T_n(.)$ we introduced before have the property that the operators $T(\tau^i \omega), i \geq 0$, commute. We now consider a case of a random sequence of operators on $H(\C)$ for which these operators do not commute:
\smallskip

\bth \label{propositioncocycleunivnoncommutintroduction}
Let $\tau : \T \to \T$ be an ergodic measure-preserving transformation on $(\T,m)$. Let $A_1, A_2$ be two disjoint Borel subsets of $[0,1)$ such that $A_1 \cup A_2 = [0,1)$ and $m(A_k) > 0$, for $k = 1,2$. Consider the map
$$
T(\omega) := \left\{ \begin{array}{ll}
    T_{\lambda, b} & \textrm{if $\omega \in A_1$}   \\
    D & \textrm{if $\omega \in A_2$}
\end{array},
\right. 
$$
where $T_{\lambda, b}$ is the operator on $H(\C)$ defined by
$$
T_{\lambda,b} f = f(\lambda \, . + b),
$$
with $\lambda \in \C$ nonzero and $b \in \C$.
Then, the random sequence $(T_n(.))_{n \geq 1}$ is topologically mixing, that is, the sequence  $(T_n(\omg))_{n \geq 1}$ is topologically mixing for almost every $\omg \in \T$.
\eth

\smallskip
The operator $T_{\lambda,b}$ is hypercyclic on $H(\C)$ precisely when $\lambda = 1$ and $b \ne 0 $ (see, for instance, \cite{BeGoMR}). An interest of Theorem \ref{propositioncocycleunivnoncommutintroduction} is that the product $T_n(\omega)$ can contain hypercyclic and non-hypercyclic operators $T(\tau^i \omega)$ that do not commute.

\smallskip
We will also consider the case of an irrational rotation and of the doubling map on $\T$. 
Given a real number $\alpha$ in $(0,1)$, the rotation of parameter $\alpha$ is defined by the transformation $R_\alpha : x \mapsto x + \alpha \pmod 1$ on $[0,1)$, where $x \pmod1$ is the fractional part of the real $x$. This map can be seen as a map on $\T$ by $\tau(x) = x + \alpha$ in $\T$, and can also be seen as a map on the unit circle by $\tau(z) = e^{2 i \pi \alpha} z$. This map is ergodic precisely when $\alpha$ is irrational. Moreover, this map is not weakly mixing. The doubling map is the map on $[0,1)$ defined by $\tau(x) = 2 x \pmod1$. This map can also be seen as a map on $\T$ by $\tau(x) = 2 x \in \T$, and also as a map on the unit circle by $\tau(z) = z^2$. Contrary to rotations, this map is mixing. These two transformations will appear many times in this article. 

\smallskip

When the ergodic transformation involved in Theorems \ref{propositionfoncenttypeexpointroduction} and \ref{propositioncocycleunivnoncommutintroduction} is an irrational rotation and when $A_1, A_2$ are intervals of $[0,1)$, the universality of the associated sequence $(T_n(\omega))_{n \geq 1}$ holds for every $\omega \in \T$. This comes from the uniform distribution mod $1$ of the sequence $(n \alpha)_{n \geq 1}$. 

Regarding the limit case $(\phi_1 \phi_2)(\D) \subset \D$ when $m(A_1) = m(A_2) = 1/2$, we studied the situation where $\phi_1 \phi_2$ is an inner function. We also generalized it to the situation where $m(A_1) \ne m(A_2)$. In this context, for the doubling map on $\T$, we obtained the following result.

\smallskip

\bth \label{cordoublmapcaspasmmesureintro}
Let $\tau $ be the doubling map on $\T$ and let $A_1$ and $A_2$ be two disjoint intervals of $[0,1)$, such that $A_1 \cup A_2 = [0,1)$ and $m(A_k) > 0$ for $k = 1,2$. Let $\phi_1,\phi_2 \in H^\infty(\D)$ be nonconstant. Suppose that one of the two following conditions holds:
\begin{enumerate}[(i)]
    \item  $\lvert \phi_1^*\rvert^{m(A_1)} \, \lvert \phi_2^*\rvert^{m(A_2)} = 1$ almost everywhere on $\T$, and one of the images $\phi_1(\D)$ or $\phi_2(\D)$ meets the unit circle $\T$.
    \item $\lvert \phi_1^*\rvert^{m(A_1)} \, \lvert \phi_2^*\rvert^{m(A_2)} = 1$ almost everywhere on $\T$, $\phi_1 \phi_2$ is not outer and either $\phi_1$ or $\phi_2$ is not inner.
\end{enumerate}
Then, the random sequence $(T_n(.))_{n \geq 1}$ is topologically weakly mixing.
\eth

\smallskip

The case of irrational rotations will be more harsh to study. This will strongly depend on the Diophantine properties of the irrational parameter $\alpha$. For example, we will establish the following result.

\smallskip

\smallskip
\bth \label{coralphairrcaspasmmesureintroduction}
Let $\alpha$ be an irrational in $(0,1)$ and let $\tau = R_\alpha$. Let $b \in (0,1)$ and let $A_1 = [0,b)$ and $A_2 = [b,1)$. Suppose that the sequence $(q_k \, b)_{k \geq 1}$ is uniformly distributed mod $1$, where $q_k$ are the denominators of the convergents of $\alpha$.
Let $\phi_1,\phi_2 \in H^\infty(\D)$ be nonconstant. Suppose that one of the two following conditions holds:
\begin{enumerate}[(i)]
    \item  $\lvert \phi_1^*\rvert^{m(A_1)} \, \lvert \phi_2^*\rvert^{m(A_2)} = 1$ almost everywhere on $\T$, and the image of either $\phi_1$ or $\phi_2$ intersects $\T$.
    \item $\lvert \phi_1^*\rvert^{m(A_1)} \, \lvert \phi_2^*\rvert^{m(A_2)} = 1$ almost everywhere on $\T$, $\phi_1 \phi_2$ is not outer and either $\phi_1$ or $\phi_2$ is not inner.
\end{enumerate}
Then, the random sequence $(T_n(.))_{n \geq 1}$ is topologically weakly mixing.
\eth

\smallskip

However, in the case $A_1 = [0,1/2)$ and $A_2 = [1/2,1)$, we do not need any assumption on the irrational $\alpha$.

\smallskip

\bth \label{corunivcasA1A2mmmesurerotirraintroduction}
Let $\alpha$ be an irrational number in $(0,1)$ and let $\tau = R_\alpha$. Let $A_1 = [0, 1/2)$ and $A_2 = [1/2, 1)$. Let $\phi_1, \phi_2 \in H^\infty(\D)$ be nonconstant, such that $\phi_1 \phi_2$ is inner and nonconstant. Suppose also that either $\phi_1$ or $\phi_2$ is not inner. Then, the random sequence $(T_n(.))_{n \geq 1}$ is topologically weakly mixing.
\eth

\smallskip

For some other intervals of $[0,1)$, it will also be possible for the random sequence $(T_n(.))_{n \geq 1}$ to be non-universal. This is, for example, the case of the following result that we obtained.

\smallskip

\bth \label{coralphabZalpharotintroduction}
Let $\alpha$ be an irrational number in $(0,1)$ and let $\tau = R_\alpha$.
Let $A_1 = [0,b)$ and $A_2 = [b,1)$, with $b \in (0,1)$. Suppose that $b \in \Z \alpha$.
Let $\phi_1, \phi_2 \in H^\infty(\D)$ be nonconstant, such that 
$$
\lvert \phi_1 ^* \rvert^{m(A_1)} \lvert \phi_2^* \rvert^{m(A_2)} = 1 \quad \textrm{almost everywhere on $\T$.}
$$
Then, for every $\omega \in \T$, the sequence $(T_n(\omega))_{n \geq 1}$ is not universal.
\eth

\medskip
These results show in particular that the random sequence $(T_n(.))_{n \geq 1}$ can be universal, even if both operators $T_1$ and $T_2$ are not hypercyclic (see Remark \ref{remarkunivsequencebutnothcop}). We will also generalize Theorems \ref{cordoublmapcaspasmmesureintro}, \ref{coralphairrcaspasmmesureintroduction} and \ref{coralphabZalpharotintroduction} for some other ergodic transformations on $(\T,m)$ (see Section \ref{sectionproduitsespaceshardy}). A main tool in our study of the linear dynamics of random sequences of products of operators will be Birkhoff's ergodic theorem.

\smallskip

\bth \label{birkhoffthintro}
Let $(\Omg, \AAA, \mu, \ta)$ be a measure-preserving system. For any $\mu$-integrable function $f \in L^1(\Omega)$, there exists a function $g \in L^1(\Omega)$ satisfying $g \circ \ta = g$ and $\displaystyle \int_\Omg g \, d\mu = \int_\Omg f \, d\mu$, such that
\begin{align}
    \frac{1}{N} \displaystyle \sum_{n = 0}^{N-1} f(\ta^{n}\omg) \underset{N\to \infty}{\longrightarrow} g(\omg) \quad \textrm{for almost every $\omg \in \Omg$.}
\end{align}
Moreover, the system is ergodic if and only if $g(\omg) = \int_\Omg f d\mu$ for almost every $\omega$.
\eth

\smallskip

We would like to point out the fact that the random variables $f(\tau^n \omg)$ involved in Theorem \ref{birkhoffthintro} are not necessarily independent. Thus, Birkhoff's ergodic theorem is a generalization of the law of large numbers. If the system is ergodic and if $A \subseteq \Omg$ is a measurable set, then, applying Theorem \ref{birkhoffthintro} to the function $f = \indic_A $, we obtain
\begin{align} \label{birkhoffpourindicatrices}
\frac{1}{N} \, \textrm{Card}\{0 \leq n \leq N-1 : \ta^{n} \omg \in A  \} \underset{N\to \infty}{\longrightarrow} \mu(A),
\end{align}
for almost every $\omega \in \Omega$. If $A_k$ is a subset of $[0,1)$, we denote by $a_k(n,\omega)$ the quantity
$$
a_k(n,\omega) := \textrm{Card}\{0 \leq i \leq n-1 : \tau^i\omega \in A_k \}.
$$

Let $f : \T \to \R$ be a real function and let $\tau$ be a measure-preserving transformation on $(\Omega, \mathcal{A},\mu)$. We denote by 
$$\mathbb{S}_n^\tau f(\omega) := \displaystyle \sum_{i=0}^{n-1} f(\tau^i \omega)
$$ 
the Birkhoff sum associated to $f$ and $\tau$. In particular, if $A_1, A_2$ are two disjoint Borel subsets of $[0,1)$ with $m(A_k) > 0$ for $k = 1,2$, we have
$$
a_1(n,\omega) - \frac{m(A_1)}{m(A_2)} a_2(n,\omega) = \mathbb{S}_n^\tau f(\omega),
$$
where $f$ is the real function on $\T$ defined by $f = \indic_{A_1} - \tfrac{m(A_1)}{m(A_2)} \indic_{A_2}$. This function $f$ is centered, that is, $\displaystyle \int_{\T} f \, dm = 0$. 

\smallskip

When $\Omega = \T$, when $\tau = R_\alpha$ with $\alpha$ an irrational and when $A = [a,b)$ with $0 \leq a < b \leq 1$, the convergence (\ref{birkhoffpourindicatrices}) holds for every $\omega \in \T$. This follows from the uniform distribution mod $1$ of the sequence $(n \alpha)_{n \geq 1}$ (see, for instance, \cite[Pages 49-50]{Pet}).

\subsection{Organization of the paper}
The behavior of Birkhoff sums associated to the function $f := \indic_{A_1} - \frac{m(A_1)}{m(A_2)} \indic_{A_2}$ on $\T$, with $A_1$ and $A_2$ two disjoint Borel subsets of $[0,1)$, will play a crucial role in determining the universality of some random sequences of products of operators $(T_n(.))_{n \geq 1}$. The boundedness of these sums is linked to the so-called coboundary equation in $L^\infty(\T)$ for the function $f$. We begin our study by presenting these results on the Birkhoff sums. We then focus on Birkhoff sums associated to this function $f$ and to the doubling map and to irrational rotations. The question of whether these Birkhoff sums satisfy a central limit theorem will be examined. This will be done in Section \ref{sectionsumbirkhoff}. After this study of the behavior of Birkhoff sums, we will study the linear dynamics of random sequences of products of operators $(T_n(.))_{n \geq 1}$, when the random operator $T(\omg)$ is equal to the adjoint of a multiplication operator $(M_{\phi_1})^*$ for every $\omg \in A_1$ and equal to the adjoint of a multiplication operator $(M_{\phi_2})^*$ for every $\omg \in A_2$, where $A_1, A_2$ are two disjoint Borel subsets of $\T$ such that $A_1 \cup A_2 = [0,1)$ and $m(A_k) > 0$ for $k = 1,2$, and where $\phi_1, \phi_2 \in H^\infty(\D)$ are nonconstant. This will be done in Section \ref{sectionproduitsespaceshardy}. Finally, we study random products of operators on the space of entire functions $H(\C)$. First, we will study the situation where the operators $T(\tau^i \omega), i \geq 0,$ are entire functions of exponential type of the derivation operator $D$, and afterwards, we will consider a case of random products where the operators $T(\tau^i \omega), i \geq 0,$ do not commute. This will be done in Section \ref{sectionproduitsespacesentieres}.

\subsection{Notation and definition} \label{sectionnotetdefinitions}

We put together here some important definitions and notations for the comprehension of this article.

\medskip

\paragraph*{{\textbf{Basic ergodic facts}}}
We denote by $\T$ the set $\T := \R / \Z$, that we will often identify as $[0,1)$, and by $m $ the normalized Lebesgue measure on $\T$. A transformation $\tau : \T \to \T$ is measure-preserving if $m(\tau^{-1}(B)) = m(B)$ for every Borel subset $B$. A measure-preserving transformation $\tau $ is said to be ergodic on $(\T,m)$ if for every Borel subset $A$ of $\T$ such that $\tau^{-1}(A) = A$, we have $m(A) = 0$ or $m(A) = 1$. Moreover, this transformation $\tau$ is mixing whenever 
$$
m(\tau^{-n}(A) \cap B) \underset{n\to \infty}{\longrightarrow} m(A) \, m(B)
$$
for every Borel sets $A,B$ of $\T$, and weakly mixing whenever
$$
\displaystyle \frac{1}{n} \sum_{j=0}^{n-1} \lvert m(\tau^{-j}(A) \cap B) - m(A) m(B) \rvert \underset{n\to \infty}{\longrightarrow} 0
$$
for every Borel subsets $A,B$ of $\T$.
Mixing is a property of asymptotic decorrelation. Every mixing system is easily seen to be weakly mixing and ergodic.  We refer to P. Walters's book \cite{Wal} and to K. Petersen's book \cite{Pet} for background on ergodic theory.

\medskip

\paragraph*{\textbf{Functions with bounded variation, continued fractions}}
We denote by $BV$ the space of real functions on $\T$ with bounded variation, and by $BV_0$ the space of real centered functions in $BV$, that is, the functions $f$ of $BV$ such that $\int_\T f \, dm = 0$. Any function defined on $\T$ can be seen as a $1$-periodic function on $\R$.  If $f \in BV_0$, we write $\lVert f \rVert_2 := \left( \displaystyle \int_0^1 f(t)^2 \, dt \right)^{1/2}$ for the variance of $f$.
For $f \in BV$, we denote by $c_r(f)$ the complex Fourier coefficients of $f$, with $r \in \Z$. When $f$ is in $BV_0$, we always have
$$
c_r(f) = \frac{\gamma_r(f)}{r}, \quad \textrm{for $r \ne 0$},
$$
where $\gamma_r(f)$ are complex numbers such that $\displaystyle \sup_{r \ne 0} \vert \gamma_r(f) \rvert \leq \frac{V(f)}{2 \pi} < \infty$. 

If $\alpha \in (0,1)$ is an irrational number, we write $\alpha = [0; a_1, a_2,...]$ its continued fraction expansion, where the $a_j$'s are positive integers. We write $(\frac{p_n}{q_n})_{n \geq 0}$ for the sequence of the convergents associated to $\alpha$, where $p_n$ and $q_n$ are positive integers. The sequences $(p_n)_{n \geq 0}$ and $(q_n)_{n \geq 0}$ satisfy the following relations:
$$
\left\{
    \begin{array}{ll}
        & p_0 = 0, \; p_1 = 1, \; q_0 = 1, \; q_1 = a_1\\
        &q_{n+1} = a_{n+1} \, q_n + q_{n-1} \\
        &p_{n+1} = a_{n+1} \, p_n + p_{n-1} \\
        &p_{n-1} \, q_n - p_n \, q_{n-1} = (-1)^n
    \end{array}\;, \; n \geq 1.
\right.
$$
We say that $\alpha$ has bounded partial quotients if the sequence $(a_j)_{j \geq 1}$ is bounded. We refer to Khinchin's book \cite{Khin} for more details on continued fractions.

\medskip

\paragraph*{\textbf{Hardy spaces}}
In what follows, we define Hardy spaces $H^p(\D)$ for $p = 2$ or $p = \infty$, but we can also define them for $1 \leq p \leq \infty$. We will just need these two cases for our study. We refer to Nikolski's book \cite{N} and to Koosis's book \cite{Koo} for background on Hardy spaces. 

\smallskip
The Hardy space $H^2(\D)$ is the space of holomorphic functions $f : \D \to \C$ such that
$$
\displaystyle \sup_{0 \leq r < 1} \frac{1}{2 \pi} \int_{0}^{2\pi} \lvert f(r e^{it})\rvert^2 dt = \lim_{r \nearrow 1} \frac{1}{2 \pi} \int_{0}^{2\pi} \lvert f(r e^{it})\rvert^2 dt  < \infty.
$$
For every $f,g \in H^2(\D)$, we have
$$
\lan f,g \ran = \lim_{r \nearrow 1} \frac{1}{2 \pi} \int_{0}^{2\pi} f(r e^{it}) \, \overline{g(r e^{it})} dt \quad \textrm{and} \quad \lVert f \rVert^2 = \displaystyle \sup_{0 \leq r < 1} \frac{1}{2 \pi} \int_{0}^{2\pi} \lvert f(r e^{it})\rvert^2 dt.
$$
Moreover, $H^2(\D)$ is isometrically isomorphic to the space $\ell_2(\N_0)$.
The space $H^\infty(\D)$ is the space of holomorphic functions $f : \D \to \C$ such that 
$$
\lVert f \rVert_\infty :=  \displaystyle \sup_{z \in \D} \lvert f(z) \rvert < \infty.  
$$
We have $H^\infty(\D) \subset H^2(\D)$. For any $\lambda \in \D$, the function $k_\lambda$ defined by 
$$
k_\lambda(z) = \frac{1}{1 - \overline{\lambda}z}
$$
satisfies
\begin{align}
    f(\lambda) = \lan f, k_\lambda \ran
\end{align}
for every $f \in H^2(\D)$, and is called the reproducing kernel of $H^2(\D)$. If $f \in H^p(\D)$ with $p = 2$ or $p = \infty$, the limits
$$
f^*(e^{it}) := \displaystyle \lim_{r \to 1} f (r e^{it})
$$
exist for almost every $e^{it} \in \T$. We call the function $f^*$ the boundary value of $f$. An analytic function $f : \D \to \C$ is said to be inner if $f$ is bounded on $\D$ and if $\lvert f^* \rvert = 1$ almost everywhere on $\T$. Any inner function $f$ satisfies $\lvert f(z) \rvert \leq 1$, for every $z \in \D$. An analytic function $F : \D \to \C$ is said to be outer if there is a constant $c \in \T$ and a positive function $h$ on $\T$ satisfying $\log h \in L^1(\T)$, such that
$$
F(z) = c \, \exp \left(\frac{1}{2\pi} \int_0^{2\pi} \frac{e^{it}+z}{e^{it} - z} \, \log h(e^{it}) \, dt \right), \quad \textrm{for every} \; z \in \D.
$$

\smallskip
If $F,G$ are outer in $H^2(\D)$, then $FG $ and $1/F$ are still outer. The product of two inner functions of $H^p(\D)$, with $p = 2$ or $p = \infty$, remains an inner function. If $f \in H^p(\D)$ with $p = 2$ or $p = \infty$, the function $Q_f$ defined on $\D$ by
$$
Q_f(z) = \exp \left( \frac{1}{2\pi} \int_{0}^{2\pi} \frac{e^{it} + z}{e^{it} - z} \, \log \lvert f^*(e^{it}) \rvert \, dt \right)
$$
is outer and is called the outer factor of $f$. An important result in Hardy spaces is the inner-outer factorization. This result states that if $f \in H^p(\D) $ with $p = 2$ or $p = \infty$, then $f$ can be written as $f = I F$, with $I$ an inner function, and $F$ an outer function in $H^p(\D)$. Moreover, we can take $F = Q_f$, and the decomposition becomes unique.
As a consequence of the inner-outer decomposition, if $f \in H^\infty(\D)$ and $1/f \in H^\infty(\D)$, then $f $ is outer.

\smallskip

We define multiplication operators on $H^2(\D)$ as follows. If $\phi$ is a function in $H^\infty(\D)$, the multiplication operator $M_\phi$ on $H^2(\D)$ is given by $M_\phi f = f \phi$. These operators satisfy $\lVert M_\phi \rVert = \lVert \phi \rVert_{\infty} = \displaystyle \sup_{z \in \D} \lvert \phi(z) \rvert$. Moreover, we have the relation
\begin{align} \label{eqvpadjointmultiplic}
    M_\phi^* \, k_\lambda = \overline{\phi(\lambda)} \, k_\lambda \quad \textrm{for every $\lambda \in \D$}.
\end{align}
It is well known that a multiplication operator $M_\phi$ is never hypercyclic on $H^2(\D)$, but that its adjoint $M_\phi^*$ is hypercyclic if and only if $\phi(\D) \cap \T \ne \emptyset$, when $\phi \in H^\infty(\D)$ is nonconstant. This result was proved in \cite[Theorem 4.5]{GodShapiro}. See also \cite[Theorem 4.42]{GEPM} and \cite[Example 1.11]{BM}.

\section{Behavior of Birkhoff sums associated to a centered real function} \label{sectionsumbirkhoff}

Birkhoff sums of centered real functions will appear a lot in Section \ref{sectionproduitsespaceshardy}. These Birkhoff sums will influence the dynamical behavior of the random sequence $(T_n(.))_{n \geq 1}$. The aim of Section \ref{sectionsumbirkhoff} is to study the behavior of such random variables. First, we state a criterion for Birkhoff sums to be unbounded on a set of measure $1$, and next we investigate central limit theorems for such sums.

\smallskip

\subsection{Boundedness of Birkhoff sums}

Let $(\Omega, \mathcal{A},\mu, \tau)$ be an ergodic measure-preserving system and let $f \in L_0^1(\Omega)$ be a real centered function. It is known that the sequence of Birkhoff sums $(\mathbb{S}_n^\tau f(\omega))_{n \geq 1}$ changes sign infinitely often for almost every $\omega \in \Omega$, in the sense that it cannot be ultimately positive or negative (\cite[Theorem 4]{Ha}). Also, we have $\displaystyle \liminf_{n \to \infty} \lvert \mathbb{S}_n^\tau f(\omega) \rvert = 0$ for almost every $\omega \in \Omega$ (\cite{Atk}, \cite{S}), and $\displaystyle \limsup_{n \to \infty} \lvert \mathbb{S}_n^\tau f(\omega) \rvert > 0$ for almost every $\omega \in \Omega$, if $f$ is nonzero (\cite[Section 4.1, Remark 1]{Kachu}). More precisely, we have the following interesting result characterizing the boundedness of the sequence $(\mathbb{S}_n^\tau f(\omega))_{n \geq 1}$ for almost every $\omega \in \Omega$.

\smallskip

\bpr[{\cite[Section 4.1, Theorem 19]{Kachu}}]  \label{Propbirkhoffsumbounded}
Let $(\Omega, \mathcal{A},\mu, \tau)$ be an ergodic measure-preserving system and let $f \in L_0^1(\Omega)$ be a real centered function. The sequence $(\mathbb{S}_n^\tau f(\omega))_{n \geq 1}$ is bounded for almost every $\omega \in \Omega$ if and only if the coboundary equation $f = h - h \circ \tau$ has a solution $h$ in $L^\infty(\Omega)$.
\epr

\smallskip

Thus, Proposition \ref{Propbirkhoffsumbounded} implies that 
$$\mu(\{ \omega \in \Omega : \displaystyle \limsup_{n \to \infty} \mathbb{S}_n^\tau f(\omega) = \infty \quad \textrm{or} \quad \liminf_{n \to \infty} \mathbb{S}_n^\tau f(\omega) = -\infty \}) = 1$$ 
when the coboundary equation $f = h - h \circ \tau$ has no solution $h \in L^\infty(\Omega)$, thanks to the ergodicity of $\tau$. However, it is not clear if this implies that 
$$\mu(\{ \omega \in \Omega : \displaystyle \limsup_{n \to \infty} \mathbb{S}_n^\tau f(\omega) = \infty \quad \textrm{and} \quad \liminf_{n \to \infty} \mathbb{S}_n^\tau f(\omega) = -\infty \}) = 1$$
in this case. The aim of the next subsection will be to give a sufficient condition to replace the "or" by an "and" in this set. But first, we would like to mention the two following important results. The first one is due to Adams and Rosenblatt.

\smallskip

\bpr[{\cite[Proposition 3.6]{AdRos}}] \label{thadamsrosenblatt}
    Let $f \in L^{\infty}(\T)$ be a real centered function. Then, there exists an ergodic transformation $\tau$ on $(\T, m)$ such that the equation $f = h - h \circ \tau$ has a solution $h$ in $L^{\infty}(\T)$.
\epr

\smallskip

In the case of a real centered step function $f$ on $\T$, the authors of \cite{AdRos} even provide a description of this transformation $\tau$ (\cite[Subsection 2.3]{AdRos}), of which they construct using Rokhlin towers. The second result studies when the sequence of Birkhoff sums associated to an irrational rotation and to a real step function is bounded at every point. If $f : \T \to \R$ is a step function, we define the number $\delta_f(x)$ for $x \in \T$ by
$$
\delta_f(x) := f(x^+) - f(x^-) = \displaystyle \lim_{\varepsilon \searrow 0} f(x + \varepsilon) - f(x - \varepsilon),
$$
which corresponds to the jump of $f$ at the point $x$. We also write $D_f := \{ x \in \T : \delta_f(x) \ne 0 \}$, which is a finite set, and we denote by $\Delta_f(x)$ the following finite number:
$$
\Delta_f(x) := \displaystyle \sum_{k= -\infty}^{\infty} \delta_f(x + k \alpha).
$$
The second result is due to I. Oren and is the following.

\smallskip

\bth [{\cite[Theorem A]{Oren}}] \label{thorenbirkhoffsumbounded}
Let $\alpha$ be an irrational number, let $\tau = R_\alpha$ be the rotation of parameter $\alpha$ and let $f$ be a real centered step function. Then, the sequence $(\mathbb{S}_n^\tau f(\omega))_{n \geq 1}$ is bounded for some/every $\omega \in \T$ if and only if $\Delta_f(.) \equiv 0$.
\eth

\smallskip

We now apply Theorem \ref{thorenbirkhoffsumbounded} to the function $f = \indic_{A_1} - \tfrac{m(A_1)}{m(A_2)} \indic_{A_2}$, with $A_1 = [0,b)$ and $A_2 = [b,1)$.

\smallskip

\bpr \label{applicationthorenfunctionf}
Let $\alpha$ be an irrational number, let $\tau = R_\alpha$ and let $f$ be the real function 
$$
f = \indic_{A_1} - \tfrac{m(A_1)}{m(A_2)} \indic_{A_2},
$$
where $A_1 = [0,b), A_2 = [b,1)$ and $b \in (0,1)$.
Then, the sequence $(\mathbb{S}_n^\tau f(\omega))_{n \geq 1}$ is bounded for some/every $\omega \in \T$ if and only if $b \in \Z \alpha$.
\epr

\smallskip

\bpf
The function $f$ has two discontinuities: one at $x = 0$ and one at $x = b$. Moreover, we have
\begin{align*}
    &\delta_f(0) = f(0^+)-f(0^-) = 1 - (-\frac{b}{1-b}) = \frac{1}{1-b} \\
    &\delta_f(b) = f(b^+) - f(b^-) = -\frac{b}{1-b} - 1 = -\frac{1}{1-b}.
\end{align*}
Suppose first that $b = l\alpha$ for some $l \in \Z$. If $x \notin \Z \alpha$, then $\Delta_f(x) = 0$. If $x \in \Z \alpha$, let's say $x = a \alpha$ for some $a \in \Z$, then $\Delta_f(x) = \delta_f(0) + \delta_f(b) = 0$. Thus $\Delta(.) \equiv 0$.

\smallskip

Suppose now that $b \notin \Z \alpha$. Then $\Delta_f(0) = \delta_f(0) \ne 0$. This proves Proposition \ref{applicationthorenfunctionf}.
\epf

\smallskip

This shows in particular that when $ b \in \Z \alpha$, for the function $f = \indic_{A_1} - \tfrac{m(A_1)}{m(A_2)} \indic_{A_2}$ with $A_1 = [0,b)$ and $A_2 = [b,1)$, the sequence of Birkhoff sums is bounded at every point. In the case $b \notin \Z \alpha$, we know that at a given point $\omega \in \T$, we have 
$$
\displaystyle \limsup_{n \to \infty} \mathbb{S}_n^\tau f(\omega) = \infty \quad \textrm{or} \quad  \limsup_{n \to \infty} \mathbb{S}_n^\tau f(\omega) = -\infty,
$$
but it does not seem clear whether both conditions can be satisfied simultaneously. 

\subsection{Central limit theorems for ergodic transformations}

In this section, we discuss central limit theorems for the sequence $(\mathbb{S}_n^\tau f)_{n \geq 1}$ of Birkhoff sums associated to a real centered function $f$ on $\T$, where $\tau$ is an ergodic transformation.

\smallskip

Given an ergodic transformation $\tau$ on $(\T, m)$ and a real centered function $f$ on $\T$, we say that the sequence $(\mathbb{S}_n^\tau f)_{n \geq 1}$ satisfies a central limit theorem (CLT) if there exists a sequence $(a_n)_{n \geq 1}$ of positive real numbers with $a_n \underset{n\to \infty}{\longrightarrow}  \infty$ such that
$$
\frac{\mathbb{S}_n^\tau f}{a_n} \underset{n\to \infty}{\longrightarrow} \mathcal{N}(0,1) \quad \textrm{in distribution},
$$
where $\mathcal{N}(0,1)$ is the random variable with density $x \mapsto \frac{1}{\sqrt{2 \pi}} e^{-x^2 / 2}$ on $\R$.
Birkhoff sums satisfying a central limit theorem always satisfy $\displaystyle \limsup_{n \to \infty} \mathbb{S}_n^\tau f(\omega) = \infty$ and $\displaystyle \liminf_{n \to \infty} \mathbb{S}_n^\tau f(\omega) = -\infty$, for almost every $\omega \in \T$.

\smallskip

\bpr \label{tcllimsupliminf}
Let $\tau$ be an ergodic measure-preserving transformation on $(\T, m)$ and let $f \in L_0^1(\T)$ be a real centered function on $\T$. Suppose that $(\mathbb{S}_n^\tau f)_{n \geq 1}$ has a subsequence satisfying a CLT. Then 
$$m(\{ \omega \in \Omega : \displaystyle \limsup_{n \to \infty} \mathbb{S}_n^\tau f(\omega) = \infty \quad \textrm{and} \quad \liminf_{n \to \infty} \mathbb{S}_n^\tau f(\omega) = -\infty \}) = 1.$$
\epr

\smallskip

\bpf
Let $(n_k)_{k \geq 1}$ be a strictly increasing sequence of positive integers such that 
$$
\frac{\mathbb{S}_{n_k}^\tau f}{a_{n_k}} \underset{k\to \infty}{\longrightarrow} \mathcal{N}(0,1) \quad \textrm{in distribution.}
$$

Observe that the set
$$
\mathcal{A}_1 := \{ \omega \in \T : \displaystyle \limsup_{n \to \infty} \mathbb{S}_{n}^\tau f(\omega) = \infty \}
$$
satisfies $m(\mathcal{A}_1) = 0$ or $m(\mathcal{A}_1) =1$. Indeed, this set is easily seen to be $\tau$-invariant.

\smallskip

Now by the convergence in distribution to $\mathcal{N}(0,1)$, there exists an integer $k_1 \geq 1$ and a constant $C > 0$ such that for every $k \geq k_1$, we have
$$
m \left (  \frac{\mathbb{S}_{n_k}^\tau f}{a_{n_k}} \geq 1  \right) \geq C.
$$
For every $k_0 \geq k_1$, let $A_{k_0} := \{ \exists \, k \geq k_0 : \frac{\mathbb{S}_{n_k}^\tau f}{a_{n_k}} \geq 1 \}$. Since the sequence $(A_{k_0})_{k_0 \geq k_1}$ is decreasing, we have $m(\{ \forall \, k_0 \geq k_1, \exists \, k \geq k_0 : \frac{\mathbb{S}_{n_k}^\tau f}{a_{n_k}} \geq 1 \}) \geq C$. 

In particular, this implies $m(\mathcal{A}_1) > 0$, and thus $m(\mathcal{A}_1) = 1$. 

\smallskip

Since
$$
-\frac{\mathbb{S}_{n_k}^\tau f}{a_{n_k}} \underset{k\to \infty}{\longrightarrow} \mathcal{N}(0,1) \quad \textrm{in distribution},
$$
we can also prove that $m(\{ \displaystyle \liminf_{n \to \infty} \mathbb{S}_n^\tau f = -\infty \}) = 1$. This concludes the proof of Proposition \ref{tcllimsupliminf}.
\epf

\smallskip
We will now give examples of such transformations $\tau$. First of all, Birkhoff sums associated to a real centered function $f$ on $\T$ and to an ergodic transformation $\tau$ do not necessarily satisfy a CLT. Indeed, if $\tau = R_\alpha$ is an irrational rotation with parameter $\alpha$ and $f$ is a real centered function on $\T$ with bounded variation, the Denjoy-Koksma inequality (\cite[Theorem 3.1]{Her}) implies that $ \lVert \mathbb{S}_{q_k}^\tau f \rVert_{\infty} \leq V(f)$ along the subsequence of denominators of the convergents of $\alpha$, where $V(f)$ is the total variation of $f$. Thus, the Birkhoff sums $(\mathbb{S}_n^\tau f)_{n \geq 1}$ cannot satisfy a CLT. In fact, for any aperiodic transformation $\tau$ on $(\T, m)$ and for any distribution $\nu$ on $\R$, there is a measurable function $f$ on $\T$ such that $\tfrac{\mathbb{S}_{n}^\tau f}{n}$ converges in distribution to $\nu$. This result is due to J-P. Thouvenot and B. Weiss (\cite{TW}). 

\smallskip

However, under certain assumptions on the function $f$ and for an irrational rotation, there exist subsequences of $(\mathbb{S}_n^\tau f)_{n \geq 1}$ satisfying CLTs. An important result on this subject is the following, due to F. Huveneers (\cite{Huv}). 

\smallskip

\bth[{\cite[Proposition 1]{Huv}}] \label{TCL F.huveneers}
Let $\alpha$ be an irrational number in $(0,1)$, let $\tau = R_\alpha$ and let $f$ be the function
$$
f = \indic_{[0,1/2)} -  \indic_{[1/2,1)}.
$$
Then, there exists a strictly increasing sequence $(n_k)_{k \geq 1}$ of positive integers such that $(\tfrac{\mathbb{S}_{n_k}^\tau f}{\sqrt{k}})_{k \geq 1}$ converges in distribution to $\mathcal{N}(0,1)$.
\eth

\smallskip

In particular:

\smallskip

\bco \label{corTCLrotationsindic01/2}
Let $\alpha$ be an irrational number in $(0,1)$, let $\tau = R_\alpha$ and let $f$ be the function
$$
f = \indic_{[0,1/2)} -  \indic_{[1/2,1)}.
$$
Then,
$$
\displaystyle \limsup_{n \to \infty} \mathbb{S}_n^\tau f(\omega) = \infty \quad \textrm{and} \quad \liminf_{n \to \infty} \mathbb{S}_n^\tau f(\omega) = -\infty 
$$
for almost every $\omega \in \T$.
\eco

\smallskip

It is important to note that Theorem \ref{TCL F.huveneers} and Corollary \ref{corTCLrotationsindic01/2} hold for any irrational $\alpha$ in $(0,1)$.

\smallskip

In \cite{CIB}, J-P. Conze, S. Isola, and S. Le Borgne also studied CLTs for rotations along subsequences of $(\mathbb{S}_n^\tau f)_{n \geq 1}$, when the parameter $\alpha$ of the rotation has unbounded partial quotients, and when $f$ is a certain real centered function with bounded variation. A few years later, J-P. Conze and S. Le Borgne continued the study of CLTs for irrational rotations and for a class of functions with bounded variation, which applies when the partial quotients of $\alpha$ are bounded (\cite{CB}). Their method is close to the method used in \cite{Huv}, and relies on an abstract CLT which is satisfied under some decorrelation conditions. It is also mentioned in \cite{CB} that when $\tau = R_\alpha$ with $\alpha$ an irrational, the behavior of the sequence $(\mathbb{S}_{n}^\tau f)_{n \geq 1}$ depends on the regularity of the function $f$ on $\T$. More precisely, under certain Diophantine conditions on $\alpha \in \R \setminus \Q$, too much regularity on $f$ will force $f$ to be a coboundary in $L^\infty(\T)$. This is one of the reasons why they considered less regular functions, but with bounded variation, as for example step functions.

\smallskip

We explore in more depth the case of irrational rotations and the case of the doubling map on $\T$. Recall that if $\alpha$ is an irrational number, we denote by $(a_j)_{j \geq 1}$ the partial quotients of $\alpha$, and by $(\frac{p_k}{q_k})_{k \geq 0}$ the sequence of convergents of $\alpha$.

\smallskip

\subsection{CLTs for an irrational rotation $R_\alpha$ with $\alpha $ having bounded partial quotients.}

In this subsection, we discuss CLTs for an irrational rotation $R_\alpha$, with $\alpha$ satisfying
\begin{align} \label{hypalphacasbpq}
    \exists \, A > 1,\, \exists \, 0 \leq p < 1/8 \,, \forall \, n \geq 1, \, a_n \leq A \, n^p.
\end{align}
In particular, the case of $\alpha$ having bounded partial quotients is obtained with $p = 0$.







\smallskip

The article \cite{CB} studies CLTs for irrational rotations, when the parameter $\alpha$ satisfies the assumption (\ref{hypalphacasbpq}). 
An important point in the approach \cite{CB} is the control of the variance $\lVert \mathbb{S}_n^\tau f \rVert_2^2$. The authors of \cite{CB} proved that CLTs hold along certain subsets of integers with density one, on which the variance is large enough. In order to give a lower bound on the variance, they make the following assumption on the real centered functions $f$ defined on $\T$.
\begin{align} \label{hypfonctionscasalphabpq}
    \exists \, M, \eta, \theta > 0 \quad \textrm{such that \; Card$\{ 0 \leq j \leq N : a_{j+1} \lvert \gamma_{q_j}(f) \rvert \geq \eta  \} \geq \theta \, N, \; \forall N \geq M.$} 
\end{align}

Under (\ref{hypfonctionscasalphabpq}), and when $\alpha$ satisfies the assumption (\ref{hypalphacasbpq}) (for example, when $\alpha$ has bounded partial quotients), we have the following CLT.

\smallskip

\bth \label{TCLrotationscasalphabpqConze}
Let $f \in BV_0$ satisfy Condition (\ref{hypfonctionscasalphabpq}). Suppose that the irrational number $\alpha$ satisfies (\ref{hypalphacasbpq}), that is,
\begin{align*}
    \exists \, A > 1,\, \exists \, 0 \leq p < 1/8 \,, \forall \, n \geq 1, \, a_n \leq A \, n^p.
\end{align*}
Then, there exists a constant $\eta_0 > 0$ such that the set 
$$
W := \left \{ n \geq 1 : \lVert \mathbb{S}_n^\tau f  \rVert_2 \geq \eta_0 \, (\log m(n))^{-\tfrac{1}{2}} m(n)^{\tfrac{1}{2}}  \right\},
$$
where $m(n)$ is the unique positive integer such that $n \in [q_{m(n)}, q_{m(n) + 1} [$, has density $1$ in $\N$. Moreover, the sequence
$$
\left (\frac{\mathbb{S}_n^\tau f}{\lVert \mathbb{S}_n^\tau f \rVert_2} \right)_{n \in W}
$$
converges in distribution to $\mathcal{N}(0,1)$ as $n \to \infty$ along $W$. 

\smallskip

In particular, if $\alpha$ has bounded partial quotients (that is, when $p = 0$), we can replace $m(n)$ by $\log(n)$.
\eth

\medskip

\begin{remark}
    \begin{enumerate}[(a)]
        \item Under the assumptions of Theorem \ref{TCLrotationscasalphabpqConze}, $m(n)$ is at least of order $\tfrac{\log(n)}{\log (\log(n))}$, up to a bounded factor (see \cite[Remark 3.3]{CB}).
        \smallskip
        \item The condition (\ref{hypalphacasbpq}) holds for almost every irrational $\alpha$ (with respect to the Lebesgue measure): for every $p > 1$, for almost every irrational $\alpha$, there is a finite constant $A(\alpha,p)$ such that $$ a_n \leq A(\alpha,p) \, n^p, \; \textrm{for every $n \geq 1$}. $$ See \cite[Lemma 3.4]{CB} for more details.
    \end{enumerate}
\end{remark}

\smallskip

Let us now apply Theorem \ref{TCLrotationscasalphabpqConze} to the function $f = \indic_{A_1} - \tfrac{m(A_1)}{m(A_2)} \indic_{A_2}$, when $A_1 = [0,b)$ and $A_2 = [b,1)$, and when $b$ is a rational number in $(0,1)$. First, let us remark the following fact.

\smallskip

\begin{fact} \label{factcongruencesq_k}
    Let $(q_j)_{j \geq 0}$ be the denominators of the convergents of an irrational $\alpha$. Let $q' \geq 2$ be an integer. For every $k \geq 0$, if $q_k \equiv 0 \pmod{q'}$, then $q_{k+1} \notequiv 0 \pmod{q'}$.
\end{fact}

\smallskip

\bpf
This follows from the fact that two consecutive $q_k$'s are coprime, since $p_{k-1} \, q_k - p_k \, q_{k-1} = (-1)^k$ for every $k \geq 1$. 
\epf

\smallskip

With this fact, we can prove the following lemma.

\smallskip

\blm \label{lem1tclrotationscasalphabpq}
Let $\alpha$ be an irrational in $(0,1)$ satisfying (\ref{hypalphacasbpq}) and let $\tau = R_\alpha$. Let $b = \frac{p'}{q'}$ be a rational in $(0,1)$ with $p'$ and $q'$ coprime, let $A_1 = [0,b)$ and let $A_2 = [b,1)$. Consider the function $f \in BV_0$ defined by
$$
f = \indic_{A_1} - \tfrac{m(A_1)}{m(A_2)} \indic_{A_2}.
$$
Then, the function $f$ satisfies the assumption (\ref{hypfonctionscasalphabpq}).
\elm

\smallskip

\bpf
The Fourier coefficients of $f$ are given by
$$
c_r(f) = \frac{1}{1-b} \, \frac{1- e^{-2 i \pi r b}}{2 i \pi r} = \frac{1}{1-b} \, \frac{\sin(\pi r \tfrac{p'}{q'})}{\pi r} \, e^{-i \pi r \tfrac{p'}{q'}}
$$
for $r \ne 0$. Thus, for $j \geq 1$,
\begin{align}
 \label{gammaqkf}   \gamma_{q_j}(f) = \frac{1}{1-b} \, \frac{\sin(\pi q_j \tfrac{p'}{q'})}{\pi } \, e^{-i \pi q_j \tfrac{p'}{q'}}.
\end{align}
Let us set $\eta := \displaystyle \min_{1 \leq j \leq q'-1} \tfrac{1}{\pi (1-b)} \lvert \sin(\pi j \tfrac{p'}{q'}) \rvert > 0$. Hence $a_{j+1} \lvert \gamma_{q_j} (f) \rvert \geq \lvert \gamma_{q_j}(f) \rvert \geq \eta$ for every $j \geq 1$ such that $q_j \notequiv 0 \pmod{q'}$. Since $q_j \equiv 0 \pmod{q'}$ implies $q_{j+1} \notequiv 0 \pmod{q'}$ for every $j \geq 1$ by Fact \ref{factcongruencesq_k}, it follows that 
\begin{align*}
    \textrm{$\frac{1}{N}$ Card$\{ 0 \leq j \leq N : a_{j+1} \lvert \gamma_{q_j}(f) \rvert \geq \eta  \} \geq \frac{1}{4}$ \;  for every $N \geq 2$.}
\end{align*}
This shows that $f$ satisfies the assumption (\ref{hypfonctionscasalphabpq}).
\epf

\smallskip

We thus deduce the following result.

\smallskip

\bco \label{corTCLrotirracasbpq}
Let $\alpha$ be an irrational number in $(0,1)$ and let $\tau = R_\alpha$. Suppose that $\alpha$ satisfies the condition (\ref{hypalphacasbpq}), that is, 
$$a_n \leq A \, n^p,\; \textrm{for every $n \geq 1$}, $$
for some $A > 0$ and $0 \leq p < 1/8$.
Let $A_1,A_2$ be the two disjoint intervals of $[0,1)$ with rational endpoints, such that $A_1 \cup A_2 = [0,1)$ and $m(A_k) > 0$ for $k = 1,2$. Let $f \in BV_0$ be the function defined by
$$
f = \indic_{A_1} - \tfrac{m(A_1)}{m(A_2)} \indic_{A_2}.
$$
Then, the sequence $\left(\frac{\mathbb{S}_n^\tau f}{\lVert \mathbb{S}_n^\tau f \rVert_2}\right)_{n \in W}$ converges in distribution to $\mathcal{N}(0,1)$. 

\smallskip

In particular,
$$
\displaystyle \limsup_{n \to \infty} \mathbb{S}_n^\tau f(\omega) = \infty \quad \textrm{and} \quad \liminf_{n \to \infty} \mathbb{S}_n^\tau f(\omega) = -\infty 
$$
for almost every $\omega \in \T$.
\eco

\smallskip

For some intervals $A_1$ and $A_2$ with irrational endpoints, it is also possible to apply Theorem \ref{TCLrotationscasalphabpqConze} to the function $f$.

\smallskip

\bpr \label{applicationTCLfonctionfunifdistbpq}
Let $\alpha$ be an irrational number in $(0,1)$ and let $\tau = R_\alpha$. Suppose that $\alpha$ satisfies Condition (\ref{hypalphacasbpq}). Let $A_1 = [0,b)$ and $A_2 = [b,1)$. Suppose that the sequence $(b \, q_k)_{k \geq 1}$ is uniformly distributed mod $1$, where the $q_k$'s are the denominators of the convergents of $\alpha$. Then, the function $f = \indic_{A_1} - \frac{m(A_1)}{m(A_2)} \indic_{A_2}$ satisfies (\ref{hypfonctionscasalphabpq}). In particular,
$$
\displaystyle \limsup_{n \to \infty} \mathbb{S}_n^\tau f(\omega) = \infty \quad \textrm{and} \quad \liminf_{n \to \infty} \mathbb{S}_n^\tau f(\omega) = -\infty 
$$
for almost every $\omega \in \T$.
\epr

\smallskip

\bpf
The proof follows from \cite[Corollary 3.11]{CB}.
\epf

\smallskip

It also mentioned in \cite[Remark 3.13]{CB} that if $b = \displaystyle \sum_{n \geq 0} b_n q_n \alpha \; \textrm{mod $1$}$  is the Ostrowski expansion of $b$ associated to the denominators of $\alpha$ and if $\displaystyle \lim_{n \to \infty} \frac{\lvert b_n \rvert}{a_{n+1}} = 0$, then $\displaystyle \lim_{k \to \infty} \lVert q_k b \rVert = 0$, where $\lVert x \rVert$ is the distance of $x$ to $\Z$ (see \cite[Proposition 1]{GP}). In particular, under this assumption on $b$ and when the partial quotients of $\alpha$ are bounded, one can check that the sequence $(a_{j+1} \, \lvert \gamma_{q_j}(f) \rvert)_{j \geq 0}$ associated to our favorite function $f$ converges to zero, and Condition (\ref{hypfonctionscasalphabpq}) is not satisfied. We thus cannot apply Theorem \ref{TCLrotationscasalphabpqConze}. This shows the limits of the techniques used in the article \cite{CB}. 

\subsection{CLTs for an irrational rotation $R_\alpha$ with $\alpha $ having unbounded partial quotients.}

In the case where the irrational number $\alpha$ does not satisfy (\ref{hypalphacasbpq}), there are also CLTs for the rotation $R_\alpha$ along subsequences, but the techniques used in \cite{CIB} are quite different. In this subsection, we focus on the case of an irrational number $\alpha$ with unbounded partial quotients, which means that the sequence $(a_j)_{j \geq 1}$ of the convergents of $\alpha$ is unbounded.

If $(t_k)_{k \geq 1}$ is an increasing sequence of positive integers, we set $L_0 := 0$ and $L_n = \displaystyle \sum_{k=1}^n q_{t_k}$ for $n \geq 1$. 
We have the following version of the CLT for an irrational rotation $R_\alpha$ when $\alpha$ is such that the sequence $(a_{j+1})_{j \geq 0}$ is growing fast enough along a subsequence.

\smallskip

\bth \label{TCLrotationsconzealphanonbpq}
Let $(t_k)_{k \geq 1}$ be a strictly increasing sequence of positive integers. Assume the growth condition: there exists $\beta > 1$ such that  $a_{t_k +1} \geq k^\beta$ for every $k \geq 1$. Then, for every function $f$ in $BV_0$ satisfying the condition
\begin{align}
    \label{conditionTCLalphanonbpqfunction} \displaystyle \liminf_{n \to \infty} \frac{1}{n} \displaystyle \sum_{k=1}^n  \displaystyle \sum_{r \ne 0} \frac{\lvert \gamma_{r q_k}(f) \rvert^2}{r^2} > 0,
\end{align}
we have $\lVert \mathbb{S}_{L_n}^\tau f  \rVert_2^2 \sim \displaystyle \sum_{k=1}^n  \displaystyle \sum_{r \ne 0} \frac{\lvert \gamma_{r q_k}(f) \rvert^2}{r^2} $ as $n \to \infty$, and the sequence $\left (\frac{\mathbb{S}_{L_n}^\tau f}{\lVert \mathbb{S}_{L_n}^\tau f  \rVert_2 } \right)_{n \geq 1}$ converges in distribution to $\mathcal{N}(0,1)$.
\eth

\smallskip

\begin{remark}
    \begin{enumerate}[(a)]
        \item If the partial quotients of $\alpha$ are not bounded, we can always find a strictly increasing sequence of positive integers $(t_k)_{k \geq 1}$ such that $a_{t_k + 1} \geq k^\beta$ with $\beta > 1$. But this doesn't mean that the condition (\ref{conditionTCLalphanonbpqfunction}) is satisfied.
        \item Condition (\ref{conditionTCLalphanonbpqfunction}) is satisfied in particular when $\displaystyle \liminf_{n \to \infty} \frac{1}{n} \displaystyle \sum_{k=1}^n \lvert \gamma_{q_{t_k}}(f) \rvert^2 > 0$.
        \item For a real function $f$ on $\T$ which is continuous and piecewise $C^1$, we know that $\lvert \gamma_{r q_n}(f) \rvert^2 = \lvert r \rvert^2 \, q_n^2 \, \lvert c_{r q_n}(f)\rvert^2 = \lvert c_{r q_n}(f') \rvert^2$, which converges to $0$ as $n \to \infty$. In particular, $f$ cannot satisfy the condition (\ref{conditionTCLalphanonbpqfunction}).
    \end{enumerate}
\end{remark}

\smallskip

We now apply Theorem \ref{TCLrotationsconzealphanonbpq} to our favorite function $f$.

\smallskip

\bpr \label{proposappliTCLrotationsalphanonbpq}
Let $\alpha$ be an irrational in $(0,1)$ and let $\tau = R_\alpha$. Let $b = \frac{p'}{q'}$ be a rational in $(0,1)$ with $p'$ and $q'$ coprime, let $A_1 = [0,b)$ and let $A_2 = [b,1)$. Consider the function $f \in BV_0$ defined by
$$
f = \indic_{A_1} - \tfrac{m(A_1)}{m(A_2)} \indic_{A_2}.
$$
Suppose that $a_{t_k+1}\underset{k\to \infty}{\longrightarrow}  \infty $ along a sequence $(t_k)_{k \geq 1}$ of positive integers such that $q_{t_k} \notequiv 0 \pmod{q'}$ for every $k \geq 1$. Then, the sequence $\left (\frac{\mathbb{S}_{n}^\tau f}{\lVert \mathbb{S}_{n}^\tau f  \rVert_2 } \right)_{n \geq 1}$ has a subsequence converging in distribution to $\mathcal{N}(0,1)$ and for which $ \lVert \mathbb{S}_{n}^\tau f  \rVert_2 \underset{n\to \infty}{\longrightarrow}  \infty$.
\epr

\smallskip

\bpf
Let $(t_k)_{k \geq 1}$ be a strictly increasing sequence of positive integers such that $a_{t_k+1}\underset{k\to \infty}{\longrightarrow}  \infty $ and $q_{t_k} \notequiv 0 \pmod{q'}$ for every $k \geq 1$. The coefficients $\gamma_{q_{t_k}}(f)$ are given by (\ref{gammaqkf}). For every $k \geq 1$, we have 
\begin{align*}
    \lvert \gamma_{q_{t_k}}(f) \rvert^2 \geq \frac{1}{(1-b)^2 \pi^2} \displaystyle \min_{1 \leq j \leq q' -1} \lvert \sin(\pi j \tfrac{p'}{q'})\rvert^2 > 0,
\end{align*}
and so Condition (\ref{conditionTCLalphanonbpqfunction}) is clearly satisfied in this case. The CLT follows from Theorem \ref{TCLrotationsconzealphanonbpq}.
\epf

\smallskip

We deduce the following consequence.

\smallskip

\bco \label{cor2TCLalphanonbpqfonctionf}
Let $\alpha$ be an irrational in $(0,1)$ and let $\tau = R_\alpha$. Let $b = \frac{p'}{q'}$ be a rational in $(0,1)$ with $p'$ and $q'$ coprime, let $A_1 = [0,b)$ and let $A_2 = [b,1)$. Consider the function $f \in BV_0$ defined by
$$
f = \indic_{A_1} - \tfrac{m(A_1)}{m(A_2)} \indic_{A_2}.
$$
Suppose that $a_{t_k+1}\underset{k\to \infty}{\longrightarrow}  \infty $ along a sequence $(t_k)_{k \geq 1}$ of positive integers such that $q_{t_k} \notequiv 0 \pmod{q'}$ for every $k \geq 1$. Then, 
$$
\displaystyle \limsup_{n \to \infty} \mathbb{S}_n^\tau f(\omega) = \infty \quad \textrm{and} \quad \liminf_{n \to \infty} \mathbb{S}_n^\tau f(\omega) = -\infty 
$$
for almost every $\omega \in \T$.
\eco

\smallskip

Exactly as in Proposition \ref{applicationTCLfonctionfunifdistbpq}, when the partial quotients of $\alpha$ are unbounded, we have the following result.

\smallskip

\bpr \label{applicationTCLfonctionfunifdistnonbpq}
Let $\alpha$ be an irrational number in $(0,1)$ and $\tau = R_\alpha$. Suppose that there exists a sequence $(t_k)_{k \geq 1}$ of positive integers such that $a_{t_k +1} \geq k^\beta$ for every $k \geq 1$, with $\beta > 1$. Let $A_1 = [0,b)$ and $A_2 = [b,1)$. Suppose that the sequence $(b \, q_{t_k})_{k \geq 1}$ is uniformly distributed mod $1$, where the $q_j$'s are the denominators of the convergents of $\alpha$. Then, the function $f = \indic_{A_1} - \frac{m(A_1)}{m(A_2)} \indic_{A_2}$ satisfies (\ref{conditionTCLalphanonbpqfunction}). In particular,
$$
\displaystyle \limsup_{n \to \infty} \mathbb{S}_n^\tau f(\omega) = \infty \quad \textrm{and} \quad \liminf_{n \to \infty} \mathbb{S}_n^\tau f(\omega) = -\infty 
$$
for almost every $\omega \in \T$.
\epr

\smallskip

\bpf
The proof follows from \cite[Subsection 2.2]{CIB}.
\epf

\smallskip

The same remark as in the case of an irrational number $\alpha$ with bounded partial quotients holds: if $b = \displaystyle \sum_{n \geq 0} b_n q_n \alpha \; \textrm{mod $1$}$  is the Ostrowski expansion of $b$ associated to the denominators of $\alpha$ and if $\displaystyle \lim_{n \to \infty} \frac{\lvert b_n \rvert}{a_{n+1}} = 0$, then Condition (\ref{conditionTCLalphanonbpqfunction}) is not satisfied for our favorite function $f$.

\smallskip

\subsection{CLTs for the doubling map}
For weakly mixing transformations, the situation is considerably different. The first one to consider CLTs for weakly mixing transformations is Kac (\cite{K}), who studied CLTs for the doubling map and for integrable functions on $\T$ whose Fourier coefficients verify a certain growth condition. More precisely, he proved the following result.

\smallskip

\bth \label{TCLKacintro}
Let $f$ be a centered function on $\T$ such that the Fourier coefficients of $f$ satisfy
$$
\exists \, \beta > 1/2, \, \forall \, n \ne 0, \; \lvert c_n(f) \rvert \leq \frac{C}{\lvert n \rvert^\beta}.
$$
Then, there is a constant $\sigma^2 \geq 0$ such that the sequence $(\frac{1}{\sqrt{n}} \displaystyle \sum_{k=0}^{n-1} f(2^k .))_{n\geq 1}$ converges in distribution to $\mathcal{N}(0,\sigma^2)$. Moreover
$$
\sigma^2 = \displaystyle \lim_{n \to \infty} \frac{1}{n}  \left \lVert \sum_{k=0}^{n-1} f(2^k .) \right \rVert^2_2,
$$
and if $\sigma^2 = 0$, then there exists a function $g$ in $L^2(\T)$ such that
\begin{align} \label{equationcoborddoublangle}
    f(t) = g(t) - g(2t), \quad \textrm{for almost every $t \in \T$}.
\end{align}
\eth

\smallskip

Many improvements regarding CLTs for weakly mixing transformations have been obtained later. For a class of piecewise monotonic weakly mixing transformations on $\T$, S. Wong proved that CLTs hold for Hölder continuous functions defined on $[0,1)$ (\cite{Wong}). G. Keller and J. Rousseau-Egele also proved the existence of CLTs for such transformations, when the function $f$ is in $BV$, by using spectral properties of the Perron-Frobenius operator associated to the weakly mixing transformation $\tau$ on $[0,1)$ (\cite{Ke}, \cite{R-E}). J. Rousseau-Egele has even investigated the rate of convergence in the CLTs associated to such piecewise monotonic transformations and also proved a local limit theorem, when the function $f$ is in $BV$. For more information on CLTs for weakly mixing transformations, we refer to Denker's article (\cite{Den}). 

\medskip

We consider now the case of the doubling map on $\T$, that is, $\tau(x) = 2 x \pmod1$. We consider the case of the two intervals $A_1 = [0,b)$ and $A_2 = [b,1)$ of $[0,1)$. We denote by $f$ the following real function in $BV_0$ given by
$$
f = \indic_{A_1} - \tfrac{m(A_1)}{m(A_2)} \indic_{A_2}.
$$
We will show, using Kac's result (Theorem \ref{TCLKacintro}), that the sequence $(\mathbb{S}_n^\tau f)_{n \geq 1}$ satisfies a CLT with $\sigma^2 > 0$, when $b\in (0,1)$. Let us notice that this function $f$ satisfies the assumptions of Theorem \ref{TCLKacintro}. The only thing to check is that $\sigma^2 > 0$, that is, we have to prove that the coboundary equation (\ref{equationcoborddoublangle}) for $f$ and $\tau$ has no solution $g \in L^2(\T)$.

\smallskip

\blm \label{lem1TCLdoublmapcasbdyadique}
Let $\tau$ be the doubling map. Let $A_1, A_2$ be the two intervals of $[0,1)$ defined by $A_1 = [0,b)$ and $A_2 = [b,1)$, where $b$ is a real number in $(0,1)$. Let $f$ be the function  
$$f = \indic_{A_1} - \tfrac{b}{1-b} \indic_{A_2}.$$
Then, the coboundary equation $f = g - g \circ \tau$ has no solution $g$ in $L^2(\T)$.
\elm

\smallskip

\bpf
Let $b \in (0,1)$ and let $q \geq 1$ be an odd integer such that $b \notin \frac{1}{q} \Z$. Suppose that there exists a function $g \in L^2(\T)$ satisfying the coboundary equation. 

\smallskip

First, let us notice that $c_{2k+1}(g \circ \tau) = 0$ for every $k \in \Z$. Indeed,
\begin{align*}
    c_{2k+1}(g \circ \tau) &= \int_0^1 g(2t) \, e^{-2 i \pi (2k+1)t} dt \\
    &=\frac{1}{2} \left( \int_0^1 g(u) \, e^{-i \pi (2k+1)u} du + \int_1^2 g(u) \, e^{-i \pi (2k+1)u} du \right) \\
    &= \frac{1}{2} \left( \int_0^1 g(u) \, e^{-i \pi (2k+1)u} du + \int_0^1 g(u) \, e^{-i \pi (2k+1)u} \, e^{-i\pi(2k+1)} du \right) \\
    &=0.
\end{align*}
We can also check that $c_{2k}(g \circ \tau) = c_k(g)$ for every $k \in \Z$. Thus, the function $g$ satisfies the equations
\begin{align}
   \label{eq1coefffouriergcasdoublangl} &c_{2k+1}(f) = c_{2k+1}(g) \\ 
    &c_{2k}(f) =c_{2k}(g) - c_k(g)  \label{eq2coeffouriergcasdoublangl}
\end{align}
for every $k \in \Z$. 

\smallskip

The Fourier coefficients of $f$ are given by
\begin{align*}
    c_n(f) &= \frac{1}{1-b} \, \frac{1 - e^{-2i \pi n b}}{2 i \pi n} \\
    &= \frac{1}{\pi (1-b)} \, \frac{\sin(\pi n b)}{n} \, e^{-i \pi n b}.
\end{align*}
for $n \ne 0$.
From equation (\ref{eq2coeffouriergcasdoublangl}), we have 
$$
c_{q. 2^j}(f) = c_{q.2^j}(g) - c_{q.2^{j-1}}(g) \quad \textrm{for every $j \geq 1$},
$$
and thus
$$
c_{q.2^k}(g) = \displaystyle \sum_{j=0}^{k} c_{q.2^j}(f) \quad \textrm{for every $k \geq 1$},
$$
since $c_q(g) = c_q(f)$.
In particular
\begin{align*}
    \Im{(c_{q.2^k}(g))} = -\frac{1}{q\pi (1-b)} \displaystyle \sum_{j=0}^{k} 2^{-j} \, (\sin(\pi q 2^j b ))^2 \leq -\frac{1}{q\pi (1-b)} (\sin(\pi q b))^2 < 0
\end{align*}
for every $k \geq 1$, 
and thus
$$
\Im{(c_{q.2^k}(g))}  \nrightarrow 0 \quad \textrm{as $k \to \infty$}.
$$
This proves that such a function $g \in L^2(\T)$ does not exist.
\epf

\smallskip

We thus deduce the following CLT for the doubling map.

\smallskip

\bco \label{corTCLdoublmapapplifcaslim}
Let $\tau$ be the doubling map on $\T$ and let $A_1,A_2$ be two disjoint intervals of $[0,1)$ such that $A_1 \cup A_2 = [0,1)$ and $m(A_k) > 0$ for $k = 1,2$. Let $f \in BV_0$ be the function
$$
f = \indic_{A_1} - \tfrac{m(A_1)}{m(A_2)} \indic_{A_2}.
$$
Then, there exists a real number $\sigma^2 > 0$ such that the sequence $(\frac{1}{\sqrt{n}} \mathbb{S}_n^\tau f)_{n \geq 1}$ converges in distribution to $\mathcal{N}(0,\sigma^2)$. 

\smallskip

In particular,
$$
\displaystyle \limsup_{n \to \infty} \mathbb{S}_n^\tau f(\omega) = \infty \quad \textrm{and} \quad \liminf_{n \to \infty} \mathbb{S}_n^\tau f(\omega) = -\infty 
$$
for almost every $\omega \in \T$.
\eco

\section{Random products of adjoint multipliers on the Hardy space $H^2(\D)$} \label{sectionproduitsespaceshardy}

Historically, the first examples of hypercyclic operators are Birkhoff's operators on the space of entire functions, MacLane's operator on the space of entire functions, and Rolewicz's operators on $X = \ell_p$ spaces, with $1 \leq p < \infty$, or on $X = c_0$. Precisely, Rolewicz's operators are the operators $T = \lambda B$, where $ \lambda \in \C$ is nonzero and $B$ is the backward shift operator on $X$ defined by
$$
B(x_1,x_2, \dotsc) = (x_2, x_3, \dotsc).
$$
These operators are hypercyclic precisely when $\lvert \lambda \rvert > 1$. In the case $\lvert \lambda \rvert \leq 1$, the operator $ T = \lambda B$ is just a contraction, so every orbit of $T$ is bounded. Consequently, the operator $T$ cannot be hypercyclic in this case. In the case $\lvert \lambda \rvert > 1$, an application of the Hypercyclicity Criterion (Theorem \ref{critunivforsequence}) with the set $\mathcal{D}_1 = c_{00}$ of finitely supported sequences as a dense subset of $X$, and with the operators $S_{k} = \lambda^{-k} S^k$, where $S$ the forward shift on $X$, shows that the operator $T = \lambda B$ is topologically mixing on $X$.
The operator $T = \lambda B$ on $X = \ell_2(\N_0)$ is topologically conjugate to the adjoint multiplier operator $(M_{\lambda z})^*$ on $H^2(\D)$. Moreover, we have seen in the introduction that the adjoint of a multiplication operator $(M_\phi)^*$ is hypercyclic if and only if $\phi(\D) \cap \T \ne \emptyset$, when $\phi$ is nonconstant. These operators are classic in linear dynamics and since their linear dynamics is well known, it is natural to investigate the universality of a sequence of random products of operators defined by adjoints of multiplication operators on $H^2(\D)$. This section is devoted to this study. We will see that the study of the universality of random sequences $(T_n(.))_{n \geq 1}$ is not this simple compared to the study of the hypercyclicity of the adjoints of multiplication operators.
\smallskip

It is known that a multiplication operator $M_\phi$ is never hypercyclic on $H^2(\D)$ (see, for instance, \cite[Section 4.4]{GEPM}). For random sequences of products of multiplication operators, we have the following remark. Let us recall that if $A_k$ is a Borel subset of $[0,1)$, we denote by $a_k(n,\omg)$ the number
$$
a_k(n,\omg) := \textrm{Card}\{0 \leq i \leq n-1 : \tau^i \omg \in A_k \}.
$$

\smallskip

\bpr \label{rmknonunivmultiplication}
Let $A_1, A_2$ be two disjoint Borel subsets of $[0,1)$ such that $A_1 \cup A_2 = [0,1)$ and $m(A_k) > 0$ for $k = 1,2$. Let $T(\omega)$ be defined by
\begin{equation} 
    T(\omega) = \left\{
    \begin{array}{ll}
        M_{\phi_1} & \mbox{if} \quad \omega \in A_1 \\
        M_{\phi_2} & \mbox{if} \quad \omega \in A_2
    \end{array},
\right.
\end{equation}
where $\phi_1, \phi_2 \in H^\infty(\D)$ are nonconstant functions. Suppose that the transformation $\tau$ is ergodic on $(\T,m)$. Then, for almost every $\omega \in \T$, the sequence $(T_n(\omega))_{n \geq 1}$ is not universal. 
\epr

\smallskip

\bpf
In this case, the operators $T_n(\omega)$ are given by
$$
T_n(\omega) = (M_{\phi_1})^{a_1(n,\omega)} (M_{\phi_2})^{a_2(n,\omega)},
$$
since the operators $M_{\phi_1}$ and $M_{\phi_2}$ commute. By Theorem \ref{birkhoffthintro} applied to the functions $\indic_{A_1}$ and $\indic_{A_2}$, there exists a subset $E$ of $\T$ with $m(E) = 1$ such that for every $\omg \in E$,
\begin{align} \label{produitaléaopmultiplicationnonunivequation}
    \frac{a_1(n,\omega)}{n} \underset{n\to \infty}{\longrightarrow} m(A_1) \quad \textrm{and} \quad \frac{a_2(n,\omega)}{n} \underset{n\to \infty}{\longrightarrow} m(A_2).
\end{align}
Suppose that there exists a function $h \in H^2(\D)$ such that the set $\{ T_n(\omega)h : n \geq 1\} = \{ \phi_1^{a_1(n,\omega)} \phi_2^{a_2(n,\omega)} h : n \geq 1 \}$ is dense in $H^2(\D)$. Then, for every $\lambda \in \D$, the set $\{ \phi_1(\lambda)^{a_1(n,\omega)} \phi_2(\lambda)^{a_2(n,\omega)} h(\lambda) : n \geq 1 \}$ is dense in $\C$, by continuity of point evaluations. In particular, the function $h$ does not vanish on $\D$. Let us notice that by (\ref{produitaléaopmultiplicationnonunivequation}), we have
\[
\begin{array}{rcl}
   \lvert \phi_1(\lambda) \rvert^{a_1(n,\omega)} \lvert \phi_2(\lambda) \rvert^{a_2(n,\omega)}  &\underset{n \to \infty}{=}& e^{n[\log(\lvert \phi_1(\lambda) \rvert^{m(A_1)} \lvert \phi_2(\lambda) \rvert^{m(A_2)}) + o(1)]}\, .
\end{array}
\]
Thus, if there exists $\lambda \in \D$ such that $\lvert \phi_1(\lambda) \rvert^{m(A_1)} \lvert \phi_2(\lambda) \rvert^{m(A_2)} < 1$, or if there exists $\mu \in \D$ such that $\lvert \phi_1(\mu) \rvert^{m(A_1)} \lvert \phi_2(\mu) \rvert^{m(A_2)} > 1$, we obtain a contradiction. 

\smallskip

Suppose now that we are in the case where $\lvert \phi_1(\lambda) \rvert^{m(A_1)} \lvert \phi_2(\lambda) \rvert^{m(A_2)} = 1$ for every $\lambda \in \D$. If $\phi_1(\D) \cap \T \ne \emptyset$, there exists $\alpha \in \T$ such that $\lvert \phi_1(\alpha) \rvert = 1$, and thus $\lvert \phi_1(\alpha)^{a_1(n,\omega)} \phi_2(\alpha)^{a_2(n,\omega)} h(\alpha) \rvert = \lvert h(\alpha) \rvert$ for every $n \geq 1$, which also contradicts the density of the set $\{ \phi_1(\alpha)^{a_1(n,\omega)} \phi_2(\alpha)^{a_2(n,\omega)} h(\alpha) : n \geq 1 \}$. We thus have $\phi_1(\D) \subset \D$ or $\phi_1(\D) \subset \C \setminus \overline{\D}$, by the open mapping theorem. 

\smallskip

Suppose that we are in the case $\phi_1(\D) \subset \D$. By the density of the set $\{ T_n(\omega)h : n \geq 1\}$ in $H^2(\D)$, we can find a strictly increasing sequence $(n_k)_{k \geq 1}$ of positive integers such that $\phi_1^{a_1(n_k, \omg)} \phi_2^{a_2(n_k, \omg)} h \underset{k\to \infty}{\longrightarrow} z $ uniformly on each compact of $\D$. Let us remark that, for every $n \geq 1$, 
$$\lvert \phi_1(0)\rvert^{a_1(n,\omega)} \lvert \phi_2(0) \rvert^{a_2(n,\omega)} = \lvert \phi_1(0) \rvert^{a_1(n,\omega) - \frac{m(A_1)}{m(A_2)} a_2(n,\omega)},$$
since $\lvert \phi_2(0) \rvert = \lvert \phi_1(0) \rvert^{-\frac{m(A_1)}{m(A_2)}}$. We recognize the Birkhoff sums 
$$a_1(n,\omega) - \frac{m(A_1)}{m(A_2)} a_2(n,\omega) = \mathbb{S}_n^{\tau}f(\omg)$$ associated to the function $f = \indic_{A_1} - \frac{m(A_1)}{m(A_2)} \indic_{A_2}$. Thus, looking at $z = 0$, the sequence $(\lvert \phi_1(0) \rvert^{\mathbb{S}_{n_k}^\tau f(\omega)})_{k \geq 1}$ converges to $0$. This implies that $(\mathbb{S}_{n_k}^\tau f(\omega))_{k \geq 1}$ diverges to $\infty$, because $\lvert \phi_1(0) \rvert < 1$. But now if we look at $z = 1/2$, we obtain a contradiction. In the same way, the case $\phi_1(\D) \subset \C \setminus \overline{\D}$ also leads to a contradiction.
\epf

\smallskip

In this section, we investigate the linear dynamics of random sequences $(T_n(.))_{n \geq 1}$, with
\begin{equation} \label{eqcocycleadjointmultiplic}
    T(\omega) = \left\{
    \begin{array}{ll}
        (M_{\phi_1})^* & \mbox{if} \quad \omega \in A_1 \\
        (M_{\phi_2})^* & \mbox{if} \quad \omega \in A_2
    \end{array},
\right.
\end{equation}
where $\phi_1, \phi_2 \in H^\infty(\D)$ are nonconstant functions, $A_1, A_2$ are two disjoint Borel subsets of $[0,1)$ such that $A_1 \cup A_2 = [0,1)$ and $m(A_k) > 0$ for $k = 1,2$, and where $\tau$ an ergodic measure-preserving transformation on $(\T,m)$. In this case, the operators $T_n(\omega)$ are given by
$$
T_n(\omega) = (M_{\phi_1}^{a_1(n,\omega)})^* \, (M_{\phi_2}^{a_2(n,\omega)})^*,
$$
because the operators $M_{\phi_1}$ and $M_{\phi_2}$ commute.
We will separate the two cases $m(A_1) = m(A_2) = 1/2$ and $m(A_1) \ne m(A_2)$, simply because the first case is simpler than the second case, and also because certain assumptions on the functions $\phi_1$ and $\phi_2$ will appear in the first case but will not appear in the second case. The behavior of the sequence $(T_n(\omega))_{n \geq 1}$ will depend very much on the ergodic transformation $\tau$, as we will see.

\subsection{Case $m(A_1) = m(A_2) = 1/2$}
In this subsection, we consider $A_1$ and $A_2$ two disjoint Borel subsets of $[0,1)$ satisfying $A_1 \cup A_2 = [0,1)$ and $m(A_1) = m(A_2) = 1/2$. 

\smallskip

The computation of $\lVert T_n(\omega) \rVert$ gives a simple obstruction to the universality of the random sequence $(T_n(.))_{n \geq 1}$.

\smallskip

\bpr \label{propcondisuffisnormt_n(omega)notunivmesegal}
Let $\tau$ be an ergodic measure-preserving transformation on $(\T,m)$. Let $\phi_1, \phi_2 \in H^\infty(\D)$ be nonconstant.
Suppose that $\lVert \phi_1 \rVert_\infty \lVert \phi_2 \rVert_\infty < 1$. Then,  for almost every $\omega \in \T$, the sequence $(T_n(\omega))_{n \geq 1}$ is not universal. In particular, the random sequence $(T_n(.))_{n \geq 1}$ is not universal. 
\epr

\smallskip

\bpf
By Theorem \ref{birkhoffthintro} applied to the functions $\indic_{A_1}$ and $\indic_{A_2}$, there exists a Borel subset $E$ of $\T$ with $m(E) = 1$ such that for every $\omg \in E$,
\begin{align} \label{equationnonunivnormecasmmmesure}
    \frac{a_1(n,\omega)}{n} \underset{n\to \infty}{\longrightarrow} 1/2 \quad \textrm{and} \quad \frac{a_2(n,\omega)}{n} \underset{n\to \infty}{\longrightarrow} 1/2.
\end{align}
Let us observe that
\begin{align} \label{eqnormphi1phi2<1}
    \lVert T_n(\omega) \rVert = \lVert \phi_1^{a_1(n,\omega)} \phi_2^{a_2(n,\omega)} \rVert_\infty \leq \lVert \phi_1 \rVert_{\infty}^{a_1(n,\omega)} \lVert \phi_2 \rVert_{\infty}^{a_2(n,\omega)}. 
\end{align}
But now since $\lVert \phi_1 \rVert_\infty \lVert \phi_2 \rVert_\infty < 1$, the upper bound of $(\ref{eqnormphi1phi2<1})$ converges to $0$, since by (\ref{equationnonunivnormecasmmmesure}), we have
\[
\begin{array}{rcl}
   \lVert \phi_1 \rVert_{\infty}^{a_1(n,\omega)} \lVert \phi_2 \rVert_{\infty}^{a_2(n,\omega)}  &\underset{n \to \infty}{=}& e^{\tfrac{1}{2}n[\log(\lVert \phi_1 \rVert_\infty \lVert \phi_2 \rVert_\infty) + o(1)]}\, .
\end{array}
\]
Thus, $\lVert T_n(\omega) \rVert $ converges to $0$, and every orbit under $(T_n(\omega))_{n \geq 1}$ is bounded. In particular, the sequence $(T_n(\omega))_{n \geq 1}$ cannot be universal, and this concludes the proof.
\epf

\smallskip

We will now see a sufficient condition for the random sequence $(T_n(.))_{n \geq 1}$ to be universal. The proof is based on the following useful criterion, which is mainly applied to operators with many eigenvalues. In what follows, if $e \in X$ is an eigenvector of $T$, we denote by $\lambda(T,e)$ the corresponding eigenvalue. The linear span of a set $\mathcal{A}$ of $X$ will be written as $\textrm{span}(\mathcal{A})$.

\smallskip

\bpr[{\cite[Theorem 7]{BeGo}}]  \label{proposcriterevalpropresuniv}
Suppose that $X$ is a separable Fréchet space and $(T_n)_{n \geq 1}$ is a sequence of operators on this space.

Suppose that there are two subsets $\mathcal{A}$ and $\mathcal{B}$ of $X$ satisfying
\begin{enumerate}[(a)]
    \item Every element of $\mathcal{A} \cup \mathcal{B}$ is an eigenvector of $T_{n}$, $n \geq 1$, satisfying $\lvert \lambda(T_{n}, a) \rvert \underset{n \to \infty}{\longrightarrow} 0$ for every $a \in \mathcal{A}$ and $\lvert \lambda(T_{n}, b) \rvert \underset{n \to \infty}{\longrightarrow} \infty$ for every $b \in \mathcal{B}$.
    \item $\textrm{span}(\mathcal{A})$ and $\textrm{span}(\mathcal{B})$ are dense in $X$.
\end{enumerate}
Then, the sequence $(T_n)_{n \geq 1}$ is topologically mixing.
\epr

\smallskip
 The sufficient condition for the universality of the random sequence is the following.

\smallskip

\bth \label{propimagephi1phi2cercle}
Let $\tau$ be an ergodic measure-preserving transformation on $(\T,m)$. Let $\phi_1, \phi_2 \in H^\infty(\D)$ be such that $\phi_1 \phi_2$ is nonconstant on $\D$ and $(\phi_1 \phi_2)(\D) \cap \T \ne \emptyset$. Then, the random sequence $(T_n(.))_{n \geq 1}$ is topologically mixing.
\eth

\smallskip

\bpf
Using the identity ($\ref{eqvpadjointmultiplic}$) on the eigenvectors for the adjoint of a multiplication operator, we get
\begin{align}
    T_n(\omega) k_z = \overline{\phi_1(z)}^{a_1(n,\omega)} \, \overline{\phi_2(z)}^{a_2(n,\omega)} \, k_z ,
\end{align}
for every $z \in \D$. Since the function $\phi_1 \phi_2$ is nonconstant on $\D$, the open mapping theorem implies that there exist $\lambda, \mu \in \D$ such that $\lvert \phi_1 \phi_2 (\lambda) \rvert > 1$ and $\lvert \phi_1 \phi_2 (\mu) \rvert < 1$. Now, it is well-known (see \cite[Lemma 4.39]{GEPM} for example) that if $\Lambda \subset \D$ is a set with an accumulation point in $\D$, then the set $\textrm{span}\{ k_\lambda \, : \lambda \in \Lambda \}$ is dense in $H^2(\D)$. Thus, the sets
\begin{align*}
    \mathcal{A}_1 :=\{ k_z : z \in \D, \lvert \phi_1(z) \phi_2(z) \rvert < 1 \} \quad \textrm{and} \quad \mathcal{A}_2 := \{ k_z : z \in \D, \lvert \phi_1(z) \phi_2(z) \rvert > 1 \}
\end{align*}
are such that $\textrm{span}(\mathcal{A}_1)$ and  $\textrm{span}(\mathcal{A}_2)$ are dense in $H^2(\D)$, since the sets $\{ z \in \D : \lvert \phi_1 \phi_2(z) \rvert < 1 \}$ and $\{ z \in \D : \lvert \phi_1 \phi_2(z) \rvert > 1 \}$ both have an accumulation point in $\D$. Finally, by Theorem \ref{birkhoffthintro}, there exists a Borel subset $E$ of $\T$ with $m(E) = 1$ such that for every $\omg \in E$,
$$
\frac{a_1(n,\omega)}{n} \underset{n\to \infty}{\longrightarrow} 1/2 \quad \textrm{and} \quad \frac{a_2(n,\omega)}{n} \underset{n\to \infty}{\longrightarrow} 1/2.
$$
We thus have
\[
\begin{array}{rcl}
   \lvert \lambda(T_n(\omega), k_z) \rvert  &\underset{n \to \infty}{=}& e^{\tfrac{1}{2}n[\log(\lvert \phi_1(z) \, \phi_2 (z)\rvert) + o(1)]}\,,
\end{array}
\]
and it easily follows that $ \lvert \lambda(T_n(\omega), k_z) \rvert \underset{n\to \infty}{\longrightarrow} 0$ for every $k_z \in \mathcal{A}_1$ and $ \lvert \lambda(T_n(\omega), k_w) \rvert \underset{n\to \infty}{\longrightarrow} \infty$ for every $k_w \in \mathcal{A}_2$. Proposition \ref{proposcriterevalpropresuniv} thus applies, and this concludes the proof of Theorem \ref{propimagephi1phi2cercle}. 
\epf

\smallskip

\begin{remark}
    If the ergodic transformation $\tau$ involved in Theorem \ref{propimagephi1phi2cercle} is an irrational rotation and if $A_1, A_2$ are intervals of $[0,1)$, the sequence $(T_n(\omega))_{n \geq 1}$ is in fact topologically mixing for every $\omega \in \T$. This comes from the uniform distribution of the sequence $(n \alpha)_{n \geq 1}$ mod $1$.
\end{remark}

\smallskip

Proposition \ref{propimagephi1phi2cercle} leads to two limit cases: the case $(\phi_1 \phi_2)(\D) \subset \D$ and the case $(\phi_1 \phi_2)(\D) \subset \C \setminus \overline{\D}$, where $\phi_1 \phi_2$ is nonconstant. We first treat some straightforward situations within these two cases.

\smallskip

\bpr \label{situationphi1phi2DcontenuedansD}
Let $\tau$ be a measure-preserving transformation on $(\T,m)$.
Let $\phi_1, \phi_2 \in H^\infty(\D)$ be such that $\phi_1(\D) \subseteq \overline{\D}$ and $\phi_2(\D) \subseteq \overline{\D}$. Then, for every $\omega \in \T$, the sequence $(T_n(\omega))_{n \geq 1}$ is not universal.
\epr

\smallskip

\bpf
This comes from the fact that $\lVert T_n(\omega) \rVert \leq 1$ for every $\omega \in \T$ and every $n \geq 1$. Every orbit under $(T_n(\omega))_{n \geq 1}$ is bounded in this case and the sequence cannot be universal.
\epf

\smallskip

As a consequence of Proposition \ref{situationphi1phi2DcontenuedansD}, we treat the second easy situation.

\smallskip

\bco \label{corcasphi1phi2exterieurdisk}
Let $\tau$ be a measure-preserving transformation on $(\T,m)$.
Let $\phi_1, \phi_2 \in H^\infty(\D)$ be such that $\phi_1(\D) \subseteq \C \setminus \D$ and $\phi_2(\D) \subseteq \C \setminus \D$. Then, for every $\omega \in \T$, the sequence $(T_n(\omega))_{n \geq 1}$ is not universal.
\eco

\smallskip

\bpf
In this case, the operators $T_n(\omega)$ are invertible, with 
$$T_n(\omega)^{-1} = (M_{\phi_1^{-1}}^*)^{a_1(n,\omega)} \, (M_{\phi_2^{-1}}^*)^{a_2(n,\omega)}. $$ 
Indeed, the functions $1/\phi_1$ and $1/\phi_2$ are well defined and belong to $H^\infty(\D)$ under the assumptions of Corollary \ref{corcasphi1phi2exterieurdisk}, and thus the operators $M_{\phi_1}$ and $M_{\phi_2}$ are invertible, with $M_{\phi_1}^{-1} = M_{1/\phi_1}$ and $M_{\phi_2}^{-1} = M_{1/\phi_2}$.

But since $(\phi_1^{-1})(\D) \subseteq \overline{\D}$ and $(\phi_2^{-1})(\D) \subseteq \overline{\D}$, the sequence $(T_n(\omega)^{-1})_{n \geq 1}$ is not universal, by Proposition \ref{situationphi1phi2DcontenuedansD}. In particular, it cannot be topologically transitive (by Proposition \ref{thbirkhoffsequence}), and the sequence $(T_n(\omega))_{n \geq 1}$ cannot be topologically transitive either. 
\epf

In light of the forthcoming results, it is worth noting that if $\phi_1, \phi_2 \in H^\infty(\D)$ are nonconstant such that $\lVert \phi_1 \rVert_\infty \leq 1$, $\lVert \phi_2 \rVert_\infty \leq 1$ and $\phi_1 \phi_2$ is inner, then $\phi_1$ and $\phi_2$ are both inner. Indeed, in this case, we have that $1 \leq \lvert \phi_1^* \rvert = \frac{1}{ \lvert \phi_2^* \rvert} \leq 1$. Thus $\lvert \phi_1^* \rvert =\lvert \phi_2^* \rvert =1$. The same results holds if $\phi_1 \phi_2$ is inner, $\phi_1(\D) \subseteq \C \setminus \D$ and  $\phi_2(\D) \subseteq \C \setminus \D$.

\smallskip

We now focus on the case where $\phi_1 \phi_2$ is an inner function. Depending on the properties of the ergodic transformation $\tau$ on $(\T,m)$, we obtain sufficient conditions for the random sequence $(T_n(.))_{n \geq 1}$ to be universal. We first make the following observation.

\smallskip

\begin{observation} \label{observationmultiplicinner}
    Let $\varphi \in H^\infty(\D)$ be an inner function. Then the multiplication operator $M_\varphi$ on $H^2(\D)$ is an isometry. Indeed, for every $f \in H^2(\D)$,
    \begin{align*}
        \lVert M_\varphi f \rVert^2 = \frac{1}{2\pi} \displaystyle \int_{0}^{2\pi} \lvert \varphi^*(e^{it}) \rvert^2 \lvert f^*(e^{it}) \rvert^2 dt =  \int_{0}^{2\pi} \lvert f^*(e^{it}) \rvert^2 dt = \lVert f \rVert^2.
    \end{align*}
    In particular, $M_\varphi^* M_\varphi = I$.
\end{observation}

\smallskip

We can now state the first sufficient condition.
\smallskip

\bth \label{propunivphi1phi2intetimagephicasegal}
Let $\tau $ be an ergodic transformation on $(\T,m)$. Let $f$ be the function $f = \indic_{A_1} - \indic_{A_2}$. Suppose that the coboundary equation $f = h - h \circ \tau$ has no solution $h \in L^\infty(\T)$. Let $\phi_1, \phi_2 \in H^\infty(\D)$ be nonconstant, such that $\phi_1 \phi_2$ is inner. 

If $\phi_1(\D) \cap \T \ne \emptyset$ and $1/\phi_1 \in H^\infty(\D)$, or $\phi_2(\D) \cap \T \ne \emptyset$ and $1/\phi_2 \in H^\infty(\D)$, then the sequence $(T_n(.))_{n \geq 1}$ is topologically weakly mixing.
\eth

\smallskip

\bpf
Suppose for example that $\phi_1(\D) \cap \T \ne \emptyset$ and that $1/\phi_1 \in H^\infty(\D)$. Then, since $\phi_1$ is nonconstant, there are two nonempty open sets $U,V$ in $\D$ such that
\begin{align*}
    \lvert \phi_1(z) \rvert < 1 \quad \textrm{if \, $z \in U$} \quad \textrm{and} \quad \lvert \phi_1(w) \rvert > 1 \quad \textrm{if \, $w \in V$.}
\end{align*}
Since $1/\phi_1 \in H^\infty(\D)$, the function $\phi_1$ does not vanish on $U$.  
\smallskip

Since the coboundary equation $f = h - h \circ \tau$ has no solution in $L^\infty(\T)$, and since the function $f$ is a real centered function, there exists (by Proposition \ref{Propbirkhoffsumbounded}) a Borel subset $E \subset \T$ with $m(E) = 1$ such that for every $\omega \in E$,
\begin{align*}
    \displaystyle \limsup_{n \to \infty} \mathbb{S}_n^\tau f(\omega) = \infty \quad \textrm{or} \quad \liminf_{n \to \infty} \mathbb{S}_n^\tau f(\omega) = -\infty .
\end{align*}

Let $\omega \in E$. Suppose first that $\displaystyle \limsup_{n \to \infty} \mathbb{S}_n^\tau f(\omega) = \infty$, that is, there is a strictly increasing sequence $(n_k)_{k \geq 1}$ of positive integers such that $$\mathbb{S}_{n_k}^\tau f(\omega) \underset{k\to \infty}{\longrightarrow} \infty.$$

We will apply the Universality Criterion using the eigenvectors $k_z$ of $T_n(\omega)$. Recall that
$$
T_n(\omega) k_z = \overline{\phi_1(z)}^{a_1(n,\omega)} \, \overline{\phi_2(z)}^{a_2(n,\omega)} \, k_z
$$
for every $z \in \D$. Since $\phi_1 \phi_2$ is inner, we have 
$$\lvert \phi_1(z) \phi_2(z) \rvert \leq 1$$ 
for every $z \in U$. In particular, we have
$$
\lvert \phi_2(z) \rvert^{a_2(n,\omega)} \leq \lvert \phi_1(z) \rvert^{-a_2(n,\omega)}
$$
for every $z \in U$, and thus
$$
\lVert  T_{n_k}(\omega) k_z \rVert \leq \lvert \phi_1(z) \rvert^{a_1(n_k, \omg) - a_2(n_k, \omg)} \, \lVert k_z \rVert
$$
for every $z \in U$. We recognize the Birkhoff sums $\mathbb{S}_n^{\tau}f(\omg)$ associated to the real centered function $f = \indic_{A_1} - \indic_{A_2}$, and we have
\begin{align} \label{ineqpartieZ1phi1phi2intv1}
    \lVert T_{n_k}(\omega) k_z \rVert \leq \lvert \phi_1(z) \rvert^{\mathbb{S}_{n_k}^\tau f(\omega)} \, \lVert k_z \rVert
\end{align}
for every $z \in U$. Let us denote by $Z_1$ the following dense subset of $H^2(\D)$
$$
Z_1 := \textrm{span}\{k_z : z \in U \}.
$$
Then by ($\ref{ineqpartieZ1phi1phi2intv1}$), we have $T_{n_k}(\omega) x \underset{k\to \infty}{\longrightarrow} 0$ for every $x \in Z_1$.

\smallskip

Let us set $\phi_1 \phi_2 = \varphi$. Since $1/\phi_1 \in H^\infty(\D)$, $M_{1/\phi_1}^*$ is a bounded operator on $H^2(\D)$, and we have that
$$
T_n(\omg) = (M_{\phi_1}^*)^{a_1(n,\omg)} (M_{1/\phi_1}^*)^{a_2(n,\omg)} (M_\varphi^*)^{a_2(n,\omg)}.
$$

Let us denote by $Z_2$ the following dense subset of $H^2(\D)$
$$
Z_2 := \textrm{span}\{k_z : z \in V \}.
$$
Using Observation \ref{observationmultiplicinner}, we define an operator $S_{n_k}(\omega)$ on $Z_2$ by setting
\begin{align*}
    S_{n_k}(\omega) k_z := \overline{\phi_1(z)}^{-\mathbb{S}_{n_k}^\tau f(\omega)} \, (M_\varphi)^{a_2(n_k,\omega)} k_z
\end{align*}
for every $z \in V$ and extending it by linearity to $Z_2$. Using the fact that $\varphi$ is inner and using the expression of $T_n(\omg)$ above, we obtain
$$
T_{n_k}(\omega) S_{n_k}(\omega) k_z = k_z
$$
for every $z \in V$, and thus,
$$
T_{n_k}(\omega) S_{n_k}(\omega) x = x
$$
for every $x \in Z_2$. Finally, since $\lVert M_\varphi \rVert \leq 1$, we have
\begin{align*}
    \lVert S_{n_k}(\omega) k_z \rVert \leq \lvert \phi_1(z) \rvert^{-\mathbb{S}_{n_k}^\tau f(\omega)} \, \lVert k_z \rVert.
\end{align*}
It follows that $S_{n_k}(\omega) x \underset{k\to \infty}{\longrightarrow} 0 $ for every $x \in Z_2$. Thus, in this case, Proposition \ref{critunivforsequence} applies.

\smallskip

Now suppose that $\displaystyle \liminf_{n \to \infty} \mathbb{S}_n^\tau f(\omega) = -\infty$, that is, there is a strictly increasing sequence $(m_k)_{k \geq 1}$ of positive integers such that $$\mathbb{S}_{m_k}^\tau f(\omega) \underset{k\to \infty}{\longrightarrow} -\infty.$$ 
We just have to exchange the roles of $Z_1$ and $Z_2$. We will have 
$$
T_{m_k}(\omega) x \underset{k\to \infty}{\longrightarrow} 0
$$
for every $x \in Z_2$. We also define a right inverse $S_{m_k}(\omega)$ on $Z_1$ in the same way as before, and we will have 
$$
S_{m_k}(\omega) x \underset{k\to \infty}{\longrightarrow} 0
$$
for every $x \in Z_1$. Observe that the definition of $S_{m_k}(\omg)$ makes sense in this case, since $\phi_1$ does not vanish on $U$. 

This concludes the proof of Theorem \ref{propunivphi1phi2intetimagephicasegal}.
\epf

\smallskip

In the case where $f = \indic_{A_1} - \indic_{A_2}$ and the sequence $(\mathbb{S}_n^\tau f(\omega))_{n \geq 1}$ has a subsequence satisfying a CLT, we can omit the assumptions on the images of $\phi_1$ and $\phi_2$. The idea is to replace the dense subsets $Z_1$ and $Z_2$ of $H^2(\D)$ for the operators $T_n(\omega)$, which are used to apply the Universality Criterion to the sequence $(T_n(\omg))_{n \geq 1}$, by another dense subset using model spaces.  It can be seen as a particular case of Theorem \ref{propunivhypphi1phi2nonoutercaspasmmesure}, that we will present afterwards.
\smallskip

\bth \label{propunivphi1phi2innerA1A2mmesuresanshypimages}
Let $\tau$ be an ergodic measure-preserving transformation on $(\T,m)$ such that, for $f = \indic_{A_1} - \indic_{A_2}$, the sequence $(\mathbb{S}_n^\tau f(\omega))_{n \geq 1}$ satisfies
$$
\displaystyle \limsup_{n \to \infty} \mathbb{S}_n^\tau f(\omega) = \infty \quad \textrm{and} \quad \liminf_{n \to \infty} \mathbb{S}_n^\tau f(\omega) = - \infty
$$
for almost every $\omega \in \T$. Let $\phi_1, \phi_2 \in H^\infty(\D)$ be two nonconstant functions, such that $\phi_1 \phi_2$ is inner. Suppose that either $1/\phi_1 \in H^\infty(\D)$ and $\phi_2$ is not outer, or $1/\phi_2 \in H^\infty(\D)$ and $\phi_1$ is not outer.

Then, the random sequence $(T_n(.))_{n \geq 1}$ is topologically weakly mixing.
\eth

\smallskip

\bpf
We will again use the Universality Criterion (Proposition \ref{critunivforsequence}).

Let us set $\phi_1 \phi_2 = \varphi$. 
Suppose for example that $1/\phi_1 \in H^\infty(\D)$ and that $\phi_2$ is not outer. Let us denote by $I_2$ the inner part of $\phi_2$, which is nonconstant, and by $F_2 \in H^\infty(\D)$ the outer part of $\phi_2$. Using a similar idea to the proof of \cite[Theorem 3.1]{GGP}, we observe that the subspace
$$
K := \displaystyle \bigcup_{n \geq 1} K_n,
$$
with $K_n$ being the orthogonal of $I_2^n H^2(\D)$ in $H^2(\D)$, is dense in $H^2(\D)$ (this follows from the inner-outer factorization). We will apply the Universality Criterion to the sequence $(T_n(\omega))_{n \geq 1}$ using this dense subspace. 

Let $x \in K$ and let $n_x \geq 1$ be an integer such that $x \in K_{n_x}$. Let us remark that $(M_{\phi_2}^*)^{n_x} x = 0$. Indeed, for every $y \in H^2(\D)$,
\begin{align} \label{equationpartKT_ncasA1A2mmesure}
	\lan (M_{\phi_2}^*)^{n_x} x, y \ran = \lan x, I_2^{n_x} F_2^{n_x} y \ran = 0.
\end{align}
We used the fact that $F_2^{n_x} \, y \in H^2(\D)$ to assert that the dot product in ($\ref{equationpartKT_ncasA1A2mmesure}$) is equal to $0$. This shows in particular that $T_n(\omega) x \underset{n\to \infty}{\longrightarrow} 0$ for every $x \in K$, for almost every $\omg \in \T$.

In order to define a right inverse of $T_n(\omega)$, we will again use the eigenvectors $k_z$ of $T_n(\omega)$. First, let us denote by $E$ a subset of $\T$ with $m(E) = 1$ such that for every $\omega \in E$, 
$$
 \displaystyle \limsup_{n \to \infty} \mathbb{S}_n^\tau f(\omega) = \infty \quad \textrm{and} \quad \liminf_{n \to \infty} \mathbb{S}_n^\tau f(\omega) = -\infty.
$$
Let us fix $\omg \in E$.

Since $\phi_1$ is nonconstant, $\phi_1(\D) \cap \D \ne \emptyset$ or $\phi_1(\D) \cap \C \setminus \overline{\D} \ne \emptyset$. Suppose, to begin with, that $\lvert \phi_1(z) \rvert < 1$ on a nonempty open set $U$ of $\D$. Let us denote by $Z$ the dense subspace of $H^2(\D)$ defined by 
$$
Z := \textrm{span}\{ k_z : z \in U \}.
$$
Let $(n_k)_{k \geq 1}$ be a strictly increasing sequence of positive integers such that $\mathbb{S}_{n_k}^\tau f(\omega) \underset{k\to \infty}{\longrightarrow} -\infty $. Exactly as in the proof of Theorem \ref{propunivphi1phi2intetimagephicasegal}, we define a right inverse $S_{n_k}(\omega)$ on $Z$ by setting
$$
S_{n_k}(\omega) k_z = \overline{\phi_1(z)}^{-\mathbb{S}_{n_k}^\tau f(\omega)} \, (M_\varphi)^{a_2(n_k,\omega)} \, k_z
$$
for every $z \in U$ and extending it by linearity to $Z$. Using the fact that $M_{1/\phi_1}^*$ is bounded on $H^2(\D)$, we have that
$$
T_n(\omg) = (M_{\phi_1}^*)^{a_1(n,\omg)} (M_{1/\phi_1}^*)^{a_2(n,\omg)} (M_{\vph}^*)^{a_2(n,\omg)}.
$$
Since $\vph$ is inner,
we obtain that
$$
T_{n_k}(\omega) S_{n_k}(\omega) x = x
$$
for every $x \in Z$. 

Finally, we have
$$
\lVert S_{n_k}(\omega) k_z \rVert \leq \lvert \phi_1(z) \rvert^{-\mathbb{S}_{n_k}^\tau f(\omega)} \, \lVert k_z \rVert
$$
for every $z \in U$, which proves that $S_{n_k}(\omega) x  \underset{k\to \infty}{\longrightarrow} 0 $ for every $x \in Z$.

If $\phi_1(\D) \cap \D = \emptyset $, we use another nonempty open set $V$ of $\D$ such that $\lvert \phi_1(z) \rvert > 1$ for $z \in V$ and we use the fact that $\displaystyle \limsup_{n \to \infty} \mathbb{S}_n^\tau f(\omega) = \infty $ to define a right inverse $S_n(\omega)$ along an appropriate subsequence on the dense subspace of $H^2(\D)$
$$
Z' := \textrm{span}\{k_z : z \in V \}.
$$
The conclusion then follows in the same way.
This shows in particular the importance to have 
$$
\displaystyle \limsup_{n \to \infty} \mathbb{S}_n^\tau f(\omega) = \infty \quad \textrm{and} \quad \liminf_{n \to \infty} \mathbb{S}_n^\tau f(\omega) = -\infty,
$$
and not just one of the two conditions. This concludes the proof of Theorem \ref{propunivphi1phi2innerA1A2mmesuresanshypimages}.    
\epf

\smallskip

We can deduce from Theorem \ref{propunivphi1phi2innerA1A2mmesuresanshypimages} the following results regarding the cases of the doubling map and the irrational rotations. The first result holds for the doubling map and follows from Corollary \ref{corTCLdoublmapapplifcaslim}.

\smallskip

\bco \label{corunivcasA1A2mmmesuredoublangl}
Let $\tau$ be the doubling map on $\T$. Let $A_1 = [0, 1/2)$ and $A_2 = [1/2,1)$. Let $\phi_1, \phi_2 \in H^\infty(\D)$ be nonconstant, such that $\phi_1 \phi_2$ is inner. Suppose that either $1/\phi_1 \in H^\infty(\D)$ and $\phi_2$ is not outer, or $1/\phi_2 \in H^\infty(\D)$ and $\phi_1$ is not outer. Then, the random sequence $(T_n(.))_{n \geq 1}$ is topologically weakly mixing.
\eco

\smallskip

The second result holds for every irrational rotation and follows from Corollary \ref{corTCLrotationsindic01/2}.

\smallskip

\bco \label{corunivcasA1A2mmmesurerotirra}
Let $\alpha$ be an irrational number in $(0,1)$ and let $\tau = R_\alpha$ on $\T$. Let $A_1 = [0, 1/2)$ and $A_2 = [1/2, 1)$. Let $\phi_1, \phi_2 \in H^\infty(\D)$ be nonconstant, such that $\phi_1 \phi_2$ is inner. Suppose that either $1/\phi_1 \in H^\infty(\D)$ and $\phi_2$ is not outer, or $1/\phi_2 \in H^\infty(\D)$ and $\phi_1$ is not outer. Then, the random sequence $(T_n(.))_{n \geq 1}$ is topologically weakly mixing.
\eco

\smallskip

\begin{remark} \label{remarkunivsequencebutnothcop}
    Obviously, Theorem \ref{propunivphi1phi2innerA1A2mmesuresanshypimages} shows that the converse of Theorem \ref{propunivphi1phi2intetimagephicasegal} is not true in general. Indeed, if the transformation $\tau$ is such that the sequence $(\mathbb{S}_n^\tau f)_{n \geq 1}$ has a subsequence satisfying a CLT, if $\phi_1(z) = \frac{3}{2} + \frac{1}{4} z $ and $\phi_2(z) = z \, (\phi_1(z))^{-1}$ for every $z \in \D$, then the random sequence $(T_n(.))_{n \geq 1}$ is topologically weakly mixing, by Theorem \ref{propunivphi1phi2innerA1A2mmesuresanshypimages}. This follows from that fact that $\phi_2$ is not outer, and that $1/\phi_1 \in H^\infty(\D)$. However, $\phi_1(\D) \cap \T = \emptyset$ and $\phi_2(\D) \cap \T = \emptyset$, so that neither $(M_{\phi_1})^*$ nor $(M_{\phi_2})^*$ is hypercyclic.
\end{remark}

\medskip

We have seen that, in the case where $\phi_1 \phi_2$ is nonconstant and such that $(\phi_1 \phi_2)(\D) \cap \T \ne \emptyset$, the sequence $(T_n(\omega))_{n \geq 1}$ is universal for almost every $\omega \in \T$. Moreover, by Theorem \ref{propunivphi1phi2intetimagephicasegal}, if $\phi_1 \phi_2$ is constant on $\D$ with $\phi_1 \phi_2 \equiv \alpha$, $\alpha \in \T$, and if $\phi_1(\D) \cap \T \ne \emptyset$, then the random sequence $(T_n(.))_{n \geq 1}$ is universal. The following proposition shows that the converse is true in this case.
\smallskip

\bpr \label{propphi1phi2constante}
Let $\tau$ be a measure-preserving transformation on $(\T,m)$.
Let $\phi_1, \phi_2 \in H^\infty(\D)$ be nonconstant, such that $\phi_1 \phi_2 \equiv \alpha$ on $\D$, with $\alpha \in \T$. If the sequence $(T_n(\omega))_{n \geq 1}$ is universal for some $\omega \in \T$, then $\phi_1(\D) \cap \T \ne \emptyset$ (and $\phi_2(\D) \cap \T \ne \emptyset$).
\epr

\smallskip

The proof of Proposition \ref{propphi1phi2constante} relies on an important result of F. León-Saavedra and V. Müller, which we state here for the convenience of the reader.

\smallskip

\bpr [{\cite[Theorem 1]{LM}}]  \label{resultLeonMuller}
Let $\mathcal{M} \subset \bx$ be a semigroup of operators on a Banach space $X$. Let $x \in X$ be such that the set $\{ \mu \, Sx : \mu \in \T, S \in \mathcal{M} \}$ is dense in $X$. Suppose that there exists an operator $T \in \bx$ such that $T^*$ has no eigenvalue and that $T$ satisfies $S T  = T S$ for every $S \in \mathcal{M}$. Then the set $\{ S x : S \in \mathcal{M} \}$ is dense in $X$.
\epr

\smallskip

In Proposition \ref{resultLeonMuller}, a semigroup of operators is a subset $\mathcal{M} \subset \bx$ containing the identity $I$, and such that $ST \in \mathcal{M}$ whenever $S,T \in \mathcal{M}$. We can now prove Proposition \ref{propphi1phi2constante}.

\smallskip

\bpf[Proof of Proposition \ref{propphi1phi2constante}]
Under the assumptions of Proposition \ref{propphi1phi2constante}, $1/\phi_1$ and $1/\phi_2$ belong to $H^{\infty}(\D)$, and $M_{\phi_1}^*$ is invertible. We have
$$
T_n(\omega) = \overline{\alpha}^{a_2(n,\omega)} \, (M_{\phi_1}^*)^{\mathbb{S}_n^\tau f(\omega)},
$$
with $f = \indic_{A_1} - \indic_{A_2}$, and with $\mathbb{S}_n^{\tau}f(\omega) = a_1(n, \omg) - a_2(n, \omg)$ which is an integer. Let $\omega \in \T$ and let $h \in H^2(\D)$ be such that the set
$$
\{ T_n(\omega) h : n \geq 1 \}
$$
is dense in $H^2(\D)$. In particular, the set
$$
\{ \mu \, (M_{\phi_1}^*)^k h : \mu \in \T,\, k \in \Z \}
$$
is dense in $H^2(\D)$. We apply Proposition \ref{resultLeonMuller}. Taking $\{ (M_{\phi_1}^*)^k : k \in \Z \}$ for the semigroup and $T = M_{\phi_1}^*$ for the operator, whose adjoint has no eigenvalue in $H^2(\D)$, we obtain that the set
$$
\{ (M_{\phi_1}^*)^k h : \, k \in \Z \}
$$
is dense in $H^2(\D)$. We will now show that there is a function $g \in H^2(\D)$ such that the set 
$$
\{ (M_{\phi_1}^*)^k g : \, k \geq 0 \}
$$
is dense in $H^2(\D)$. Let us denote by $O_+(h)$ the positive orbit of $h$ and by $O_-(h)$ the negative orbit of $h$, that is,
$$
O_+(h) := \{(M_{\phi_1}^*)^k h : \, k \geq 0 \} \quad \textrm{and} \quad O_-(h) := \{(M_{\phi_1}^*)^k h : \, k \leq 0 \}.
$$

Since the set
$$
\{ (M_{\phi_1}^*)^k h : \, k \in \Z \}
$$
is dense, there is a sequence $(n_k)_{k \geq 1}$ of integers with $\lvert n_k \rvert \underset{k\to \infty}{\longrightarrow} \infty$, such that
$$
(M_{\phi_1}^*)^{n_k} h \underset{k\to \infty}{\longrightarrow} h.
$$
For every integer $j \in \Z$, one has
$$
(M_{\phi_1}^*)^{n_k + j} h \underset{k\to \infty}{\longrightarrow} (M_{\phi_1}^*)^{j} h.
$$
But there are infinitely many positive integers $n_k$ or infinitely many negative integers $n_k$. So,
$$
\{ (M_{\phi_1}^*)^k h : \, k \in \Z \} \subset \overline{O_+(h)} \quad \textrm{or} \quad \{ (M_{\phi_1}^*)^k h : \, k \in \Z \} \subset \overline{O_-(h)}. 
$$
Suppose that we are in the second case. Let $U,V $ be two non-empty open sets in $H^2(\D)$. There exist then two integers $i<j<0$ such that $(M_{\phi_1}^*)^i h \in U$ and $(M_{\phi_1}^*)^j h \in V$. Hence, $(M_{\phi_1}^*)^{j-i}(U) \cap V \ne \emptyset$.

This proves that $M_{\phi_1}^*$ is topologically transitive, and there exists a function $g \in H^2(\D)$ such that the set 
$$
\{ (M_{\phi_1}^*)^k g : \, k \geq 0 \}
$$
is dense in $H^2(\D)$. But now, it is well known that this implies $\phi_1(\D) \cap \T \ne \emptyset$, because the adjoint of a multiplication operator $(M_\phi)^*$ associated to a nonconstant function $\phi \in H^\infty(\D)$ is hypercyclic if and only if $\phi(\D) \cap \T \ne \emptyset$.

This concludes the proof of Proposition \ref{propphi1phi2constante}.
\epf

\smallskip

In the case where $m(A_1) \ne m(A_2)$, we will see that it is more difficult to identify limit cases, compared to the case where $m(A_1) = m(A_2) = 1/2$. However, we can still obtain sufficient conditions, similar to those of Theorems \ref{propunivphi1phi2intetimagephicasegal} and \ref{propunivphi1phi2innerA1A2mmesuresanshypimages}, for the sequence $(T_n(\omega))_{ n \geq 1}$ to be universal for almost every $\omega \in \T$.

\smallskip

\subsection{Case where possibly $m(A_1) \ne m(A_2)$}
In this subsection, we consider $A_1$ and $A_2$ two disjoint Borel subsets of $[0,1)$ satisfying $A_1 \cup A_2 = [0,1)$ and $\min\{m(A_1), m(A_2)\} > 0$. We study the dynamics of $(T_n(\omega))_{n \geq 1}$ for almost every $\omega \in \T$, when $T(\omega)$ is given by $(\ref{eqcocycleadjointmultiplic})$.

\smallskip

Exactly as in Proposition \ref{propcondisuffisnormt_n(omega)notunivmesegal}, the computation of $\lVert T_n(\omega) \rVert$ provides a simple obstruction to the universality of random sequences $(T_n(.))_{ n \geq 1}$.

\smallskip

\bpr \label{propnonunivnormcasA1A2pasmmmesure}
Let $\tau$ be an ergodic measure-preserving transformation on $(\T,m)$. Let $\phi_1, \phi_2 \in H^\infty(\D)$ be nonconstant.
Suppose that $\lVert \phi_1 \rVert_\infty^{m(A_1)} \lVert \phi_2 \rVert_\infty^{m(A_2)} < 1$. Then, the random sequence $(T_n(.))_{n \geq 1}$ is not universal.
\epr

\smallskip

We also have an analog of Theorem \ref{propimagephi1phi2cercle} in the case where $m(A_1)$ is not necessarily equal to $m(A_2)$. This is the following result. 

\smallskip

\bth \label{propimagephi1phi2cerclecasmesnonegales}
Let $\tau$ be an ergodic measure-preserving transformation on $(\T,m)$. Let $\phi_1, \phi_2 \in H^\infty(\D)$ be nonconstant on $\D$ and suppose that there exist $\lambda,\mu \in \D$ such that
\begin{align} 
   \label{eq1critvpscocyclemesnonegales} &\lvert \phi_1(\lambda) \rvert^{m(A_1)} \lvert \phi_2(\lambda)\rvert^{m(A_2)} < 1 \\
 \label{eq2critvpscocyclemesnonegales}  \textrm{and} \quad  &\lvert \phi_1(\mu) \rvert^{m(A_1)} \lvert \phi_2(\mu)\rvert^{m(A_2)} > 1.
\end{align}
Then, the random sequence $(T_n(.))_{n \geq 1}$ is topologically mixing.
\eth

\smallskip

\bpf
The proof is very similar to the proof of Theorem \ref{propimagephi1phi2cercle}. Let $E$ be a subset of $\T$ with $m(E) = 1$, such that
$$
\frac{a_1(n,\omega)}{n} \underset{n\to \infty}{\longrightarrow} m(A_1) \quad \textrm{and} \quad \frac{a_2(n,\omega)}{n} \underset{n\to \infty}{\longrightarrow} m(A_2)
$$
for every $\omega \in E$. We again use the relation on the eigenvectors 
$$
T_n(\omg) k_z = \overline{\phi_1(z)}^{a_1(n,\omg)}  \overline{\phi_2(z)}^{a_2(n,\omg)} k_z
$$
for every $z \in \D$.
We thus have for $n \geq 1,\, z \in \D $ and $\omega \in E$,
\[
\begin{array}{rcl}
   \lvert \lambda (T_n(\omega), k_z) \rvert  &\underset{n \to \infty}{=}& e^{n[\log(\lvert \phi_1(z) \rvert^{m(A_1)} \, \lvert \phi_2 (z)\rvert^{m(A_2)}) + o(1)]}\, .
\end{array}
\]
By (\ref{eq1critvpscocyclemesnonegales}) and (\ref{eq2critvpscocyclemesnonegales}), the following sets
\begin{align*}
    &\mathcal{A}_1 := \{ k_z : z \in \D, \, \lvert \phi_1(z) \rvert^{m(A_1)} \lvert \phi_2(z)\rvert^{m(A_2)} < 1   \}\\
  \textrm{and} \quad  &\mathcal{A}_2 := \{ k_z : z \in \D, \, \lvert \phi_1(z) \rvert^{m(A_1)} \lvert \phi_2(z)\rvert^{m(A_2)} > 1   \}
\end{align*}
are such that $\textrm{span}(\mathcal{A}_1)$ and $\textrm{span}(\mathcal{A}_2)$ are dense in $H^2(\D)$. Finally, $\lvert \lambda(T_n(\omega), k_z) \rvert \underset{n\to \infty}{\longrightarrow} 0 $ for every $k_z \in \mathcal{A}_1$ and $\lvert \lambda(T_n(\omega), k_w) \rvert\underset{n\to \infty}{\longrightarrow} \infty$ for every $k_w \in \mathcal{A}_2 $. Thus Proposition \ref{proposcriterevalpropresuniv} applies, and this concludes the proof of Theorem \ref{propimagephi1phi2cerclecasmesnonegales}.
\epf

\smallskip

\begin{remark}
    Again, if the ergodic transformation $\tau$ involved in Theorem \ref{propimagephi1phi2cerclecasmesnonegales} is an irrational rotation and if $A_1, A_2$ are intervals of $[0,1)$, then the sequence $(T_n(\omg))_{n \geq 1}$ is topologically mixing for every $\omg \in \T$. 
\end{remark}

\smallskip

We would like to point out that Theorem \ref{propimagephi1phi2cerclecasmesnonegales} can also be deduced from \cite[Theorem 3.4]{CH08} in the case where $\phi_1(z) = \lambda_1 z$ and $\phi_2(z) = \lambda_2 z$ for every $z \in \D$, in the setting of i.i.d products.

In the case where $m(A_1) = m(A_2)$, if $\phi_1 \phi_2$ is nonconstant, Conditions ($\ref{eq1critvpscocyclemesnonegales}$) and ($\ref{eq2critvpscocyclemesnonegales}$) are equivalent to $(\phi_1 \phi_2)(\D) \cap \T \ne \emptyset$, by the open mapping theorem. In the situation where $m(A_1) \ne m(A_2)$, it is difficult to find a similar statement. The following two results are the analogs of Theorems \ref{propunivphi1phi2intetimagephicasegal} and \ref{propunivphi1phi2innerA1A2mmesuresanshypimages} respectively. The first statement requires an assumption on the images of $\phi_1$ and of $\phi_2$, as in Theorem \ref{propunivphi1phi2intetimagephicasegal}. The only difference with Theorem \ref{propunivphi1phi2intetimagephicasegal} is that here, we require both $\displaystyle \limsup_{n \to \infty} \mathbb{S}_n^{\tau}f(\omg) = \infty $ and $\displaystyle \liminf_{n \to \infty} \mathbb{S}_n^{\tau}f(\omg) = -\infty $ for almost every $\omg \in \T$, and not just one of the two conditions for almost every $\omega \in \T$.

\smallskip

\bth \label{propunivhypimagesphi1phi2caspasmmesure}
Let $\tau$ be an ergodic measure-preserving transformation on $(\T,m)$, such that the sequence $(\mathbb{S}_n^\tau f(\omega))_{n \geq 1}$ satisfies
$$
\displaystyle \limsup_{n \to \infty} \mathbb{S}_n^\tau f(\omega) = \infty \quad \textrm{and} \quad \liminf_{n \to \infty} \mathbb{S}_n^\tau f(\omega) = - \infty
$$
for almost every $\omega \in \T$, with $f = \indic_{A_1} - \frac{m(A_1)}{m(A_2)} \indic_{A_2}$. Let $\phi_1,\phi_2 \in H^\infty(\D)$ be nonconstant, such that
\begin{align}
 \label{eqfonctionprsqintcaspasmmmes}   \lvert \phi_1^*\rvert^{m(A_1)} \, \lvert \phi_2^*\rvert^{m(A_2)} = 1 \quad \textrm{almost everywhere on $\T$}.
\end{align}
 
If $\phi_1(\D) \cap \T \ne \emptyset$ or $\phi_2(\D) \cap \T \ne \emptyset$, then the random sequence $(T_n(.))_{n \geq 1}$ is topologically weakly mixing.
\eth 

\smallskip

\bpf
Let us denote by $\phi_1 = I_1 F_1$ and $\phi_2 = I_2 F_2$ the inner-outer decomposition of $\phi_1$ and $\phi_2$, with $I_1, I_2$ inner and $F_1, F_2$ outer. We can take
$$
F_1(z) =  \exp \left( \frac{1}{2\pi} \int_{0}^{2\pi} \frac{e^{it} + z}{e^{it} - z} \, \log \lvert {\phi_1}^*(e^{it}) \rvert \, dt \right),
$$
and similarly for $F_2$ with $\phi_2$. Let us remark that 
\begin{align}
\label{eqfonctionsextcaspasmmesure}    \lvert F_1 \rvert^{m(A_1)} \lvert F_2 \rvert^{m(A_2)} = 1
\end{align}
on $\D$. Indeed, this follows from the expression of $F_1$ and $F_2$, and from (\ref{eqfonctionprsqintcaspasmmmes}). Moreover, let us remark that
\begin{align}
 \label{ineqpresqintpasmmesure}   \lvert \phi_1(\lambda) \rvert^{m(A_1)} \, \lvert \phi_2(\lambda) \rvert^{m(A_2)} \leq 1
\end{align}
for every $\lambda \in \D$. This is a consequence of the inner-outer decomposition of $\phi_1$ and $\phi_2$, and of (\ref{eqfonctionsextcaspasmmesure}). 

\smallskip

Suppose for example that $\phi_1(\D) \cap \T \ne \emptyset$. Let us denote by $U$ and $V$ two nonempty open sets in $\D$ such that $\lvert \phi_1(z) \rvert < 1$ if $z \in U$ and $\lvert \phi_1(w) \rvert > 1$ if $w \in V$. Since the zeros of $\phi_1$ in $U$ are isolated, we can replace $U$ by $U \setminus \{ z \in \D : \phi_1(z) \ne 0\}$ and suppose that $\phi_1$ does not vanish on $U$. 

\smallskip

We will again use the Universality Criterion and the eigenvectors $k_z$ of $T_n(\omega)$. Let $E$ be a subset of $\T$ with $m(E) = 1$ such that for every $\omega \in E$,
\begin{align*}
    \displaystyle \limsup_{n \to \infty} \mathbb{S}_n^\tau f(\omega) = \infty \quad \textrm{and} \quad \liminf_{n \to \infty} \mathbb{S}_n^\tau f(\omega) = -\infty .
\end{align*}
Let us fix $\omg \in E$.
Let $(n_k)_{k \geq 1}$ be a strictly increasing sequence of positive integers such that $\mathbb{S}_{n_k}^\tau f(\omega) \underset{k\to \infty}{\longrightarrow} \infty $. 

By (\ref{ineqpresqintpasmmesure}) and since
$$
T_n(\omg) k_z = \overline{\phi_1(z)}^{a_1(n,\omg)}  \overline{\phi_2(z)}^{a_2(n,\omg)} k_z,
$$
we have
$$
\lVert T_{n_k}(\omega)  k_z \rVert \leq \lvert \phi_1(z) \rvert^{a_1(n_k, \omg) - \frac{m(A_1)}{m(A_2)} a_2(n_k,\omg)} \, \lVert k_z \rVert
$$
for every $z \in U$, that is
$$
\lVert T_{n_k}(\omega)  k_z \rVert \leq \lvert \phi_1(z) \rvert^{\mathbb{S}_{n_k}^\tau f(\omega)} \, \lVert k_z \rVert
$$
for every $z \in U$, where
$$
a_1(n_k, \omg) - \frac{m(A_1)}{m(A_2)} a_2(n_k, \omg) := \mathbb{S}_{n_k}^{\tau}f(\omg)
$$
are the Birkhoff sums associated to the function $f$. In particular, we have $T_{n_k}(\omega) x \underset{k\to \infty}{\longrightarrow} 0$ for every $x$ in the dense subspace of $H^2(\D)$
$$
Z_1 := \textrm{span}\{k_z : z \in U \}.
$$

We define a right inverse $S_{n_k}(\omega)$ on the dense subspace of $H^2(\D)$
$$
Z_2 := \textrm{span}\{k_z : z \in V \}
$$
by setting
$$
S_{n_k}(\omega) k_z = \overline{F_1(z)}^{-a_1(n_k,\omega)} \, \overline{F_2(z)}^{-a_2(n_k,\omega)} \, (M_{I_1})^{a_1(n_k,\omega)}\, (M_{I_2})^{a_2(n_k,\omega)}\, k_z  
$$
for every $z \in V$ and extending it by linearity to $Z_2$. 
This makes sense, since the operator $T_n(\omg)$ is given by
$$
T_n(\omg) = (M_{F_1}^*)^{a_1(n,\omg)}  (M_{F_2}^*)^{a_2(n,\omg)}  (M_{I_1}^*)^{a_1(n,\omg)}  (M_{I_2}^*)^{a_2(n,\omg)}.  
$$
Since the operators $M_{I_1}$ and $M_{I_2}$ are contractions, and by (\ref{eqfonctionsextcaspasmmesure}), we have
$$
\lVert S_{n_k}(\omega) k_z \rVert \leq \lvert F_1(z) \rvert^{-\mathbb{S}_{n_k}^\tau f(\omega)} \lVert k_z \rVert
$$
for every $z \in V$. Now, we also have $\lvert F_1(z) \rvert > 1$ on $V$, since $\lvert \phi_1(z) \rvert > 1$. Thus $S_{n_k}(\omega) y \underset{k\to \infty}{\longrightarrow} 0$ for every $y \in Z_2$.

\smallskip

In the case where $\phi_1(\D) \cap \T = \emptyset$, that is where $\phi_2(\D) \cap \T \ne \emptyset$, the proof is the same. This is allowed because, if we denote by $g := \indic_{A_2} - \frac{m(A_2)}{m(A_1)} \indic_{A_1}$, then $-\mathbb{S}_{n}^\tau f(\omega) = \frac{m(A_1)}{m(A_2)}\, \mathbb{S}_{n}^\tau g(\omega)$, so the sequence $(\mathbb{S}_{n}^\tau g(\omega))_{n \geq 1}$ has the same limsup and liminf as $(\mathbb{S}_{n}^\tau f(\omega))_{n \geq 1}$, for every $\omega \in E$.

\smallskip

This concludes the proof of Theorem \ref{propunivhypimagesphi1phi2caspasmmesure}.
\epf

\smallskip

The second theorem does not require any hypothesis on the images of $\phi_1$ and of $\phi_2$, as in Theorem \ref{propunivphi1phi2innerA1A2mmesuresanshypimages}, but requires $\phi_1 \phi_2$ not to be outer. It is a more general situation than Theorem \ref{propunivphi1phi2innerA1A2mmesuresanshypimages}, since any function $\phi \in H^\infty(\D)$ such that $1/\phi \in H^\infty(\D)$ is outer.

\smallskip

\bth \label{propunivhypphi1phi2nonoutercaspasmmesure}
Let $\tau$ be an ergodic measure-preserving transformation of $(\T,m)$, such that the sequence $(\mathbb{S}_n^\tau f(\omega))_{n \geq 1}$ satisfies 
$$
\displaystyle \limsup_{n \to \infty} \mathbb{S}_n^\tau f(\omega) = \infty \quad \textrm{and} \quad \liminf_{n \to \infty} \mathbb{S}_n^\tau f(\omega) = - \infty
$$
for almost every $\omega \in \T$, with $f = \indic_{A_1} - \frac{m(A_1)}{m(A_2)} \indic_{A_2}$. Let $\phi_1,\phi_2 \in H^\infty(\D)$ be nonconstant, such that $\phi_1 \phi_2$ is not outer, either $\phi_1$ or $\phi_2$ is not inner, and 
\begin{align}
 \notag   \lvert \phi_1^*\rvert^{m(A_1)} \, \lvert \phi_2^*\rvert^{m(A_2)} = 1 \quad \textrm{almost everywhere on $\T$}.
\end{align} 
Then, the random sequence $(T_n(.))_{n \geq 1}$ is topologically weakly mixing.
\eth 

\smallskip

\bpf
The proof is very similar to the proof of Theorem \ref{propunivphi1phi2innerA1A2mmesuresanshypimages}. 

Let us consider the inner-outer decomposition $\phi_1 = I_1  F_1$ and $\phi_2 = I_2 F_2$ of $\phi_1$ and $\phi_2$, with $I_1, I_2$ inner and $F_1, F_2$ outer. Since $\phi_1 \phi_2$ is not outer, either $\phi_1$ or $\phi_2$ is not outer. Suppose for example that $I_1$ is nonconstant. As in the proof of Theorem \ref{propunivphi1phi2innerA1A2mmesuresanshypimages}, the subspace
$$
K := \displaystyle \bigcup_{n \geq 1} K_n,
$$
with $K_n$ being the orthogonal of $I_1^n H^2(\D)$ in $H^2(\D)$, is dense in $H^2(\D)$, and $T_n(\omega) x \underset{n\to \infty}{\longrightarrow} 0$ for every $x \in K$, for almost every $\omg \in \T$.

\smallskip

Let us suppose first that $\phi_1$ is not inner. In particular, we have that $F_1(\D) \cap \D \ne \emptyset$ or $F_1(\D) \cap \C \setminus \overline{\D} \ne \emptyset$. Suppose for example that we are in the first case, that is, $\lvert F_1(z) \rvert < 1$ on a nonempty open set $U$ of $\D$. 

Let $E$ be a subset of $\T$ with $m(E) = 1$ such that for every $\omega \in E$,
\begin{align*}
    \displaystyle \limsup_{n \to \infty} \mathbb{S}_n^\tau f(\omega) = \infty \quad \textrm{and} \quad \liminf_{n \to \infty} \mathbb{S}_n^\tau f(\omega) = -\infty .
\end{align*}
Let us fix $\omg \in E$. Let $(n_k)_{k \geq 1}$ be a strictly increasing sequence of positive integers such that $\mathbb{S}_{n_k}^\tau f(\omega) \underset{k\to \infty}{\longrightarrow} -\infty $. 

We define a right inverse on the dense subspace of $H^2(\D)$
$$
Z := \textrm{span}\{ k_z : z \in U \}
$$
by setting
$$
S_{n_k}(\omega) k_z = \overline{F_1(z)}^{-a_1(n_k,\omega)} \, \overline{F_2(z)}^{-a_2(n_k,\omega)} \, (M_{I_1})^{a_1(n_k,\omega)}\, (M_{I_2})^{a_2(n_k,\omega)}\, k_z  
$$
for every $z \in U$, and by linear extension to $Z$.
This makes sense, since the operator $T_n(\omg)$ is given by
$$
T_n(\omg) = (M_{F_1}^*)^{a_1(n,\omg)}  (M_{F_2}^*)^{a_2(n,\omg)}  (M_{I_1}^*)^{a_1(n,\omg)}  (M_{I_2}^*)^{a_2(n,\omg)}.  
$$
Since the operators $M_{I_1}$ and $M_{I_2}$ are contractions and since $\lvert F_1\rvert^{m(A_1)} \lvert F_2 \rvert^{m(A_2)} = 1$ on $\D$, we obtain
$$
\lVert S_{n_k}(\omega) k_z \rVert \leq \lvert F_1(z) \rvert^{-\mathbb{S}_{n_k}^\tau f(\omega)} \lVert k_z \rVert
$$
for every $\omega \in E$ and $z \in U$, and thus $S_{n_k}(\omega) y \underset{k\to \infty}{\longrightarrow} 0 $ for every $y \in Z$.

If we have $F_1(\D) \cap \D = \emptyset$, that is, if we have $F_1(\D) \cap \C \setminus{\overline{\D}} \ne \emptyset$, we use the fact that $\displaystyle \limsup_{n \to \infty} \mathbb{S}_n^\tau f(\omega) = \infty$. Finally, if $\phi_1$ is inner and thus $\phi_2$ is not inner, we reason in the same way, using the function $F_2$ and the Birkhoff sums $\mathbb{S}_{n}^\tau g(\omega)$, with $g := \indic_{A_2} - \frac{m(A_2)}{m(A_1)} \indic_{A_1}$.

\smallskip

This concludes the proof of Theorem \ref{propunivhypphi1phi2nonoutercaspasmmesure}.
\epf

\smallskip

We can apply Theorems \ref{propunivhypimagesphi1phi2caspasmmesure} and \ref{propunivhypphi1phi2nonoutercaspasmmesure} to the doubling map and to some irrational rotations. For the doubling map, this is the following result, which follows from Corollary \ref{corTCLdoublmapapplifcaslim}.

\smallskip

\bco \label{cordoublmapcaspasmmesure}
Let $\tau $ be the doubling map on $\T$ and let $A_1$ and $A_2$ be two disjoint intervals of $[0,1)$, such that $A_1 \cup A_2 = [0,1)$ and $m(A_k) > 0$ for $k = 1,2$. Let $\phi_1,\phi_2 \in H^\infty(\D)$ be nonconstant. Suppose that one of the two following conditions holds:
\begin{enumerate}[(i)]
    \item  $\lvert \phi_1^*\rvert^{m(A_1)} \, \lvert \phi_2^*\rvert^{m(A_2)} = 1$ almost everywhere on $\T$, and one of the images $\phi_1(\D)$ or $\phi_2(\D)$ meets the unit circle $\T$.
    \item $\lvert \phi_1^*\rvert^{m(A_1)} \, \lvert \phi_2^*\rvert^{m(A_2)} = 1$ almost everywhere on $\T$, $\phi_1 \phi_2$ is not outer and either $\phi_1$ or $\phi_2$ is not inner.
\end{enumerate}
Then, the random sequence $(T_n(.))_{n \geq 1}$ is topologically weakly mixing.
\eco

\smallskip

For irrational rotations, we first have the following result, which follows from Corollary \ref{corTCLrotirracasbpq} and holds for example if the irrational $\alpha$ has bounded partial quotients (in particular, for almost every irrational $\alpha$). 

\smallskip

\bco \label{coralphabpqrotcaspasmmesure}
Let $\alpha$ be an irrational number satisfying (\ref{hypalphacasbpq}), that is,
$$
a_n \leq A \, n^p, \quad \textrm{for every $n \geq 1$ },
$$
for some $A > 0$ and $0 \leq p < 1/8$. Let $\tau = R_\alpha $. Let $A_1$, $A_2$ be two disjoint intervals of $[0,1)$ with rational endpoints, such that $A_1 \cup A_2 = [0,1)$ and $m(A_k) > 0$ for $k = 1,2$. Let $\phi_1,\phi_2 \in H^\infty(\D)$ be nonconstant. Suppose that one of the two following conditions holds:
\begin{enumerate}[(i)]
    \item  $\lvert \phi_1^*\rvert^{m(A_1)} \, \lvert \phi_2^*\rvert^{m(A_2)} = 1$ almost everywhere on $\T$, and the image of either $\phi_1$ or $\phi_2$ intersects $\T$.
    \item  $\lvert \phi_1^*\rvert^{m(A_1)} \, \lvert \phi_2^*\rvert^{m(A_2)} = 1$ almost everywhere on $\T$, $\phi_1 \phi_2$ is not outer and either $\phi_1$ or $\phi_2$ is not inner.
\end{enumerate}
Then, the random sequence $(T_n(.))_{n \geq 1}$ is topologically weakly mixing.
\eco

\smallskip

The case where $\alpha$ as unbounded partial quotients is given by the following result, which follows from Corollary \ref{cor2TCLalphanonbpqfonctionf}.

\smallskip

\bco \label{coralphanonbpqrotcaspasmmesure}
Let $\alpha$ be an irrational in $(0,1)$ and let $\tau = R_\alpha$. Let $b = \frac{p'}{q'}$ be a rational in $(0,1)$ with $p'$ and $q'$ coprime. Let $A_1 = [0,b)$ and let $A_2 = [b,1)$. Suppose that $a_{t_k+1}\underset{k\to \infty}{\longrightarrow}  \infty $ along a sequence $(t_k)_{k \geq 1}$ of positive integers such that $q_{t_k} \notequiv 0 \pmod{q'}$ for every $k \geq 1$.
Let $\phi_1,\phi_2 \in H^\infty(\D)$ be nonconstant. Suppose that one of the two following conditions holds:
\begin{enumerate}[(i)]
    \item  $\lvert \phi_1^*\rvert^{m(A_1)} \, \lvert \phi_2^*\rvert^{m(A_2)} = 1$ almost everywhere on $\T$, and the image of either $\phi_1$ or $\phi_2$ intersects $\T$.
    \item $\lvert \phi_1^*\rvert^{m(A_1)} \, \lvert \phi_2^*\rvert^{m(A_2)} = 1$ almost everywhere on $\T$, $\phi_1 \phi_2$ is not outer and either $\phi_1$ or $\phi_2$ is not inner.
\end{enumerate}
Then, the random sequence $(T_n(.))_{n \geq 1}$ is topologically weakly mixing.
\eco

\smallskip

Finally, we have the following results for some irrational rotations and for some intervals of $[0,1)$, which follows from Propositions \ref{applicationTCLfonctionfunifdistbpq} and \ref{applicationTCLfonctionfunifdistnonbpq}.

\smallskip

\bco \label{coralphaequidistribdynamics}
Let $\alpha$ be an irrational in $(0,1)$ and let $\tau = R_\alpha$. Let $b \in (0,1)$. Let $A_1 = [0,b)$ and let $A_2 = [b,1)$. Suppose that either $a_{t_k+1}\underset{k\to \infty}{\longrightarrow}  \infty $ along a sequence $(t_k)_{k \geq 1}$ of positive integers such that $(b q_{t_k})_{k \geq 1}$ is uniformly distributed mod $1$, or $\alpha$ has bounded partial quotients and the sequence $(b q_k)_{k \geq 1}$ is uniformly distributed mod $1$, where the $q_k$'s are the denominators of the convergents of $\alpha$.
Let $\phi_1,\phi_2 \in H^\infty(\D)$ be nonconstant. Suppose that one of the two following conditions holds:
\begin{enumerate}[(i)]
    \item  $\lvert \phi_1^*\rvert^{m(A_1)} \, \lvert \phi_2^*\rvert^{m(A_2)} = 1$ almost everywhere on $\T$, and the image of either $\phi_1$ or $\phi_2$ intersects $\T$.
    \item $\lvert \phi_1^*\rvert^{m(A_1)} \, \lvert \phi_2^*\rvert^{m(A_2)} = 1$ almost everywhere on $\T$, $\phi_1 \phi_2$ is not outer and either $\phi_1$ or $\phi_2$ is not inner.
\end{enumerate}
Then, the random sequence $(T_n(.))_{n \geq 1}$ is topologically weakly mixing.
\eco

\smallskip

As we can see, the dynamics of the random sequence $(T_n(.))_{n \geq 1}$ is very dependent on the properties of the ergodic measure-preserving transformation $\tau$ on $(\T,m)$. We will now see that there exist ergodic transformations on $(\T,m)$ for which the random sequence $(T_n(.))_{n \geq 1}$ is never universal.

\smallskip

\bpr \label{lemnonunivpourtransfoerg}
Let $\tau$ be a measure-preserving transformation on $(\T,m)$ such that the coboundary equation $f = h - h \circ \tau$ has a solution $h \in L^\infty(\T)$ for $f = \indic_{A_1} - \frac{m(A_1)}{m(A_2)} \indic_{A_2}$.
Let $\phi_1,\phi_2 \in H^\infty(\D)$ be nonconstant, such that
\begin{align}
 \notag   \lvert \phi_1^*\rvert^{m(A_1)} \, \lvert \phi_2^*\rvert^{m(A_2)} = 1 \quad \textrm{almost everywhere on $\T$}.
\end{align}
Then, for almost every $\omega \in \T$, the sequence $(T_n(\omega))_{ n \geq 1}$ is not universal.
\epr

\smallskip

\bpf
Under the assumptions of Proposition \ref{lemnonunivpourtransfoerg}, there exists a Borel subset $E$ of $\T$ with $m(E) = 1$ such that the sequence $(\mathbb{S}_n^\tau f(\omega))_{n \geq 1}$ is bounded when $\omega \in E$. Indeed, one has $\mathbb{S}_n^\tau f(\omega) = h(\omega) - h(\tau^n \omega)$ for every $n \geq 1$. Moreover, we have
$$
T_n(\omega) = (M_{F_1}^*)^{a_1(n,\omega)}  (M_{F_2}^*)^{a_2(n,\omega)} (M_{I_1}^*)^{a_1(n,\omega)} (M_{I_2}^*)^{a_2(n,\omega)},
$$
where $\phi_1 = I_1 F_1$ and $\phi_2 = I_2 F_2$, with $I_1, I_2$ inner and $F_1,F_2$ outer. We have seen in the proof of Theorem \ref{propunivhypimagesphi1phi2caspasmmesure} that
$$
\lvert F_1 \rvert^{m(A_1)} \lvert F_2 \rvert^{m(A_2)} = 1
$$
on $\D$. In particular, this implies that
\begin{align*}
    \lVert F_2 \rVert_{\infty}^{-\frac{m(A_2)}{m(A_1)}} \leq \lvert F_1(z) \rvert \leq  \lVert F_1 \rVert_{\infty}, 
\end{align*}
for every $z \in \D$, and similarly for $F_2$. For $\omega \in E$, we have
\begin{align*}
    \lVert T_n(\omega) \rVert &\leq \displaystyle \sup_{z \in \D} \{ \lvert F_1(z) \rvert^{a_1(n,\omega)} \lvert F_2(z) \rvert^{a_2(n,\omega)} \} \\
    &\leq \displaystyle \sup_{z \in \D} \{ \lvert F_1(z) \rvert^{\mathbb{S}_n^\tau f(\omega)}\}.
\end{align*}
Since the quantity $\lvert F_1(z) \rvert^{\mathbb{S}_n^\tau f(\omega)}$ is bounded from above for every $z \in \D,\, \omega \in E$ and $n \geq 1$ by
$$
\exp(\max\{\mathbb{S}_n^\tau f(\omega) \log(\lVert F_1 \rVert_\infty) , \ - \tfrac{m(A_2)}{m(A_1)} \,\mathbb{S}_n^\tau f(\omega) \log(\lVert F_2 \rVert_\infty )  \}),
$$
it follows that $\lVert T_n(\omega) \rVert$ is also bounded from above for every $n \geq 1$ and $\omega \in E$, since the sequence $(\mathbb{S}_n^\tau f(\omega))_{n \geq 1}$ is bounded.

\smallskip

In particular, every orbit under $(T_n(\omega))_{n \geq 1}$ is bounded and thus cannot be dense in $H^2(\D)$. This concludes the proof of Proposition \ref{lemnonunivpourtransfoerg}.
\epf

\smallskip

\begin{remark} \label{remarkN^2action}
    It is clear that if the random sequence $(T_n(.))_{n \geq 1}$ associated to an ergodic transformation $\tau$ on $\T$ and to two operators $T_1$ and $T_2$ is universal, then the $\N_0^2$-action generated by $T_1, T_2$ is hypercyclic, as soon as $T_1$ and $T_2$ commute. That is, there exists a vector $x \in X$ such that the set $\{ T_1^n \, T_2^m x : n,m \geq 0 \}$ is dense in $X$. The converse is not true in general. Indeed, let us consider $A_1 = [0,1/2), \, A_2 = [1/2,1)$, and let $\tau$ be an ergodic transformation such that the sequence of Birkhoff sums $(a_1(n,\omg) - a_2(n,\omg))_{n \geq 1}$ is bounded for almost every $\omg \in \T$. Let $T_1 = (M_{\phi_1})^*$ and $T_2 = (M_{\phi_2})^*$ on $H^2(\D)$, where $\phi_1(z) = e^z$ and $\phi_2(z) = e^{-z}$ for every $z \in \D$. Then, for every $n,m \geq 0$, $T_1^n \, T_2^m = (M^*_{\phi_1})^{n-m}$. In particular, since the operator $(M_{\phi_1})^*$ is hypercyclic on $H^2(\D)$, the $\N_0^2$-action generated by $T_1$ and $T_2$ is hypercyclic. However, by Proposition \ref{lemnonunivpourtransfoerg}, the random sequence $(T_n(.))_{n \geq 1}$ is not universal.  
\end{remark}

\smallskip

In certain cases, for an irrational rotation, the sequence $(T_n(\omg))_{n \geq 1}$ can be non-universal, for every $\omg \in \T$.

\smallskip

\bco \label{coralphabZalpharot}
Let $\alpha$ be an irrational number in $(0,1)$ and let $\tau = R_\alpha$.
Let $A_1 = [0,b)$ and $A_2 = [b,1)$, with $b \in (0,1)$. Suppose that $b \in \Z \alpha$.
Let $\phi_1, \phi_2 \in H^\infty(\D)$ be nonconstant, such that 
$$
\lvert \phi_1 ^* \rvert^{m(A_1)} \lvert \phi_2^* \rvert^{m(A_2)} = 1 \quad \textrm{almost everywhere on $\T$.}
$$
Then, for every $\omega \in \T$, the sequence $(T_n(\omega))_{n \geq 1}$ is not universal.
\eco

\smallskip

\bpf
This follows from Proposition \ref{lemnonunivpourtransfoerg} since in this case, the sequence $(\mathbb{S}_n^\tau f(\omega))_{n \geq 1}$, where $f = \indic_{A_1} - \tfrac{m(A_1)}{m(A_2)} \indic_{A_2}$, is bounded for every $\omega \in \T$, thanks to Proposition \ref{applicationthorenfunctionf}.
\epf

\smallskip

Thanks to Theorem \ref{thadamsrosenblatt}, we also deduce the following result.

\smallskip

\bco \label{corr1applideadamsrosental}
Let $A_1, A_2$ be two disjoint Borel subsets of $[0,1)$ such that $A_1 \cup A_2 = [0,1)$ and $m(A_k) > 0$ for $k = 1,2$.
Let $\phi_1,\phi_2 \in H^\infty(\D)$ be nonconstant, such that
\begin{align}
 \notag   \lvert \phi_1^*\rvert^{m(A_1)} \, \lvert \phi_2^*\rvert^{m(A_2)} = 1 \quad \textrm{almost everywhere on $\T$}.
\end{align}

There exists an ergodic measure-preserving transformation on $(\T,m)$ such that for almost every $\omega \in \T$, the sequence $(T_n(\omega))_{n \geq 1}$ is not universal.
\eco

\smallskip

In the same line, we have the following result.

\smallskip

\bpr \label{corr2appliadamsrosental}
Let $A_1, A_2$ be two disjoint Borel subsets of $[0,1)$ such that $A_1 \cup A_2 = [0,1)$ and $m(A_k) > 0$ for $k = 1,2$.
Let $\phi_1,\phi_2 \in H^\infty(\D)$ be nonconstant, such that
\begin{align}
 \notag   \lvert \phi_1(z)\rvert^{m(A_1)} \, \lvert \phi_2(z)\rvert^{m(A_2)} > 1 \quad \textrm{for every  $z \in \D$}.
\end{align}

There exists an ergodic measure-preserving transformation on $(\T,m)$ such that for almost every $\omega \in \T$, the sequence $(T_n(\omega))_{n \geq 1}$ is not universal.
\epr

\smallskip

\bpf
Let us observe that the functions $\phi_1, \phi_2$ are bounded below and above on $\D$, and that $1/\phi_1, 1/\phi_2$ belong to $H^\infty(\D)$. Indeed, we have
$$
\lVert \phi_1 \rVert_\infty^{-1} \leq \lvert \phi_1(z) \rvert^{-1} \leq \lVert \phi_2 \rVert_\infty ^{\frac{m(A_2)}{m(A_1)}}
$$
for every $z \in \D$, and similarly for $\phi_2$.

By Proposition \ref{thadamsrosenblatt}, there exists an ergodic measure-preserving transformation $\tau$ on $(\T,m)$ such that the sequence $(\mathbb{S}_n^\tau f(\omega))_{n \geq 1}$ is bounded for every $\omega$ in a certain Borel subset $E$ of $\T$ with $m(E) = 1$, where $f := \indic_{A_1} - \frac{m(A_1)}{m(A_2)} \indic_{A_2}$.

For this transformation $\tau$, the operators $T_n(\omega)$ are invertible, with 
$$
T_n(\omega)^{-1} = (M_{1/\phi_1}^*)^{a_1(n,\omega)}  (M_{1/\phi_2}^*)^{a_2(n,\omega)}.
$$
But
$$
\lVert T_n(\omega)^{-1} \rVert \leq \displaystyle \sup_{z \in \D} \{ \lvert \phi_1(z) \rvert^{-\mathbb{S}_n^\tau f(\omega)} \}.
$$
The same arguments as in the proof of Proposition \ref{lemnonunivpourtransfoerg} show that each orbit under $(T_n(\omega)^{-1})_{n \geq 1}$ is bounded when $\omega \in E$, and this concludes the proof of Proposition \ref{corr2appliadamsrosental}.
\epf

\smallskip

\begin{remark}
    Let $\tau$ be an ergodic measure-preserving transformation on $(\T,m)$. Let $A_1, A_2$ be two disjoint Borel subsets of $[0,1)$ such that $A_1 \cup A_2 = [0,1)$ and $m(A_k) > 0$ for $k=1,2$. Let $\phi_1,\phi_2 \in H^\infty(\D)$ be nonconstant, such that
\begin{align}
    \lvert \phi_1^*\rvert^{m(A_1)} \, \lvert \phi_2^*\rvert^{m(A_2)} = 1 \quad \textrm{almost everywhere on $\T$}.
\end{align}
We have seen in Theorems \ref{propunivhypimagesphi1phi2caspasmmesure} and \ref{propunivhypphi1phi2nonoutercaspasmmesure} that under some assumptions on the sequence of Birkhoff sums associated to the function $f := \indic_{A_1} - \frac{m(A_1)}{m(A_2)} \indic_{A_2}$ and under some assumptions on the functions $\phi_1$ and $\phi_2$, the random sequence $(T_n(.))_{n \geq 1}$ is topologically weakly mixing. Since $\displaystyle \liminf_{n \to \infty} \lvert \mathbb{S}_n^\tau f(\omg) \rvert = 0$ for almost every $\omg \in \T$, there is a subset $E$ of $\T$ with $m(E) = 1$ such that the sequence $(\mathbb{S}_n^\tau f(\omg))_{n \geq 1}$ has a bounded subsequence for every $\omg \in E$. In particular, for this subsequence $(\mathbb{S}_{n_k}^\tau f(\omg))_{k \geq 1}$, the sequence $(T_{n_k}(\omg))_{k \geq 1}$ is not universal by Proposition \ref{lemnonunivpourtransfoerg}, and in particular, the sequence $(T_n(\omg))_{n \geq 1}$ is not topologically mixing (see \cite[Exercise 3.4.5]{GEPM}). 
\end{remark}


\section{Other examples of random products on the space of entire functions}\label{sectionproduitsespacesentieres}

In this section, we study some natural examples of random products of operators defined on the space of entire functions. In particular, we will construct examples of universal random sequences $(T_n(.))_{n \geq 1}$ for which the operators $T(\tau^i \omega)$, $i \geq 0$, do not commute. 

\smallskip

We denote by $H(\C)$ the space of entire functions, equipped with the topology of uniform convergence on compact subsets of $\C$. This topology is generated by the seminorms $(p_n)_{n \geq 1}$, where $p_n(f) = \displaystyle \sup_{ \lvert z \rvert \leq n } \lvert f(z) \rvert$.
We denote by $D$ the derivation operator on $H(\C)$ defined by $D : f \mapsto f'$. Using Cauchy estimates, one can show that $D$ is indeed an operator on $H(\C)$, that is, a continuous linear map on $H(\C)$. 

\smallskip

We will first consider random products generated by entire functions of exponential type of $D$. Let us recall some facts regarding entire functions of exponential type. 

\smallskip

An entire function $\varphi \in H(\C)$ is said to be of exponential type if there exist constants $M, A > 0$ such that 
$$
\lvert \varphi(z) \rvert \leq M e^{A \lvert z \rvert}, \quad \textrm{for every $z \in \C$.}
$$
If $\varphi : \C \to \C$ is an entire function of exponential type, with $\varphi(z) = \displaystyle \sum_{n = 0}^{\infty} a_n z^n$, the expression 
$$\varphi(D) := \displaystyle \sum_{n \geq 0} a_n D^n $$
defines an operator on $H(\C)$ (\cite[Proposition 4.19]{GEPM}). In particular, if $a \ne 0$ is a complex number and if we denote by $T_a : f \mapsto f(. + a)$ the translation by $a$, we have $T_a = \varphi(D)$ with $\varphi(z) = e^{az}$. It is also known that if $T$ is an operator on $H(\C)$ commuting with $D$, then there exists an entire function $\varphi$ of exponential type such that $T = \varphi(D)$. 

\smallskip

The linear dynamics of the operators $D$ and $T_a$, $a \ne 0$, is well known (see, for instance, \cite[Examples 2.20 and 2.21]{GEPM}). It is thus natural to have a look at the linear dynamics of random sequences $(T_n(.))_{n \geq 1}$ composed of random products of the operators $D$ and $T_a$ with $a \ne 0$. This is the aim of the following result.

\smallskip

\bth \label{propositionfoncenttypeexpo}
Let $\ta : \T \to \T$ be an ergodic measure-preserving transformation on $(\T,m)$. Let $A_1, A_2$ be two disjoint Borel subsets of $[0,1)$ such that $A_1 \cup A_2 = [0,1)$ and $m(A_k) > 0$ for $k = 1,2$. Let $\vph_1$ and $\vph_2$ be two nonconstant entire functions of exponential type. Consider the map
$$
T(\omega) := \left\{ \begin{array}{ll}
    \varphi_1(D) & \textrm{if $\omega \in A_1$}   \\
    \varphi_2(D) & \textrm{if $\omega \in A_2$}
\end{array}.
\right. 
$$
Suppose that there exist $z, w \in \C$ such that
\begin{align}
    \lvert \varphi_1(z)\rvert^{m(A_1)} \lvert \varphi_2( z)\rvert^{m(A_2)} &< 1 \label{eqfonctypeexpo1} \\
    \textrm{and} \quad  \lvert \varphi_1(w)\rvert^{m(A_1)} \lvert \varphi_2( w)\rvert^{m(A_2)} &> 1. \label{eqfonctypeexpo2}
\end{align}
Then, the random sequence $(T_n(.))_{n \geq 1}$ is topologically mixing.
\eth

\smallskip

\bpf
Let us set $\varphi_{1}(z) = \displaystyle \sum_{n \geq 0} a_n z^n$ and $\vph_2(z) = \displaystyle \sum_{n \geq 0} b_n z^n$ for $z \in \C$. Since $\vph_1(D)$ and $\vph_2(D)$ commute, one has
$$
T_n(\omega) = (\vph_1(D))^{a_1(n, \omg)} (\vph_2(D))^{a_2(n, \omg)}.
$$
Our goal is to apply Proposition \ref{proposcriterevalpropresuniv}.
It is well known that the functions $e_{\lambda} : z \mapsto e^{\lambda z}$ are eigenvectors for the operator $D$. It is thus natural to use these functions for the proof. More precisely, we have for every $\lambda \in \C$
$$T_n(\omg) e_\lambda = (\vph_1(\lambda))^{a_1(n,\omg)} (\vph_2(\lambda))^{a_2(n,\omg)} e_\lambda.$$
Let us denote by $E$ a Borel subset of $\T$ with $m(E) = 1$, such that 
$$
\frac{a_1(n,\omega)}{n} \underset{n\to \infty}{\longrightarrow} m(A_1) \quad \textrm{and} \quad \frac{a_2(n,\omega)}{n} \underset{n\to \infty}{\longrightarrow} m(A_2)
$$
for every $\omega \in E$.
We thus have,
\[
\begin{array}{rcl}
    \lvert (\vph_1(\lambda))^{a_1(n,\omg)} (\vph_2(\lambda))^{a_2(n,\omg)} \rvert &=& e^{a_1(n,\omg) \log(\lvert \vph_1(\lambda)\rvert) + a_2(n,\omg) \log(\lvert \vph_2(\lambda) \rvert)} \\
    &\underset{n \to \infty}{=}& e^{n [\log(\lvert \vph_1(\lambda) \rvert^{m(A_1)} \lvert \vph_2(\lambda) \rvert^{m(A_2)} ) + o(1)]}.
\end{array}
\]
Under the assumptions of Theorem \ref{propositionfoncenttypeexpo}, the sets
\begin{align*}
    &\Lambda_1 := \{ \lambda \in \C : \lvert \vph_1(\lambda) \rvert^{m(A_1)} \lvert \vph_2(\lambda) \rvert^{m(A_2)} > 1  \} \\ 
  \textrm{and} \quad  &\Lambda_2 := \{ \lambda \in \C : \lvert \vph_1(\lambda) \rvert^{m(A_1)} \lvert \vph_2(\lambda) \rvert^{m(A_2)} < 1  \}
\end{align*}
are nonempty and both have an accumulation point in $\C$. Now let us denote by $\mathcal{A}_1$ and $\mathcal{A}_2$ the sets
$$
\mathcal{A}_1 := \textrm{span}\{ e_\lambda : \lambda \in \Lambda_1 \} \quad \textrm{and} \quad \mathcal{A}_2 := \textrm{span}\{ e_\lambda : \lambda \in \Lambda_2 \}.
$$
These sets are dense in $H(\C)$. Moreover, denoting by $\lambda(T_n(\omg), f)$ the eigenvalue of $T_n(\omg)$ associated to the eigenfunction $f$, we have $\lvert \lambda(T_n(\omega), f) \rvert  \underset{n\to \infty}{\longrightarrow} \infty $ for every $f \in \mathcal{A}_1$ and $\lvert \lambda(T_n(\omega), g) \rvert  \underset{n\to \infty}{\longrightarrow} 0 $ for every $g \in \mathcal{A}_2$, for every $\omega \in E$.

This proves that the sequence $(T_n(\omega))_{n \geq 1}$ is topologically mixing for almost every $\omega \in \T$ and this concludes the proof of Theorem \ref{propositionfoncenttypeexpo}.
\epf

\smallskip

In the case where $m(A_1) = m(A_2) = 1/2$ and where $\varphi_1$ and $\varphi_2$ are two nonconstant functions of exponential type, both Conditions \ref{eqfonctypeexpo1} and \ref{eqfonctypeexpo2} are simultaneously satisfied if and only if $\varphi_1 \varphi_2$ is nonconstant. Indeed, this comes from the fact that the image of a nonconstant entire function is dense in $\C$. In the situation where $\varphi_1 \varphi_2$ is constant and nonzero, we must have $\varphi_1(z) = e^{\alpha z + \beta}$ for some $\alpha \ne 0$ and $\beta \in \C$. Indeed, suppose that $\lvert \varphi_1(z) \rvert \leq M e^{A \lvert z \rvert}$ for every $z \in \C$, where $A,M > 0$. Since $\vph_1$ does not vanish on $\C$, we can write $\vph_1 = e^{\varphi}$ with $\vph$ an entire function. We thus obtain that $\Re(\vph(z)) \leq \log(M) + A \lvert z \rvert$ for every $z \in \C$, and thus $\vph$ is a polynomial of degree at most one. 

In this latter case, we have that $\varphi_1(D) = e^\beta \, T_\alpha $ and $\varphi_2(D) = c \, e^{-\beta} \, T_{-\alpha} $  for some $c \ne 0$. We treat a case in the situation where $\varphi_1 \varphi_2$ is constant on $\C$. This is again linked to the behavior of Birkhoff sums. This shows in particular that a random sequence $(T_n(.))_{n \geq 1}$ can be non-universal even if $T(\omg)$ is hypercyclic for every $\omg \in \T$.

\smallskip

\bpr \label{propexamplecasfoncentieres}
Let $\ta : \T \to \T$ be an ergodic measure-preserving transformation on $(\T,m)$. Let $A_1, A_2$ be two disjoint Borel subsets of $[0,1)$ such that $A_1 \cup A_2 = [0,1)$ and $m(A_k) > 0$ for $k = 1,2$. Let $\vph_1(z) = e^z$ and $\vph_2(z) = e^{-z}$ for every $z \in \C$. Consider the map
$$
T(\omega) := \left\{ \begin{array}{ll}
    \varphi_1(D) & \textrm{if $\omega \in A_1$}   \\
    \varphi_2(D) & \textrm{if $\omega \in A_2$}
\end{array}.
\right. 
$$
\begin{enumerate}
    \item If $m(A_1) \ne m(A_2)$, then $(T_n(\omg))_{n \geq 1}$ is topologically mixing for almost every $\omg \in \T$;
    \item If $m(A_1) = m(A_2) = 1/2$, then $(T_n(\omg))_{n \geq 1}$ is universal for almost every $\omg \in \T$ if and only if the sequence $(\mathbb{S}_n^{\tau}f(\omg))_{n \geq 1}$ is unbounded for almost every $\omg \in \T$, where $f := \indic_{A_1} - \indic_{A_2}$.
\end{enumerate}
\epr

\smallskip

\bpf
In the context of Proposition \ref{propexamplecasfoncentieres}, the operators $T_n(\omg)$ are given by
\[
T_n(\omg) = T_{a_1(n,\omg) - a_2(n,\omg)}.
\]
The first point follows from Theorem \ref{propositionfoncenttypeexpo}. Let us prove the second point. We again use the eigenvectors $e_\lambda$ to obtain that
$$
T_n(\omg) e_\lambda = e^{\lambda \mathbb{S}_n^{\tau}f(\omg)} e_\lambda
$$
for every $\lambda \in \C$. Suppose first that $(\mathbb{S}_n^{\tau}f(\omg))_{n \geq 1}$ is unbounded for almost every $\omg \in \T$. Let us consider a subset $E$ of $\T$ with $m(E) = 1$ such that for every $\omg \in E$, we have that
$$
\displaystyle \liminf_{n \to \infty} \mathbb{S}_n^{\tau}f(\omg) = -\infty \quad \textrm{or} \quad \displaystyle \limsup_{n \to \infty} \mathbb{S}_n^{\tau}f(\omg) = \infty.
$$
Let $\omg \in E$, and suppose, for example, that $\displaystyle \limsup_{n \to \infty} \mathbb{S}_n^{\tau}f(\omg) = \infty$. Let $(n_k)_{k \geq 1}$ be a strictly increasing sequence of positive integers such that $\mathbb{S}_{n_k}^{\tau}f(\omg) \underset{k \to \infty}{\longrightarrow} \infty$. We define a right inverse of $T_{n_k}(\omg)$ on the eigenvectors $e_\lambda$ by setting $S_{n_k}(\omg) e_\lambda =  e^{-\lambda \mathbb{S}_{n_k}^{\tau}f(\omg)} e_\lambda$ for $\Re{\lambda} > 0$ and extending it by linearity to the dense subspace $\textrm{span}\{ e_\lambda : \Re{\lambda} > 0 \}$ of $H(\C)$. We thus have $T_{n_k}(\omg) g\underset{k \to \infty}{\longrightarrow}  0 $ for every $g \in \textrm{span}\{ e_\lambda : \Re{\lambda} < 0 \}$ and $S_{n_k}(\omg) g\underset{k \to \infty}{\longrightarrow}  0 $ for every $g \in \textrm{span}\{ e_\lambda : \Re{\lambda} > 0 \}$, and the Universality Criterion applies. In the case where $\displaystyle  \liminf_{n \to \infty} \mathbb{S}_n^{\tau}f(\omg) = -\infty$, we proceed in the same way. 

Suppose now that $(\mathbb{S}_n^{\tau}f(\omg))_{n \geq 1}$ is bounded for almost every $\omg \in \T$. Let $E$ be a subset of $\T$ with $m(E) = 1$ such that $(\mathbb{S}_n^{\tau}f(\omg))_{n \geq 1}$ is bounded for every $\omg \in E$. Let $\omg \in E$ and suppose that there exists a function $h \in H(\C)$ such that the set $\{ T_n(\omg) h : n \geq 1 \}$ is dense in $H(\C)$. Since $T_n(\omg) h = h(. + \mathbb{S}_n^{\tau}f(\omg))$ and since the map $g \mapsto g(0)$ is continuous on $H(\C)$, the set $\{ h(\mathbb{S}_n^{\tau}f(\omg)) : n \geq 1 \}$ is dense in $\C$. This gives a contradiction, since the set $\{ h(\mathbb{S}_n^{\tau}f(\omg)) : n \geq 1 \}$ is bounded in $\C$. This concludes the proof of Proposition \ref{propexamplecasfoncentieres}.
\epf

We now study an example of random sequences $(T_n(.))_{n \geq 1}$ for which the operators $T(\tau^i \omega),\, i \geq 0$, do not commute.

\smallskip

If $\lambda \in \C$ is nonzero and if $b \in \C$, we denote by $T_{\lambda, b}$ the operator on $H(\C)$ defined by $T_{\lambda, b} \,  f = f (\lambda\, . + b)$. Let us remark that 
\begin{align}
    T_{\lambda,b}^n \, f = f(\lambda^n \, . + b \displaystyle \sum_{k=0}^{n-1} \lambda^k)
\end{align}
for every $n \geq 1$.

\smallskip

\bth \label{propmelangeTlambdabD}
Let $\tau : \T \to \T$ be an ergodic measure-preserving transformation on $(\T,m)$. Let $A_1, A_2$ be two disjoint Borel subsets of $[0,1)$ such that $A_1 \cup A_2 = [0,1)$ and $m(A_k) > 0$ for $k = 1,2$. Consider the map
$$
T(\omega) := \left\{ \begin{array}{ll}
    T_{\lambda, b} & \textrm{if $\omega \in A_1$}   \\
    D & \textrm{if $\omega \in A_2$}
\end{array},
\right. 
$$
with $\lambda \in \C$ nonzero and $b \in \C$.
Then, the random sequence $(T_n(.))_{n \geq 1}$ is topologically mixing.
\eth

\smallskip

\begin{remark}
It is known that the operator $T_{\lambda,b}$ is hypercyclic if and only if $\lambda = 1$ and $b \ne 0$, by \cite[Proposition 2.4]{BeGoMR}.

\smallskip

Thus, Theorem \ref{propmelangeTlambdabD} gives an example of a map $\omg \mapsto T(\omega)$ which is hypercyclic for every $\omg \in A_2$ and not hypercyclic for every $\omg \in A_1$, for which the random sequence $(T_n(.))_{n \geq 1}$ is universal, and for which the maps $T(\tau^i \omega), \, i \geq 0$, do not commute.
\end{remark}

\smallskip

\bpf [Proof of Theorem \ref{propmelangeTlambdabD} ]
The operators $D$ and $T_{\lambda,b}$ do not commute, but we have the relation
$$
D  T_{\lambda, b} = \lambda \,  T_{\lambda, b}  D.
$$
Thus, for $n \geq 1,\, \omega \in \T $ and $f \in H(\C) $, we have
\begin{align}
    T_n(\omega)f = \lambda^{c(n,\omega)} f^{(a_2(n,\omega))}\left(\lambda^{a_1(n,\omega)} \, . + b \displaystyle \sum_{k=0}^{a_1(n,\omega) - 1} \lambda^k \right),
\end{align}
where $c(n, \omega)$ is a random integer such that $0 \leq c(n,\omega) \leq n$ for every $n \geq 1$ and $\omega \in \T$.
Let us denote by $r_n(\omega)$ the following number
$$r_n(\omega) := b \displaystyle \sum_{k=0}^{a_1(n,\omega) - 1} \lambda^k, $$ 
and by $E \subseteq \T$ a Borel set of measure $1$ such that 
$$\frac{a_1(n,\omega)}{n} \underset{n\to \infty}{\longrightarrow} 
 m(A_1) \quad \textrm{and} \quad \frac{a_2(n,\omega)}{n} \underset{n\to \infty}{\longrightarrow} 
 m(A_2)$$
for every $\omg \in E$. 

\smallskip
 
 Our strategy will be to apply the Universality Criterion with respect to the full sequence $(n)$, with $Z := \C[z]$ as a dense subset of $H(\C)$. As in the proof of \cite[Theorem 1]{FH}, we define a right inverse $S_n(\omega)$ on $Z$ by setting for every $k \geq 0$ and $n \geq 1$
 \begin{align*}
    S_n(\omega) z^k = \frac{k!}{(k + a_2(n, \omega))!} \, \lambda^{-c(n,\omega)} \lambda^{-k \, a_1(n, \omega)} \, \displaystyle \sum_{j=0}^{k} \binom{k+ a_2(n,\omega)}{j} z^{k+ a_2(n,\omega)-j} (-r_n(\omega))^j
 \end{align*}
 if $b \ne 0$, and by
 \begin{align*}
     S_n(\omega) z^k =  \frac{k!}{(k + a_2(n, \omega))!} \, \lambda^{-c(n,\omega)} \lambda^{-k \, a_1(n, \omega)} \, z^{k+ a_2(n, \omega)}
 \end{align*}
if $b = 0$, and extending it by linearity to $Z$.

\smallskip

We clearly have, for every $f \in Z$, $T_n(\omega) f = 0$ when $n$ is sufficiently large, since $a_2(n, \omega)  \underset{n\to \infty}{\longrightarrow}  \infty $. For every $k \geq 0, \,n \geq 1$ and $\omega \in E$, we have
\begin{align*}
    T_n(\omega) S_n(\omega) z^k &= \lambda^{-k \, a_1(n,\omega)} \ \displaystyle \sum_{j=0}^{k} \binom{k}{j} (\lambda^{a_1(n,\omega)} z + r_n(\omega))^{k-j} (-r_n(\omega))^j \\
    &= \lambda^{-k \, a_1(n,\omega)} (\lambda^{a_1(n,\omega)} z)^k \\
    &= z^k,
\end{align*}
and thus $T_n(\omega) S_n(\omega) f = f$ for every $f \in Z$. We now prove that $S_n(\omega) f \underset{n\to \infty}{\longrightarrow}  0 $ for every $\omega \in E$. To this aim, we will distinguish the cases $ \lvert \lambda \rvert < 1$ and $\lvert \lambda \rvert \geq 1$.

\smallskip

\begin{case}
Suppose first that $\lvert \lambda \rvert \geq 1$. 
\end{case}

Let $k \geq 0, \, R > 0, \, n \geq 1$ and $0 \leq j \leq k$. Let $z \in \C$ such that $\lvert z \rvert < R$. We have
\begin{align}
    \label{ineq1presqueaffreuse} &\left \lvert \frac{k!}{(k + a_2(n, \omega))!} \, \lambda^{-c(n,\omega)} \lambda^{-k \, a_1(n, \omega)} \binom{k+ a_2(n,\omega)}{j} z^{k+ a_2(n,\omega)-j} (-r_n(\omega))^j \right\rvert \\ \notag
    &\leq \frac{k! \, \lvert r_n(\omega) \rvert^j \, R^{k+a_2(n,\omega)-j}}{j! \, (k+a_2(n,\omega)-j)!} \lvert \lambda \rvert^{-k \, a_1(n,\omega)} \\ \notag
    &\leq \frac{k! \, R^{k+a_2(n,\omega)-j}}{j! \, (k+a_2(n,\omega)-j)!} \lvert b \rvert^j \, n^j \, \lvert \lambda \rvert^{n j} \lvert \lambda \rvert^{-k \, a_1(n, \omega)} \\ \notag
    &\leq C_1 \, \frac{R^{k+a_2(n,\omega)-j}}{(k+ a_2(n,\omega) -j)!} (a_2(n,\omega))^j \left (\frac{n}{a_2(n,\omega)} \right)^j \lvert \lambda \rvert^{k (n- a_1(n,\omega))} \\ \notag
    &\leq C_1 \, \frac{R^{k+a_2(n,\omega)-j}}{(k+ a_2(n,\omega) -j)!} (a_2(n,\omega))^j \left (\frac{n}{a_2(n,\omega)} \right)^j \lvert \lambda \rvert^{k \, a_2(n, \omega)},
\end{align}
where $C_1$ is a constant which depends on $k,b$ and $j$, but which is independent of $n$ and of $\omega \in E$. We used the fact that $\lvert \lambda \rvert \geq 1$ in the first inequality to bound $\lvert \lambda\rvert ^{-c(n,\omega)}$, and we used that $\lvert \lambda \rvert \geq 1 $ and that $a_1(n,\omega) \leq n$ in the second inequality to obtain that $\lvert r_n(\omega) \rvert \leq n \, \lvert b \rvert \, \lvert \lambda \rvert^n $.
Since
\begin{align*}
    \frac{R^{k+n-j}}{(k+n-j)!} \lvert \lambda \rvert^{kn} n^j \underset{n\to \infty}{\longrightarrow}  0, \quad a_2(n,\omega) \underset{n\to \infty}{\longrightarrow}  \infty \quad \textrm{and} \quad \left(\frac{n}{a_2(n,\omega)} \right)^j \underset{n\to \infty}{\longrightarrow}  m(A_2)^{-j}
\end{align*}
when $\omega \in E$, the upper bound of (\ref{ineq1presqueaffreuse}) converges to $0$ as $n \to \infty$, for every $\omega \in E$. Now summing over $0 \leq j \leq k$ the left hand side of (\ref{ineq1presqueaffreuse}), we obtain that $S_n(\omega) z^k \underset{n\to \infty}{\longrightarrow}  0$ for every $k \geq 0$, when $b \ne 0$. If $b = 0$, the convergence also holds because this case follows by taking $j = 0$ in the previous inequalities, using the convention that $a^0 = 1$ for every real number $a \geq 0$.

\smallskip
\begin{case}
    We now consider the case where $\lvert \lambda \rvert<1$.
\end{case}
Let $k \geq 0, \, R > 0, \, n \geq 1$ and $0 \leq j \leq k$. Let $z \in \C$ such that $\lvert z \rvert < R$. We have
\begin{align}
   \notag &\left \lvert \frac{k!}{(k + a_2(n, \omega))!} \, \lambda^{-c(n,\omega)} \lambda^{-k \, a_1(n, \omega)} \binom{k+ a_2(n,\omega)}{j} z^{k+ a_2(n,\omega)-j} (-r_n(\omega))^j \right\rvert \\ \notag
    &\leq C_2 \; \lvert \lambda \rvert^{-(k+1)n} \, \frac{R^{k+a_2(n,\omega)-j}}{(k+a_2(n,\omega)-j)!} \\ \notag
    &\leq C_3 \; \lvert \lambda \rvert^{-(k+1)n} \, \frac{R^{a_2(n,\omega)}}{a_2(n,\omega)!},
\end{align}
where $C_2$ and $C_3$ are two constants depending on $k, b$ and $j$, but that are independent of $\omega \in E$ and $n \geq 1$.
For the first inequality, we used the fact that $\lvert \lambda \rvert < 1,\, c(n,\omega) \leq n$ and $a_1(n,\omega) \leq n$, and we also used the fact that $\lvert \lambda \rvert < 1$ to assert that
$$\lvert r_n(\omega) \rvert^j \leq \lvert b \rvert^j \left( \displaystyle \sum_{l \geq 0} \lvert \lambda \rvert^l \right)^j \leq (1-\lvert \lambda \rvert)^{-j} \, \lvert b \rvert^j.$$

By Stirling's formula, there is a constant $C > 0$ such that
\begin{align} \label{eqstirlingcocycl}
    \frac{1}{n!} \leq C \, \left (\frac{e}{n} \right)^n \, n^{-1/2}
\end{align}
for every $n \geq 1$. Using (\ref{eqstirlingcocycl}), we deduce that there exists a constant $C_4$ which is independent of $\omg \in E$ and $n \geq 1$, such that for every $n \geq 1,\, k \geq 0, \, 0 \leq j \leq k$, and $\lvert z \rvert < R$,
\begin{align}
   \label{ineq1affreuse} &\left \lvert \frac{k!}{(k + a_2(n, \omega))!} \, \lambda^{-c(n,\omega)} \lambda^{-k \, a_1(n, \omega)} \binom{k+ a_2(n,\omega)}{j} z^{k+ a_2(n,\omega)-j} (-r_n(\omega))^j \right\rvert \\ \notag
    &\leq C_4 \; \lvert \lambda \rvert^{-(k+1)n} \, R^{a_2(n,\omega)} \left (\frac{e}{a_2(n,\omega)}\right)^{a_2(n,\omega)} (a_2(n,\omega))^{-1/2} \\ \notag
    &\leq C_4 \; \exp{(u_n(\omega))},
\end{align}
where 
\begin{align*}
    u_n(\omega) :=\; &a_2(n,\omega)(1- \log(a_2(n,\omega))) - \tfrac{1}{2} \log(a_2(n,\omega)) \\ &- (k+1) \, n \, \log\lvert \lambda \rvert + a_2(n,\omega) \log(R)
\end{align*}
for every $n \geq 1$ and $\omega \in E$. Finally, we have
\[
\begin{array}{rcl}
    u_n(\omega) &\underset{n \to +\infty}{=}& (n \, m(A_2) + o(n))(1- \log(n \, m(A_2) + o(n)))- \tfrac{1}{2} \log(n \, m(A_2) + o(n)) \\ && - n(k+1) \log(\lvert \lambda \rvert) + n \, m(A_2) \log(R) + o(n) \\ \\
    &\underset{n \to +\infty}{=}& n \,(m(A_2) + m(A_2) \log(R) - (k+1) \log(\lvert \lambda \rvert)) + o(n) \\ && - \log(n \, m(A_2) + o(n))(\tfrac{1}{2} + n \, m(A_2) + o(n)).
\end{array}
\]
From this, one can check that
$$
u_n(\omega) \underset{n \to \infty}{\sim} -n \, m(A_2) \, \log(n \, m(A_2)),
$$
which in particular proves that $\exp(u_n(\omg))$ converges to $0$ as $n \to \infty$ when $\omg \in E$.
It follows that the left hand side of (\ref{ineq1affreuse}) converges to $0$ as $n \to \infty$. Summing over $0 \leq j \leq k$ the left hand side of (\ref{ineq1affreuse}), we obtain that $S_n(\omega) z^k \underset{n\to \infty}{\longrightarrow}  0$ for every $k \geq 0$, when $b \ne 0$. If $b = 0$, the convergence also holds because this case follows by taking $j = 0$ in the previous inequalities, using again the convention that $a^0 = 1$ for every real number $a \geq 0$. 

This proves that the sequence $(T_n(\omega))_{n \geq 1}$ satisfies Proposition \ref{critunivforsequence} for every $\omega \in E$. In particular, the sequence $(T_n(\omega))_{n \geq 1}$ is topologically mixing for every $\omega
\in E$.
\epf

\begin{remark}
Again, the conclusions of Theorems \ref{propositionfoncenttypeexpo} and \ref{propmelangeTlambdabD} hold for every $\omega \in \T$ when the ergodic transformation $\tau$ is an irrational rotation and when $A_1, A_2$ are intervals of $[0,1)$, that is, the sequence $(T_n(\omg))_{n \geq 1}$ is topologically mixing for every $\omg \in \T$.
\end{remark}

\section{Comments and open problems}\label{Questions}
We end this article with some comments and open questions related to our previous results. 

\smallskip

But before that, it is worth noting that many of the results stated in this paper remain valid for any measure-preserving ergodic dynamical system. This is the case, for example, of Propositions \ref{tcllimsupliminf}, \ref{rmknonunivmultiplication}, \ref{propcondisuffisnormt_n(omega)notunivmesegal}, \ref{propphi1phi2constante}, \ref{propnonunivnormcasA1A2pasmmmesure}, \ref{lemnonunivpourtransfoerg}, as well as Theorems \ref{propimagephi1phi2cercle}, \ref{propunivphi1phi2intetimagephicasegal}, \ref{propunivphi1phi2innerA1A2mmesuresanshypimages}, \ref{propimagephi1phi2cerclecasmesnonegales}, \ref{propunivhypimagesphi1phi2caspasmmesure}, \ref{propunivhypphi1phi2nonoutercaspasmmesure}, \ref{propositionfoncenttypeexpo} and \ref{propmelangeTlambdabD}.

\smallskip

\subsection{A general remark on the universality of the random sequence $(T_n(.))_{n \geq 1}$}
Let us observe that in all of our results, the sequence of operators $(T_n(\omg))_{n \geq 1}$ is either universal for almost every $\omg \in \T$, or not universal for almost every $\omg \in \T$. Most of the time, this comes from Birkhoff's Theorem \ref{birkhoffthintro}, or this comes from the fact that the sequence of Birkhoff sums satisfies a CLT. We in fact have the following $0-1$ law.

\smallskip

\bpr \label{0-1lawrandomsequence}
Let $\tau : \Omg \to \Omg$ be an ergodic measure-preserving transformation on a Polish probability space $(\Omg, \textrm{Bor}(\Omg), m)$, where $\textrm{Bor}(\Omg)$ is the Borel $\sigma$-algebra of $\Omg$ and $m$ is a Borel probability measure on $\Omg$. Let $X$ be a separable Banach space. Let $T : \Omg \to \bx$ be strongly measurable, which means that for every $x \in X$, the map $\omg \mapsto T(\omg) x$ is $(\textrm{Bor}(\Omg), \textrm{Bor}(X))$-measurable for the Borel $\sigma$-algebras of $\Omg$ and $X$. Then, either the sequence $(T_n(\omg))_{n \geq 1}$ is universal for almost every $\omg \in \Omg$, or the sequence $(T_n(\omg))_{n \geq 1}$ is not universal for almost every $\omg \in \Omg$.  
\epr

\smallskip

\bpf
Let us denote by $A$ the set 
$$
A := \{ \omg \in \Omg : (T_n(\omg))_{n \geq 1} \; \textrm{is universal} \}.
$$
Let $(V_j)_{j \geq 1}$ be a countable basis of open sets of $X$.
Let us remark that the set $A$ can be written as
$$
A = \textrm{proj}_\Omega  (\displaystyle \bigcap_{j \geq 1} \bigcup_{n \geq 1} \{ (\omg, x) \in \Omg \times X : T_n(\omg) x \in V_j \} ),
$$
where $\textrm{proj}_\Omega$ is the projection onto $\Omega$. Now since the Banach space $X$ is separable, the strong measurability of the map $T : \Omg \to \bx$ is equivalent to the $(\textrm{Bor}(\Omg) \otimes \textrm{Bor}(X), \textrm{Bor}(X))$-measurability of the map $(\omg, x) \mapsto T(\omg) x$. Moreover, the composition of strongly measurable maps is strongly measurable (see \cite[Appendix A]{GTQ} for these remarks). In particular, for each $n \geq 1$, the map $(\omg, x) \mapsto T_n(\omg) x$ is measurable. Thus, the set $A$ is the projection of a Borel subset of $\Omg \times X$ and is an analytic subset of $\Omg$ (see \cite{Kechris} for background on analytic subsets). In particular, the set $A$ is universally measurable (\cite[Theorem 21.10, page 155]{Kechris}) and belongs to the completion of the Borel $\sigma$-algebra of $\Omg$, that we denote by $\Sigma^*$. We also denote by $m^*$ the completion of the measure $m$ on $\Omg$. Every element of $\Sigma^*$ can be written as $B \cup N$, with $B$ a Borel subset of $\Omg$ and $N$ a subset of $\Omg$ such that $N \subset C$ with $C$ a Borel subset satisfying $m(C) = 0$. Moreover, $m^*(B \cup N) = m(B)$. We can easily show that the transformation $\tau$ is measure-preserving for $m^*$. Let $D = B \cup N$ be in $\Sigma^*$ and satisfy $\tau^{-1}(D) = D$, with $N \subset C$ and $m(C) = 0$. From $B \cup N = \tau^{-1}(B) \cup \tau^{-1}(N)$, it follows that $B \subset \tau^{-1}(C) \cup \tau^{-1}(B)$ and $\tau^{-1}(B) \subset C \cup B$, and thus $m(\tau^{-1}(B) \triangle B) = 0$. Since $\tau$ is ergodic for $m$, it follows that $m(B) = m^*(D) = 0$ or $m(B) = m^*(D) = 1$, which proves that $\tau$ is ergodic for $m^*$.

Now let us notice that $\tau(A) \subset A$, since $T_{n+1}(\omg) = T_n(\tau \omg) T(\omg)$ for every $n \geq 1$ and every $\omg \in \Omg$. By ergodicity, we thus have $m^*(A) = 0$ or $m^*(A) = 1$, which in particular implies the conclusion of Proposition \ref{0-1lawrandomsequence}.
\epf

\subsection{Comments on random products of adjoints of multiplication operators}
Let us consider the sequence $(T_n(\omega))_{n \geq 1}$ given by 
\begin{equation} \label{eqcocycleadjointmultipliccomments}
    T(\omega) = \left\{
    \begin{array}{ll}
        (M_{\phi_1})^* & \mbox{if} \quad \omega \in A_1 \\
        (M_{\phi_2})^* & \mbox{if} \quad \omega \in A_2
    \end{array},
\right.
\end{equation}
where $\phi_1, \phi_2 \in H^\infty(\D)$ are nonconstant functions, $A_1, A_2$ are two disjoint Borel subsets of $[0,1)$ such that $A_1 \cup A_2 = [0,1)$ and $m(A_1) = m(A_2) =1/2 $, and with $\tau$ an ergodic measure-preserving transformation on $(\T,m)$. For example, let us suppose that $A_1 = [0,1/2)$ and $A_2 = [1/2,1)$.

\smallskip

We know that, thanks to Theorem \ref{propimagephi1phi2cercle}, if $\phi_1 \phi_2$ is nonconstant and such that $(\phi_1 \phi_2)(\D) \cap \T \ne \emptyset$, then the random sequence $(T_n(.))_{n \geq 1}$ is universal. This situation thus leads to two cases that we should call limit cases: the case $(\phi_1 \phi_2)(\D) \subset \D$ and the case $(\phi_1 \phi_2)(\D) \subset \C \setminus \overline{\D}$. We have also studied some cases that lead to one of these two limit cases (Proposition \ref{situationphi1phi2DcontenuedansD} and Corollary \ref{corcasphi1phi2exterieurdisk}). These two sub-cases can be seen as trivial situations of the two limit cases. We have also studied the situation where $\phi_1 \phi_2$ is inner, which is included in the case $(\phi_1 \phi_2)(\D) \subset \D$ (Theorems \ref{propunivphi1phi2intetimagephicasegal} and \ref{propunivphi1phi2innerA1A2mmesuresanshypimages}). We can sum up our work on the two limit cases with the following table.

\medskip

\renewcommand{\arraystretch}{1.2}
\setlength{\tabcolsep}{6pt}

\begin{center}
\begin{tabular}{|m{7cm}|m{7cm}|}
\hline
\begin{minipage}[c][3cm][c]{\linewidth}
\centering
$\phi_1 \phi_2$ inner.\\
$\phi_1 (\D) \cap \T \ne \emptyset$ and $1/\phi_1 \in H^\infty(\D)$, or $\phi_2 (\D) \cap \T \ne \emptyset$ and $1/\phi_2 \in H^\infty(\D)$.\\
$\displaystyle \limsup_{n \to \infty} \mathbb{S}_n^{\tau}f = \infty$ or  $\displaystyle \liminf_{n \to \infty} \mathbb{S}_n^{\tau}f = -\infty$.\\
\textcolor{red}{Universality of the random sequence}
\end{minipage}
&
\begin{minipage}[c][3cm][c]{\linewidth}
\centering
$\phi_1 \phi_2$ inner and nonconstant.\\
$\phi_1$ or $\phi_2$ not inner. \\
$\displaystyle \limsup_{n \to \infty} \mathbb{S}_n^{\tau}f = \infty$ and  $\displaystyle \liminf_{n \to \infty} \mathbb{S}_n^{\tau}f = -\infty$.\\
\textcolor{red}{Universality of the random sequence}
\end{minipage}
\\ \hline
\begin{minipage}[c][3cm][c]{\linewidth}
\centering
$\phi_1(\D) \subset \D$ and $\phi_2(\D) \subset \D$.\\
\textcolor{red}{Non-universality of the random sequence}
\end{minipage}
&
\begin{minipage}[c][3cm][c]{\linewidth}
\centering
$\phi_1(\D) \subset \C \setminus \overline{\D}$ and $\phi_2(\D) \subset \C \setminus \overline{\D}$.\\
\textcolor{red}{Non-universality of the random sequence}
\end{minipage}
\\ \hline
\end{tabular}
\end{center}

\medskip

Certain situations are missing in our work. For example, we don't know if the sequence $(T_n(.))_{n \geq 1}$ is universal when $\phi_1, \phi_2 \in H^\infty(\D)$ are nonconstant such that $\phi_1(\D) \subset \D$, $\phi_2(\D) \subset \C \setminus \overline{\D}$, $(\phi_1 \phi_2)(\D) \subset \D$, $\phi_1 \phi_2$ is nonconstant and not inner, and $\lVert \phi_1 \rVert_\infty \lVert \phi_2 \rVert_\infty \geq 1$.
However, for certain types of functions $\phi_1$ and $\phi_2$ in this last case, we can observe that the sequence $(T_n(.))_{n \geq 1}$ is not universal. This is the case of the following example.

\smallskip

\begin{example} \label{examplecommentscaslimfacil}
    Let $\phi_1(z) = e^{1-z}$ and $\phi_2(z) = e^{1/2(z-1)}$ for every $z \in \D$. In this case, we have $\phi_1(\D) \subset \C \setminus \overline{\D}$, $\phi_2(\D) \subset \D$, $\lVert \phi_1 \rVert_{\infty} = e^2$, $\lVert \phi_2 \rVert_{\infty} = 1$, and $\lVert \phi_1 \rVert_{\infty} \lVert \phi_2 \rVert_{\infty}= e^2$. Also, $(\phi_1 \phi_2)(\D) \subseteq \C \setminus \overline{\D}$. In this case, the operator $T_n(\omega)$ is the adjoint of the multiplication operator by the function $e^{(1-z)(a_1(n,\omega) - \frac{1}{2}a_2(n,\omega))}$. Moreover, the map $T_n(\omega)$ is invertible and its inverse $S_n(\omega)$ is the adjoint of the multiplication operator by the function $e^{(z-1)(a_1(n,\omega) - \frac{1}{2}a_2(n,\omega))}$. Thus,
    $$
    \lVert S_n(\omega) \rVert = \displaystyle \sup_{z \in \D} \{ e^{(\Re(z)-1)(a_1(n,\omega) - \frac{1}{2}a_2(n,\omega))} \} = 1
    $$
    when $n\geq 1$ is sufficiently large, for almost every $\omega \in \T$, because $\frac{1}{n}(a_1(n,\omega) - \frac{1}{2} a_2(n,\omega))$ converges to $1/4$ for almost every $\omega \in \T$. This shows that both sequences $(T_n(\omega))_{n \geq 1}$ and $(S_n(\omega))_{n \geq 1}$ are non-universal, for almost every $\omega \in \T$. 
\end{example}

\smallskip

Example \ref{examplecommentscaslimfacil} gives a nontrivial situation where the random sequence $(T_n(.))_{n \geq 1}$ is non-universal, and where the limit cases hold. But this example is quite simple because we could use the norm of the inverse of  $\lVert T_n(\omega) \rVert$ to prove the non-universality. Unfortunately, it is not always possible to do it in general, since it depends heavily on the expression of the functions $\phi_1$ and $\phi_2$. Also, let us observe that the functions $\phi_1$ and $\phi_2$ in Example \ref{examplecommentscaslimfacil} are both outer. More generally, in the limit case $(\phi_1 \phi_2)(\D) \subset \C \setminus \overline{\D}$ and $\phi_1 \phi_2$ nonconstant, both functions $\phi_1$ and $\phi_2$ are outer since $1/\phi_1$ and $1/\phi_2$ belong to $H^\infty(\D)$. In this case, the inner parts $I_1$ and $I_2$ of $\phi_1$ and $\phi_2$ respectively are constant, and we thus cannot use the subspaces $K := \displaystyle \bigcup_{n \geq 1} (I_1^n H^2(\D))^\perp$ and $K' := \displaystyle \bigcup_{n \geq 1} (I_2^n H^2(\D))^\perp $ as in the proof of Theorem \ref{propunivphi1phi2innerA1A2mmesuresanshypimages} to apply the Universality Criterion. This also is what makes the limit case $(\phi_1 \phi_2)(\D) \subset \C \setminus \overline{\D}$ hard to solve.

\smallskip

Regarding the limit case $(\phi_1 \phi_2)(\D) \subset \D$ and $\phi_1 \phi_2$ nonconstant, we believe that a good understanding of the following example will help us to solve the limit case $(\phi_1 \phi_2)(\D) \subset \D$. 

\smallskip

\begin{example}
Let $\phi_1, \phi_2$ be the two functions defined by $\phi_1(z) = e^{-1} z$ and $\phi_2(z) = e^{z}$ for every $z \in \D$. 
In this case, $(\phi_1 \phi_2)(z) = z  e^{z-1}$, and thus $(\phi_1 \phi_2)(\D) \subset \D$. Moreover, we have $\phi_2(\D) \cap \T \ne \emptyset$, $\phi_1(\D) \subset \D$, $\lVert \phi_1 \rVert_{\infty} = e^{-1}$, $\lVert \phi_2 \rVert_{\infty} = e$, and $\lVert \phi_1 \rVert_{\infty} \lVert \phi_2 \rVert_{\infty}= 1$.The operator $T_n(\omg)$ is the adjoint of the multiplication operator by the function $e^{-a_1(n,\omg)} z^{a_1(n,\omg)} e^{a_2(n,\omg) z}$. The operator $T_n(\omega)$ is not invertible, but we can use the subspace $\C[z]$ of the complex polynomials as a dense subspace for the Universality Criterion, since $T_n(\omg) f$ converges to $0$ for every complex polynomial $f$, for almost every $\omg \in \T$. However, defining a right inverse $S_n(\omg)$ on the eigenvectors $k_z$ as in the proof of Theorem \ref{propunivphi1phi2innerA1A2mmesuresanshypimages} is not enough to apply the Universality Criterion here. Indeed, let us define a right inverse $S_n(\omg)$ on the eigenvectors $k_z$ by
$$
S_n(\omg) k_z = e^{a_1(n,\omg)} \, e^{-a_2(n,\omg) \, \overline{z}} M_z^{a_1(n,\omg)} k_z
$$
and extending it by linearity to the dense subspace $Z := \textrm{span}\{ k_z : z \in \D \}$ of $H^2(\D)$. We can observe that $$\lVert S_n(\omg) k_z \rVert = e^{a_1(n,\omg) - a_2(n,\omg) \Re{z}} \lVert k_z \rVert,$$
which diverges to $\infty$ for every $z \in \D$, for almost every $\omg \in \T$. This makes this case not that easy to solve.

Finally, let us observe that
$$
\lVert T_n(\omega) \rVert = \displaystyle \sup_{z \in \D} \{ e^{a_2(n,\omega) \Re(z)} e^{-a_1(n,\omega)} \lvert z\rvert ^{a_1(n,\omega)} \} = e^{a_2(n,\omega) - a_1(n,\omega)}.
$$
If the sequence $(\mathbb{S}_n^{\tau}f(\omg))_{n \geq 1}$ of Birkhoff sums $\mathbb{S}_n^{\tau}f(\omg) = a_2(n,\omg) - a_1(n,\omg)$ is bounded for almost every $\omg \in \T$, then the random sequence $(T_n(.))_{n \geq 1}$ associated to this transformation $\tau$ is not universal. But in the case of an irrational rotation or the doubling map, it is not clear if this remains true. 
\end{example}

\smallskip

The following problem is thus still open.

\smallskip

\begin{question}
    Let $\tau$ be an irrational rotation or the doubling map on $\T$. Let $A_1 = [0,1/2)$ and $A_2 = [1/2, 1)$. Let $(T_n(.))_{n \geq 1}$ be the random sequence of operators defined by the map 
    \begin{equation} \notag
    T(\omega) = \left\{
    \begin{array}{ll}
        (M_{\phi_1})^* & \mbox{if} \quad \omega \in A_1 \\
        (M_{\phi_2})^* & \mbox{if} \quad \omega \in A_2
    \end{array},
\right.
\end{equation}
where $\phi_1, \phi_2 \in H^\infty(\D)$ are nonconstant functions. Suppose that we are in one of the two limit cases: $\phi_1 \phi_2$ is nonconstant and either $(\phi_1 \phi_2)(\D) \subset \D$ or $(\phi_1 \phi_2)(\D) \subset \C \setminus \overline{\D}$, and that no trivial situation on $\phi_1$ and $\phi_2$ holds. Suppose also that $\phi_1 \phi_2$ is not inner. Is it true that the random sequence $(T_n(.))_{n \geq 1}$ is not universal?
\end{question}

\medskip

When $A_1, A_2$ are two disjoint Borel subsets of $[0,1)$ such that $A_1 \cup A_2 = [0,1)$, Theorems \ref{propunivphi1phi2innerA1A2mmesuresanshypimages}, \ref{propunivhypimagesphi1phi2caspasmmesure} and \ref{propunivhypphi1phi2nonoutercaspasmmesure} require the Birkhoff sums $(\mathbb{S}_n^{\tau}f(\omg))_{n \geq 1}$ associated to the function $f = \indic_{A_1} - \frac{m(A_1)}{m(A_2)} \indic_{A_2}$  to satisfy
$$
\displaystyle \limsup_{n \to \infty} \mathbb{S}_n^{\tau}f(\omg) = \infty \quad \textrm{and} \quad \displaystyle \liminf_{n \to \infty} \mathbb{S}_n^{\tau}f(\omg) = -\infty
$$
for almost every $\omega \in \T$, in order to make the second part of the Universality Criterion on the right inverses work. This raises the following question.

\smallskip

\begin{question}
    Do Theorems \ref{propunivphi1phi2innerA1A2mmesuresanshypimages}, \ref{propunivhypimagesphi1phi2caspasmmesure} and \ref{propunivhypphi1phi2nonoutercaspasmmesure} work if we replace the condition on the Birkhoff sums by the condition
    $$
\displaystyle \limsup_{n \to \infty} \mathbb{S}_n^{\tau}f(\omg) = \infty \quad \textrm{or} \quad \displaystyle \liminf_{n \to \infty} \mathbb{S}_n^{\tau}f(\omg) = -\infty
    $$
for almost every $\omega \in \T$?
\end{question}

\smallskip

Corollaries \ref{cordoublmapcaspasmmesure}, \ref{coralphabpqrotcaspasmmesure} and \ref{coralphanonbpqrotcaspasmmesure} hold only for certain intervals of $[0,1)$, since their proofs make use of the computation of the Fourier coefficients of the function $f = \indic_{A_1} -\tfrac{m(A_1)}{m(A_2)} \indic_{A_2}$ in order to apply CLTs for the doubling map and for the irrational rotations. Thus, the following question is natural.

\smallskip

\begin{question}
    Do Corollaries  \ref{cordoublmapcaspasmmesure}, \ref{coralphabpqrotcaspasmmesure} and \ref{coralphanonbpqrotcaspasmmesure} still hold if we replace the two intervals by disjoint Borel subsets of $[0,1)$ such that $A_1 \cup A_2 = [0,1)$ and $m(A_k) > 0$ for $k = 1,2$?
\end{question}

\smallskip

In the same vein, we have seen in Corollary \ref{coralphabZalpharot} that for some irrational rotations, the random sequence $(T_n(.))_{n \geq 1}$ won't be universal. Thus, the following question is also natural.

\smallskip

\begin{question}
    Let $\alpha$ be an irrational in $(0,1)$ and let $\tau = R_\alpha$. Let $A_1 = [0,b)$ and $A_2 = [b,1)$ with $b \in (0,1)$. Let $(T_n(.))_{n \geq 1}$ be the random sequence of operators defined by the map 
    \begin{equation} \notag
    T(\omega) = \left\{
    \begin{array}{ll}
        (M_{\phi_1})^* & \mbox{if} \quad \omega \in A_1 \\
        (M_{\phi_2})^* & \mbox{if} \quad \omega \in A_2
    \end{array},
\right.
\end{equation}
where $\phi_1, \phi_2 \in H^\infty(\D)$ are nonconstant functions. Do Corollaries \ref{coralphabpqrotcaspasmmesure} and \ref{coralphanonbpqrotcaspasmmesure} hold for these sets $A_1$ and $A_2$, if $b \notin \Z \alpha$?
\end{question}

\smallskip

In our work, we only use two disjoint Borel subsets $A_1$ and $A_2$ of $[0,1)$ such that $A_1 \cup A_2 = [0,1) $ and $m(A_k) > 0$ for $k = 1,2$. It is therefore natural to be interested in the case of more than two Borel subsets of $[0,1)$.

\smallskip

\begin{problem} \label{openproblemfinitenumberborels}
    Generalize our work on the random products of adjoints of multiplication operators for a finite number of disjoint Borel subsets $A_k$ of $[0,1)$ covering $[0,1)$ such that $m(A_k) > 0$ for every $k \geq 1$, when the operator $T(\omg)$ is the adjoint of the multiplication operator $(M_{\phi_k})^*$ for every $\omg \in A_k$, with $\phi_k \in H^\infty(\D)$ nonconstant.
\end{problem}

\smallskip

Obviously, we can imagine similar results to Proposition \ref{propnonunivnormcasA1A2pasmmmesure} and Theorem \ref{propimagephi1phi2cerclecasmesnonegales} in the context of Problem \ref{openproblemfinitenumberborels}. Finally, one can also look at the linear dynamics of random products $(T_n(.))_{n \geq 1}$ when the operator $T(\omg)$ is another type of operator on a space of holomorphic functions.

\smallskip

\begin{problem}  \label{openproblemproductsofotheroperators}
    Study the linear dynamics of random sequences $(T_n(.))_{n \geq 1}$ when the operator $T(\omg)$ is given by
    \begin{equation} \notag
    T(\omega) = \left\{
    \begin{array}{ll}
        T_1 & \mbox{if} \quad \omega \in A_1 \\
        T_2 & \mbox{if} \quad \omega \in A_2
    \end{array},
\right.
\end{equation}
where $T_1, T_2 $ are two operators belonging to some specific families of operators on separable Fréchet spaces, where $\tau$ is either the doubling map or an irrational rotation on $\T$, and where $A_1, A_2$ are two disjoint Borel subsets of $[0,1)$ such that $A_1 \cup A_2 = [0,1)$ and $m(A_k) > 0$ for $k = 1,2$.
\end{problem}

\smallskip

For example, it would be natural to consider composition operators for $T_1$ and $T_2$ on a space of holomorphic functions and to study the linear dynamics of the corresponding sequence $(T_n(.))_{n \geq 1}$ in the situation of Problem \ref{openproblemproductsofotheroperators}, as well as weighted shifts on $\ell_p$-spaces, since the linear dynamics of these operators is rather well known and of interest for many problems in linear dynamics.

\subsection{The case of a non-ergodic transformation}

Almost all of our previous results in this work hold for an ergodic measure-preserving transformation on $(\T,m)$. Without the assumption on the ergodicity of the transformation, Birkhoff's theorem (Theorem \ref{birkhoffthintro}) still holds, but the limit is not necessarily constant. Thus, our previous main results won't hold in this context, except Proposition \ref{situationphi1phi2DcontenuedansD}, Corollary \ref{corcasphi1phi2exterieurdisk} and Propositions \ref{propphi1phi2constante} and \ref{lemnonunivpourtransfoerg}. In order to investigate the problem of universality of random sequences of products of operators $(T_n(.))_{n \geq 1}$ for non-ergodic transformations on $(\T,m)$, it is natural to start the investigation with the case of a rational rotation.

\smallskip

Let $\alpha = \frac{p}{q}$ be a rational number in $(0,1)$, with $p$ and $q$ coprime. Let $\tau$ be the rotation with parameter $\alpha$, which is non-ergodic. In this case, for $f$ a real centered bounded function on $\T$, the coboundary equation $f = h - h \circ \tau$ has a solution $h \in L^\infty(\T)$ if and only if 
\begin{align} \label{equationrotrationellefin}
    \displaystyle \sum_{j=0}^{q-1} f(x + j \alpha) = 0
\end{align}
for almost every $x \in [0,1)$. In this case, the function 
$$
h(x) := -\frac{1}{q} \displaystyle \sum_{j=0}^{q-1} (j+1) f(x + j \alpha) 
$$
provides a bounded solution of the coboundary equation. Thus, when Equation \ref{equationrotrationellefin} has a solution for $f = \indic_{A_1} - \frac{m(A_1)}{m(A_2)} \indic_{A_2}$, the Birkhoff sums associated to the function $f$ are bounded. In particular, Proposition \ref{lemnonunivpourtransfoerg} applies in this case.

\end{document}